\documentclass{article}

\usepackage{amsmath}
\usepackage{amsthm}
\usepackage[utf8]{inputenc}
\usepackage{amsfonts}
\usepackage{hyperref}
\hypersetup{
	colorlinks = true,
	linkcolor = purple, 
	filecolor = purple, 
	citecolor = purple,
	urlcolor = purple
}
\usepackage{amssymb}
\usepackage{graphicx}
\usepackage{stackrel}
\usepackage[english]{babel}
\usepackage{url}
\usepackage{xcolor}
\usepackage{tikz-cd}
\usepackage{pdfpages}
\usepackage{csquotes}
\allowdisplaybreaks

\newtheorem{theorem}{Theorem}[section]
\newtheorem{introtheorem}{Theorem}
\newtheorem{introprop}{Proposition}
\newtheorem*{theorem*}{Theorem}
\newtheorem{prop}[theorem]{Proposition}
\newtheorem*{prop*}{Proposition}
\newtheorem{lem}[theorem]{Lemma}
\newtheorem{cor}[theorem]{Corollary}

\newtheorem*{conj*}{Conjecture}
\theoremstyle{definition}
\newtheorem{defn}[theorem]{Definition}
\newtheorem{nota}[theorem]{Notation}
\newtheorem{constr}[theorem]{Construction}
\newtheorem{exa}[theorem]{Example}
\newtheorem{rem}[theorem]{Remark}

\begin{document}
	\title{The quadratic Euler characteristic of a smooth projective same-degree complete intersection}
	\author{Anna M. Viergever}
	\maketitle 
	
	\begin{abstract}
		We find an algorithm to compute the quadratic Euler characteristic of a smooth projective complete intersection of hypersurfaces of the same degree, generalizing the argument of \cite{LevineECHWHH}. As an example, we compute the quadratic Euler characteristic of a smooth projective complete intersection of two generalized Fermat hypersurfaces. 
	\end{abstract}
	
	\phantomsection
	\addcontentsline{toc}{section}{Introduction}

	\section*{Introduction}
	The quadratic Euler characteristic of a smooth projective scheme over a perfect field of characteristic not equal to $2$ is a refined or motivic analogue of the usual topological Euler characteristic, and of Euler characteristics defined using \'etale or De Rham cohomology. It comes from a very general definition of a categorical Euler characteristic associated to a dualizable object of a symmetric monoidal category, living in the endomorphism ring of the unit. Motivic homotopy theory, introduced by Morel and Voevodsky, constructs the stable motivic homotopy category $\text{SH}(k)$ of a field $k$; a symmetric monoidal category in which a smooth projective scheme over $k$ has an image which is dualizable. A deep theorem by Morel (see \cite[Theorem 6.4.1]{MorelIAHT}) states that if $k$ is perfect and not of characteristic $2$, then the endomorphism ring of the unit in $\text{SH}(k)$ is isomorphic to the Grothendieck-Witt ring $\text{GW}(k)$ of~$k$. This is the group completion of the ring of all isometry classes of nondegenerate quadratic forms over~$k$. Therefore, we obtain the \textit{quadratic Euler characteristic} $\chi(X/k)\in \text{GW}(k)$ of any smooth projective scheme $X$ over $k$, which is a (virtual) quadratic form. \\
	These quadratic Euler characteristics carry a lot of information within them: if $k\subset\mathbb{R}$ then the rank of $\chi(X/k)$ is equal to the topological Euler characteristic of the $\mathbb{C}$-points $X(\mathbb{C})$, while the signature of $\chi(X/k)$ is the topological Euler characteristic of the real points $X(\mathbb{R})$. Quadratic Euler characteristics are often used in the fast-growing field of refined enumerative geometry, which aims to obtain ``quadratic enrichments" of results in classical enumerative geometry. However, they are in general hard to compute. \\
	The motivic Gauss-Bonnet Theorem (see \cite{LevineGB}) proven by Levine and Raksit gives a rather explicit way to compute quadratic Euler characteristics. Namely, consider a smooth projective scheme $X$ over a perfect field $k$ which is not of characteristic $2$ as before. For $a\in k^*$, let $\langle a \rangle$ be the quadratic form $x\mapsto ax^2\in\text{GW}(k)$. Then we can compute $\chi(X/k)\in\text{GW}(k)$ as follows: 
	\begin{itemize}
		\item If $\dim(X)$ is odd, then $\chi(X/k) = C\cdot H$ for some $C\in \mathbb{Z}$, where $H$ is the hyperbolic form $\langle 1 \rangle + \langle -1 \rangle$. 
		\item If $\dim(X) = 2n$ is even, then $\chi(X/k) = C\cdot H + Q$ for some $C\in \mathbb{Z}$, where $Q$ is the quadratic form given by the composition 
		$$H^n(X,\Omega_X^n)\times H^n(X,\Omega_X^n)\xrightarrow{\cup}H^{2n}(X,\Omega_X^{2n})\xrightarrow{\text{Trace}}k.$$
		Here, $\Omega_X$ denotes the sheaf of differential forms on $X$, the first map is the cup product on cohomology and we write $\Omega_X^q = \wedge^q\Omega_X$. 
	\end{itemize}
	The constant $C$ can be computed in practice. Therefore, one can compute the quadratic Euler characteristic of a smooth projective scheme if one understands the form $Q$ in the even dimensional case. This form has been computed successfully in the case of hypersurfaces by Levine, Lehalleur and Srinivas in \cite{LevineECHWHH}. Given a smooth projective hypersurface $X = V(F)\subset\mathbb{P}^n$, the authors use inspiration from the paper \cite{CarlsonIVHS} by Carlson and Griffiths to describe an isomorphism from the primitive cohomology of $\Omega_X^q$ to certain graded pieces of the Jacobian ring $$J_X = k[X_0,\cdots, X_n]/\left(\frac{\partial F}{\partial X_0},\cdots, \frac{\partial F}{\partial X_n}\right).$$ This is then applied to compare the cup product on cohomology with the usual ring multiplication of $J_X$. In doing this, the authors represent the result on an open cover of $X$ and compare this with a representation of $c_1(\mathcal{O}(m))^n$ to compute the trace, which then allows them to compute the form $Q$ explicitly. \\
	The purpose of this paper is to show a similar result for smooth projective complete intersections of hypersurfaces which are of the same degree. For this, some work in the style of \cite{CarlsonIVHS} has already been done by Konno in the paper \cite{KonnoVTPCI} and Terasoma in the paper \cite{TerasomaIVH}. Using inspiration from those papers, given a smooth projective complete intersection $X=V(F_0,\cdots, F_r)\subset \mathbb{P}^n$ where the $F_i$ are of the same degree $m>1$ and $n\geq r+2$, we consider the hypersurface $$\mathcal{X} = V(F)\subset\mathbb{P}^r\times\mathbb{P}^n$$ where $F = Y_0F_0 + \cdots + Y_rF_r$ and show that one can compute $\chi(X/k)$ from~$\chi(\mathcal{X}/k)$. After this, we study isomorphisms from the primitive cohomology groups of $\Omega_{\mathcal{X}}^q$ to graded parts of the Jacobian ring 
	$$J = k[Y_0,\cdots, Y_r, X_0,\cdots, X_n]/\left(F_0,\cdots, F_r, \frac{\partial F}{\partial X_0},\cdots, \frac{\partial F}{\partial X_n}\right).$$
	Unlike the Jacobian ring for the hypersurface, this Jacobian ring is infinite dimensional over $k$. One can link the cup product on cohomology to the usual product in the Jacobian ring, namely, following \cite{TerasomaIVH}, we show the following result. 
	\begin{introprop}[See Corollary \ref{corollary definition map Jrho to cohomology}]
		Consider the bidegree $$\rho = (n-r-1, (n+r+1)m - 2(n+1)).$$ There is a surjective homomorphism $\phi: J^\rho\to H^{n+r}(\mathbb{P}^r\times\mathbb{P}^n, \Omega_{\mathbb{P}^r\times\mathbb{P}^n}^{n+r})\cong k$, such that the diagram 
		\[
		\begin{tikzcd} 
		H^q(\mathcal{X}, \Omega_{\mathcal{X}}^p)_{prim}\otimes H^p(\mathcal{X}, \Omega_{\mathcal{X}}^q)_{prim} \arrow[r, "i_*\circ \cup"] & H^{n+r}(\mathbb{P}^r\times\mathbb{P}^n, \Omega_{\mathbb{P}^r\times\mathbb{P}^n}^{n+r}) \\
		J^{q-r,(q+1)m - (n+1)} \otimes J^{p-r,(p+1)m - (n+1)}  \arrow[r] \arrow[u]  & J^\rho \arrow[u, "\phi"] 
		\end{tikzcd} 
		\]
		commutes for all $p,q\in\mathbb{Z}_{\geq 0}$ such that $p+q = n+r-1$. 
	\end{introprop} 
	A slight generalization of an argument from \cite{KonnoVTPCI} shows that the map $\phi$ is in fact an isomorphism unless if $\mathcal{X}$ is odd dimensional, $r=1$ and $m=2$. The one exception will not matter for our purposes, because we know from the Motivic Gauss-Bonnet Theorem that the quadratic Euler characteristic of $\mathcal{X}$ is hyperbolic in this case. Furthermore, we will study a slight variant of the Jacobian ring, given by 
	$$\tilde{J} = k[Y_0,\cdots, Y_r, X_0,\cdots, X_n]/(Y_0F_0,\cdots, Y_rF_r, X_0\bar{F}_0,\cdots, X_n\bar{F}_n)$$ and show that $\tilde{J}^{\rho + (r+1,n+1)}$ is one dimensional. \\
	If we make some extra assumptions, we can compute the trace map. 
	\begin{introtheorem}[See Lemma \ref{construction of ctilde} and Theorem \ref{theorem trace}]
		Assume that $m+1$ is invertible in $k$, that $V(F_i)$ is smooth for all $i\in\{0,\cdots, r\}$ and that $V(F_0,\cdots, F_r)$ is smooth and of codimension $r+1$. Assume moreover that these assumptions remain true after setting any subset of the $X_i$  equal to zero. For $A\in J^{q-r,(q+1)m - (n+1)}$ and $B\in J^{p-r,(p+1)m - (n+1)}$, write $\omega_A\in H^q(\mathcal{X}, \Omega_{\mathcal{X}}^p)_{prim}$ and $\omega_B\in H^p(\mathcal{X}, \Omega_{\mathcal{X}}^q)_{prim}$ for their images. Write 
		$$G_0 = Y_0F_0, \cdots, G_r =  Y_rF_r, G_{r+1} = X_0\bar{F}_0,\cdots, G_{n+r+1} = X_n\bar{F}_n$$ and $Z_0 = Y_0,\cdots, Z_r = Y_r, Z_{r+1} = X_0,\cdots, Z_{n+r+1} = X_n$. Let $M$ be the Jacobian matrix given by~$\{\frac{\partial G_i}{\partial Z_j} \}_{i,j}$. Let $M_{i|j}$ be its minor with the $i$'th row and the $j$'th column missing. \\
		There exists a unique element $\tilde{C}\in k[Y_0,\cdots, Y_r, X_0,\cdots, X_n]$ such that $$(m+1)Y_iX_j\tilde{C} =(-1)^j \det(M_{0|j+r+1})Y_i + (-1)^{r+i}\det(M_{0|i})X_j$$ for all $i\in\{1,\cdots, r\}$ and $j\in \{r+1,\cdots, n+r+1\}$. Assume that we are not in the situation that $\dim(\mathcal{X})$ is odd, $r=1$ and $m=2$. Then the map
		$$\psi: J^\rho\to \tilde{J}^{\rho + (r+1,n+1)}, D\mapsto D\prod_{i=0}^rY_i\prod_{j=0}^nX_j$$ is an isomorphism. Therefore, we have that $\tilde{C} = \psi(C)$ for a unique $C\in J^\rho$. Write $AB = \lambda C$ in $J^\rho$ for some $\lambda\in k^*$. Then $$\text{Tr}(\omega_A\cup\omega_B) = (-1)^{r+1}m^{n+1}\binom{n+r}{r}\lambda.$$ 
	\end{introtheorem}
	Even though this does not give a completely explicit formula to compute the quadratic Euler characteristic, it may provide a useful algorithm in concrete cases.\\
	Also, we will see from the proof of the above theorem that if $\binom{n+r}{r}$ is invertible in $k$, we have that $Cm^{-n}\binom{n+r}{r}^{-1}$ has trace $1$. We call this the \textit{Scheja-Storch generator}, and conjecture that there is a way to define this without the assumption that $\binom{n+r}{r}$ is invertible in $k$. \\	
	As an application, we compute the quadratic Euler characteristic of a complete intersection of two generalized Fermat hypersurfaces.  
	\begin{introtheorem}[See Corollary \ref{qec of X}]
		Let $F_0 = \sum_{i=0}^na_iX_i^m, F_1 = \sum_{i=0}^nb_iX_i^m$ be two generalized Fermat hypersurfaces in $\mathbb{P}^n$. Assume that $a_i,b_i\in k^*$ are such that $a_ib_j-a_jb_i\neq 0$ for all $i\neq j$. Let $X = V(F_0,F_1)$. Then 
		$$\chi(X/k) = \begin{cases}
		B_{n,m}\cdot H &\text{ if } n \text{ is odd}\\
		B_{n,m}\cdot H + \langle 1 \rangle&\text{ if } n \text{ is even, } m \text{ odd}\\
		B_{n,m}\cdot H + \langle 1 \rangle + \sum_{i=0}^n\langle \prod_{j\neq i}(a_ib_j-a_jb_i) \rangle &\text{ if } n,m\text{ are even} 
		\end{cases} $$
		where $B_{n,m}\in\mathbb{Z}$ is given by
		$$B_{n,m} = \begin{cases}
		\frac{1}{2}\deg(c_{n-2}(T_{X})) &\text{ if $n$ odd} \\
		\frac{1}{2}\deg(c_{n-2}(T_{X})) - 1 &\text{ if $n$ even, $m$ odd} \\
		\frac{1}{2}\deg(c_{n-2}(T_{X})) - n - 1 &\text{ if } n,m \text{ even} 
		\end{cases} $$
	\end{introtheorem}
	\noindent The results presented here also form a chapter in the authors thesis, which was submitted on May 30'th, 2023. \\
	We note that in the paper \cite{Cox} by Cox and Batyrev, there is a more general construction of an isomorphism between primitive Hodge cohomology groups and certain graded parts of a Jacobian ring in the setting of toric varieties. This follows the methods of \cite{CarlsonIVHS} as do we, but they do not consider the multiplicative structure. \\
	We also point out that recently, Villaflor has independently proven a result which is similar to our Theorem \ref{theorem trace}, in the setting of toric varieties, see~\cite{Loyola}. Villaflor's result is more general than ours in the sense that it includes hypersurfaces of different degrees as well, and more limited in the sense that the arguments and constructions are made in the setting of varieties over the complex numbers. \\
	It would be interesting to extend the results in the present paper to the case where the hypersurfaces do not necessarily have the same degrees. This might be possible by extending the above arguments to the situation where $\mathbb{P}^r$ is a replaced by a weighted $r$-dimensional projective space, and will be explored in future work. 
	
	\subsection*{Structure} 
	Section \ref{section introduction to qec's} contains the definition of a quadratic Euler characteristic and some of its basic properties, after which we give a more detailed summary of the results in \cite{LevineECHWHH}. Section \ref{section cohomology} is devoted to some results on cohomology of differential forms and first Chern classes which will be needed later on. Then in Section \ref{section isomorphism from jacobian ring} we construct isomorphisms from bigraded parts of the Jacobian ring to primitive cohomology. In Section \ref{section two products} we compare the cup product and the product in the Jacobian ring, and we show that $J^\rho$ and $\tilde{J}^{\rho + (r+1,n+1)}$ are one dimensional, after which we compute the trace in Section~\ref{section computing trace}. Finally, we work out the example of intersecting two generalized Fermat hypersurfaces of the same degree in Section~\ref{section Fermat}. 
	
	\subsection*{Acknowledgements} 
	I deeply thank my advisor Marc Levine for suggesting this project to me and for all of his invaluable help over the last years. I also want to thank my coadvisor Jochen Heinloth, for all of his time and useful ideas, and for a lot of good advice. And I wish to thank V. Srinivas for helpful discussions. I also wish to thank Lukas Br\"oring for helping me straighten out my thoughts on some of the more technical arguments. \\
	This work was funded by the RTG Graduiertenkolleg 2553. 
	
	\subsection*{Notation}
	Throughout, let $k$ be a perfect field which is not of characteristic $2$. 
	
	\tableofcontents 
	
	\section{Quadratic Euler characteristics}\label{section introduction to qec's}
	In this section we give the definition and some basic properties of quadratic Euler characteristics and then discuss the computation of the quadratic Euler characteristic of a smooth projective hypersurface from \cite{LevineECHWHH}.  
	
	\subsection{Quadratic Euler characteristics}
	We give the definition of a quadratic Euler characteristic following the one by Levine in \cite[Section 1]{LevineAEGQF}. The quadratic Euler characteristic will be a particular case of a more general definition of Euler characteristic, introduced by Dold and Puppe in \cite{DoldPuppe}. \\
	Let $\mathcal{C}$ be a symmetric monoidal category, and denote $\otimes: \mathcal{C}\times\mathcal{C}\to\mathcal{C}$ for the tensor product and $1\in\mathcal{C}$ for the unit. Let $\tau$ be the symmetry isomorphism from the tensor product $\otimes$ to $\otimes\circ t$, where $t:\mathcal{C}\times\mathcal{C}\to\mathcal{C}\times\mathcal{C}$ is the usual symmetry given by $t(a,b) = (b,a)$. The following definition is taken from \cite[Definition 1.2 and Theorem 1.3]{DoldPuppe}. 
	\begin{defn}
		An object $X\in\mathcal{C}$ is \textit{strongly dualizable} if there exists an object $X^\vee\in\mathcal{C}$ and morphisms $\delta_X:1\to X\otimes X^\vee$ and $\text{ev}_X: X^\vee\otimes X\to 1$ in $\mathcal{C}$ such that the compositions 
		$$X\cong 1\otimes X \xrightarrow{\delta_X\otimes\text{Id}} X\otimes X^\vee\otimes X \xrightarrow{\text{Id}\otimes\text{ev}_X}X\otimes 1\cong X$$
		and 
		$$X^\vee\cong X^\vee\otimes 1 \xrightarrow{\text{Id}\otimes\delta_X} X^\vee\otimes X\otimes X^\vee \xrightarrow{\text{ev}_X\otimes\text{Id}}1\otimes X^\vee\cong X^\vee$$
		are the identity morphisms. 
	\end{defn}
	\begin{rem} 
		If $X$ is strongly dualizable, the triple $(X^\vee, \delta_X, \text{ev}_X)$ is unique up to unique isomorphism. We usually call $X^\vee$ the \textit{dual} of $X$, with the morphisms $\delta_X$ and $\text{ev}_X$ being understood.
	\end{rem} 
	Now let $X\in\mathcal{C}$ be a strongly dualizable object. The following definition is a special case of \cite[Definition 4.1]{DoldPuppe}.
	\begin{defn}
		The \textit{categorical Euler characteristic} of $X$ is the composition 
		$$1\xrightarrow{\delta_X}X\otimes X^\vee \xrightarrow{\tau} X^\vee\otimes X \xrightarrow{\text{ev}_X}1.$$
	\end{defn}
	To $k$, we can associate the \textit{motivic stable homotopy category} $\text{SH}(k)$, see for instance Morel's book \cite{MorelATF} or Hoyois' paper \cite{HoyoisSOEMHT} for its construction and properties. We have that $\text{SH}(k)$ is a symmetric monoidal category, with the ``smash product" as its tensor product. For a smooth projective scheme $X$ over $k$, we have the suspension spectrum $\Sigma_T^\infty X_+\in\text{SH}(k)$ and this is a strongly dualizable object, see for instance \cite[Theorem 5.22 and Corollary 6.13]{HoyoisSOEMHT}.  
	\begin{defn}\label{definition GW}
		The \textit{Grothendieck-Witt ring} $\text{GW}(k)$ of $k$ is the group completion of the monoid (under orthogonal direct sum) of isometry classes of nondegenerate quadratic forms over~$k$. 
	\end{defn}
	\begin{rem}\label{remark GW} 
		One can think of $\text{GW}(k)$ as the group generated by the forms $$\langle a \rangle: x\mapsto ax^2$$ for $a\in k^*$ modulo the relations
		\begin{itemize}
			\item $\langle ab^2 \rangle  =\langle a \rangle$ for $a,b\in k^*$. 
			\item $\langle a\rangle + \langle b \rangle = \langle a+b\rangle + \langle ab(a+b)\rangle$ for $a,b,a+b\in k^*$. 
			\item $\langle a \rangle + \langle -a \rangle = \langle 1 \rangle + \langle -1 \rangle$ for $a\in k^*$. 
		\end{itemize}
		Note that $\langle a \rangle \langle b \rangle =\langle ab \rangle$ for $a,b\in k^*$, by definition.\\
		This presentation originally goes back to Witt, see \cite[Section 1]{WittQF}. In the form above, it is \cite[Lemma 2.9]{MorelATF}, where the result is deduced from the statement for Witt rings, see \cite[Lemma (1.1) in Chapter 4]{MilnorSBF}. 
	\end{rem} 
	\begin{defn}
		The form $H = \langle 1 \rangle + \langle -1 \rangle$ is called the \textit{hyperbolic form}. 
	\end{defn}
	By a deep result of Morel (see \cite[Theorem 6.4.1]{MorelIAHT}) we have that $$\text{End}(1_{\text{SH}(k)})\cong \text{GW}(k).$$ 
	Combining the above, the quadratic Euler characteristic of a smooth projective scheme over $k$ can now be defined as follows. 
	\begin{defn}
		The \textit{quadratic Euler characteristic} $\chi(X/k)\in\text{GW}(k)$ of a smooth projective scheme $X$ over $k$ is the categorical Euler characteristic of $X$ in $\text{SH}(k)$. 
	\end{defn} 
	Quadratic Euler characteristics satisfy several nice relations, which one can find in e.g. \cite{LevineAEGQF}. One property which we will need is the following. 
	\begin{prop}[See \cite{LevineAEGQF}, Proposition 1.4(3)]\label{cut and paste proposition}
		Let $X$ be a smooth projective scheme over $k$ and let $Z$ be a smooth closed subscheme of pure codimension $c$ with complement $U$. Then 
		$$\chi(X/k) = \chi(U/k) + \langle -1 \rangle^c\chi(Z/k).$$
	\end{prop} 
	\begin{rem}
		Even though $U$ in the above statement is not projective, it has a quadratic Euler characteristic in $\text{GW}(k)$. Namely, it has been proven by Riou in \cite{RiouAECTAF} that a quasi-projective scheme is dualizable in the stable motivic homotopy category if we invert the characteristic of $k$. If $\text{char}(k) = p>0$ we have an injective morphism $\text{GW}(k)\to \text{GW}(k)[p^{-1}]$ and one can show that the categorical Euler characteristic always lands in $\text{GW}(k)$. See \cite[Remark 1.1.2]{LevineAEGQF} for more details. 
	\end{rem}
	\begin{exa}[See \cite{LevineAEGQF}, Proposition 1.4(4)] \label{Euler char of projective space}
		We have that $$\chi(\mathbb{P}^n/k) = \sum_{i=0}^{n}\langle -1\rangle^i.$$ One way to prove this is to use Proposition \ref{cut and paste proposition} together with induction on $n$ and the fact that $\mathbb{A}^n$ is equivalent to a point in $\text{SH}(k)$. Note that this is a multiple of $H$ if $n$ is odd. Also, the rank of this form is $n+1$ (which is the topological Euler characteristic of complex projective space) and its signature is either $0$ or $1$, depending on the parity of $n$ (which is the topological Euler characteristic of real projective space). 
	\end{exa}	
	\begin{rem} 
		This is true in general: if $k\subset\mathbb{R}\subset\mathbb{C}$ then we have that the rank of $\chi(X/k)$ is equal to the topological Euler characteristic of $X(\mathbb{C})$ while the signature is equal to the topological Euler characteristic of $X(\mathbb{R})$. See \cite[Remark 1.4.1]{LevineAEGQF}. 
	\end{rem} 
	\begin{rem} 
		For a smooth quasi-projective scheme $U$ over $k$ we have that $\chi(\mathbb{P}^n\times U/k) = \chi(\mathbb{P}^n/k)\chi(U/k)$ by \cite[Proposition 1.4(4)]{LevineAEGQF}. 
	\end{rem} 
	
	\subsection{The quadratic Euler characteristic of a smooth projective hypersurface}
	\begin{nota}
		For a scheme $X$ over $k$, we denote the sheaf of differential forms on $X$ over $k$ by $\Omega_X$. We write $\Omega_X^q = \wedge^q\Omega_X$ for $q\in\mathbb{Z}_{\geq 0}$. 
	\end{nota}
	Note that by \cite[Exercise II.4.5]{HartshorneAG}, for a smooth projective scheme $Y$ we have that $H^1(Y,\mathcal{O}_Y^*)\cong \text{Pic}(Y)$, where $\text{Pic}(Y)$ is the Picard group of $Y$, i.e. the group of line bundles on $Y$. There is the canonical ``dlog morphism" 
	$$\mathcal{O}_{XY}^*\to \Omega_{Y}, f\mapsto \frac{df}{f}$$ inducing the map 
	$$c_1: H^1(Y, \mathcal{O}_Y^*)\to H^1(Y,\Omega_Y).$$  
	\begin{defn}
		The \textit{first Chern class} of a line bundle $L$ on a scheme $Y$ is the image $c_1(L)\in H^1(Y,\Omega_Y)$ of the class of $L$ under the above map. 
	\end{defn}
	\begin{nota} 
		We write $c_1(L)^i\in H^i(Y, \Omega_{Y}^i)$ for the $i$-fold cup product of $c_1(L)$ with itself. 
	\end{nota}  
	An important computational tool is the motivic Gauss-Bonnet Theorem for $\text{SL}$-oriented cohomology theories proven by Levine and Raksit in their paper \cite{LevineGB}. A more general motivic Gauss-Bonnet theorem has been proven by D\'eglise, Jin and Khan in \cite{DegliseFCMHT}, and the theorem of Levine-Raksit can be viewed as a special case of their statement. A version where the scheme does not need to be smooth has been proven by Azouri, see \cite{AzouriMCCSS}. \\
	We do not state the theorem in all of its generality here, but rather one of its applications which provides a way to compute a quadratic Euler characteristic in practice. 
	\begin{theorem}[See \cite{LevineGB}, Corollary 8.7]\label{Motivic Gauss Bonnet}
		Let $X$ be a smooth projective scheme over~$k$. Then: 
		\begin{itemize}
			\item If $\dim(X)$ is odd, then $\chi(X/k) = C\cdot H$ for some $C\in\mathbb{Z}$,
			where $H$ is the hyperbolic form. 
			\item If $\dim(X) = 2n$ is even, then $\chi(X/k) = C\cdot H + Q$ for some $C\in\mathbb{Z}$, where $Q$ is the quadratic form given by the composition 
			$$H^n(X,\Omega_X^n)\times H^n(X,\Omega_X^n)\xrightarrow{\cup}H^{2n}(X,\Omega_X^{2n})\xrightarrow{\text{Trace}}k.$$
			Here the first map is the cup product on cohomology. 
		\end{itemize}
	\end{theorem}  
	\begin{rem}
		By \cite[Theorem 5.3]{LevineGB} the rank of $\chi(X/k)$ is equal to the degree of $c_n(T_X)$ where $T_X$ is the tangent bundle of $X$. This gives a way to determine the constant $C$ in the above theorem in practice. There is also a formula for $C$ in terms of dimensions of cohomology groups of $\Omega_X^q$ given in \cite[Corollary 8.7]{LevineGB}. 
	\end{rem}
	The form $Q$ has been computed successfully in the case of hypersurfaces by Levine, Lehalleur and Srinivas in \cite{LevineECHWHH}, using inspiration from the paper \cite{CarlsonIVHS} by Carlson and Griffiths. We now summarize their strategy. \\
	Consider a smooth hypersurface $X = V(F)\subset\mathbb{P}^n$ where $F\in k[X_0,\cdots, X_n]$ is a homogeneous polynomial of degree $m\in\mathbb{Z}_{\geq 2}$. Assume that the characteristic of $k$ is coprime to $m$. 
	\begin{defn} 
		The \textit{Jacobian ring of $X$} is $$J_X = k[X_0,\cdots, X_n]/\left(\frac{\partial F}{\partial X_0},\cdots, \frac{\partial F}{\partial X_n}\right).$$ 
	\end{defn}
	Note that $J_X$ has a natural grading induced by the grading of $k[X_0,\cdots, X_n]$. Furthermore, $J_X$ is a finite dimensional $k$-algebra. The top nonzero graded part is $J_X^{(m-2)(n+1)}$, which is a one dimensional vector space over $k$; for a proof, see \cite[Lemma 4]{KassWickelgren}, where the result is deduced from the fact that $J_X$ is Gorenstein together with the proof of \cite[(4.7) Korrolar]{SchejaStorchSF}. 
	\begin{constr} \label{construction SS generator}
		There is a canonical choice of generator $e_F$ of $J_X^{(m-2)(n+1)}$ called the \textit{Scheja-Storch generator}. Namely, as $m\geq 2$, for $i\in\{0,\cdots, n\}$ we can write $$\frac{\partial F}{\partial X_i} = \sum_{j=0}^n a_{ij}X_j$$ for some (non-unique) $a_{ij}\in k[X_0,\cdots, X_n]$. One defines $e_F = \det((a_{ij})_{i,j})$. One can show that this is independent of the choice of $a_{ij}$, see \cite[(1.2)($\alpha$)]{SchejaStorchSF}. 
	\end{constr} 
	\begin{exa}
		Let $F = \sum_{i=0}^na_iX_i^m$ where $a_0,\cdots, a_n\in k^*$. Then $X$ is a \textit{generalized Fermat hypersurface}. We have that $\frac{\partial F}{\partial X_i} = ma_iX_i^{m-1}$ and so one computes $e_F$ as $$e_F = m^{n+1}\prod_{i=0}^na_iX_i^{m-2}.$$  
	\end{exa}
	Let $i: X\to\mathbb{P}^n$ be the natural inclusion. This induces a pushforward map $i_*: H^q(X,\Omega_X^{p})\to H^{q+1}(\mathbb{P}^n,\Omega_{\mathbb{P}^n}^{p+1})$ for all $p,q\in\mathbb{Z}_{\geq 0}$ as defined by Srinivas in~\cite{SrinivasGMCSHC}.
	\begin{defn} \label{definition primitive cohomology hypersurface}
		The \textit{primitive cohomology of $X$ with respect to $p,q\in\mathbb{Z}_{\geq 0}$ such that $p+q = n-1$} is defined by $H^q(X,\Omega_X^{p})_{prim} = \ker(i_*)$.
	\end{defn} 
	\begin{rem}\label{remark primitive cohomology}
		Let $c_1(\mathcal{O}(1))\in H^1(X,\Omega_X)$ be the first Chern class of $\mathcal{O}(1)$. The Hard Lefschetz Theorem tells us that for $0<i \leq n-1$, the map 
		$$(-)\cup c_1(\mathcal{O}(1))^{i}:\bigoplus_{p+q = n-i}H^q(X,\Omega_X^p)\to \bigoplus_{p+q = n+i}H^q(X,\Omega_X^p)$$
		is an isomorphism. Classically, for $0\leq i \leq n-1$, the primitive cohomology of $X$ is defined to be the kernel of the morphism 
		$$(-)\cup c_1(\mathcal{O}(1))^{i+1}:\bigoplus_{p+q = n-i}H^q(X,\Omega_X^p)\to \bigoplus_{p+q = n+i + 2}H^q(X,\Omega_X^p).$$
		For $i=0$ and $p,q$ such that $p+q = n-1$, this definition coincides with the one above. To see this, note that multiplication with $c_1(\mathcal{O}(1))$ takes $H^q(X,\Omega_X^p)$ to $H^{q+1}(X,\Omega_X^{p+1})$. The pullback $i^*: H^{q+1}(\mathbb{P}^n,\Omega_{\mathbb{P}^n}^{p+1})\to H^{q+1}(X,\Omega_X^{p+1})$ is an isomorphism by the Weak Lefschetz Theorem, so we can view multiplication with $c_1(\mathcal{O}(1))$ as a morphism to $H^{q+1}(\mathbb{P}^n,\Omega_{\mathbb{P}^n}^{p+1})$. Now as $i^*i_*$ is multiplication with $c_1(\mathcal{O}(1))$, we see that the kernel coincides with $\ker(i_*)$. 
	\end{rem}
	One can show that $H^q(X,\Omega_X^{p})_{prim}=H^q(X,\Omega_X^{p})$ whenever $p\neq q$. In \cite{LevineECHWHH}, the authors prove the following result. 
	\begin{prop}[\cite{LevineECHWHH}, Proposition 3.2]
		For each $q\geq 0$, there is a canonical isomorphism $\psi_q: J_X^{(q+1)m-n-1}\to H^q(X,\Omega_X^{n-1-q})_{prim}$. 
	\end{prop} 
	This result originally goes back to Dolgachev, see \cite{DolgachevWPV}, and the characteristic zero case is due to Griffiths, see \cite{GriffithsOPOCRI}. In \cite{LevineECHWHH}, there is then a comparison of the cup product on cohomology with the usual ring multiplication of $J_X$, leading up to the following result. 
	\begin{prop}[\cite{LevineECHWHH}, Proposition 3.7]
		Consider $p,q\in\mathbb{Z}_{\geq 0}$ be such that $p+q = n-1$ and let $A\in J_X^{(q+1)m-n-1}$ and $B\in J_X^{(p+1)m-n-1}$. Let $$\omega_A = \psi_q(A)\in H^q(X,\Omega_X^{p})_{prim} \text{ and } \omega_B = \psi_p(B)\in H^p(X,\Omega_X^{q})_{prim}$$ be their images. Write $F_i = \frac{\partial F}{\partial X_i}$. Cover $\mathbb{P}^n$ by the open cover $\mathcal{U} = \{U_0,\cdots, U_n\}$ where $U_i = \{ F_i\neq 0\}$, and let $C^i(\mathcal{U},\Omega_{\mathbb{P}^n}^n)$ denote the $i$'th group in the \v{C}ech complex corresponding to $\mathcal{U}$. Furthermore, let $\bar{\omega} = \sum_{i=0}^n(-1)^iX_idX^i$ be the generator of $\Omega_{\mathbb{P}^n}^n(n+1)$, where we write $dX^i = dX_0\cdots dX_{i-1}dX_{i+1}\cdots dX_n$. \\
		Then the element $i_*(\omega_A\cup\omega_B)\in H^n(\mathbb{P}^n,\Omega_{\mathbb{P}^n}^n)$ is represented by 
		$$\frac{-mAB\bar{\omega}}{F_0\cdots F_n}\in C^n(\mathcal{U},\Omega_{\mathbb{P}^n}^n).$$
	\end{prop} 
	They then compare this with a representation of $c_1(\mathcal{O}(m))^n$ on the same cover to compute the trace, which yields the following result.  
	\begin{theorem}[\cite{LevineECHWHH}, Theorem 3.9]\label{Final theorem Levine}
		In the situation of the above proposition, suppose that $AB = \lambda e_F$ for $\lambda\in k^*$. Then 
		$$\text{Tr}(\omega_A\cup\omega_B) = -m\lambda.$$
	\end{theorem} 
	\begin{exa}
		In the case of a generalized Fermat hypersurface $X$ as before, if $n = 2p+1$ is odd, we need to calculate the form $Q$. One can show that $H^p(X,\Omega_X^p) = H^p(X,\Omega_X^p)_{prim}\oplus c_1(\mathcal{O}(1))^p$ and compute that $ c_1(\mathcal{O}(1))^p$ contributes a form $\langle m\rangle$. For the primitive cohomology, we evaluate the form on basis elements of $H^p(X,\Omega_X^p)_{prim}$, i.e. on the corresponding parts of the Jacobian ring. If $AB = \lambda e_F$ for some $\lambda\in k^*$ and two distinct basis elements $A$ and $B$, we also have that $BA = \lambda e_F$ and one can check that this yields a hyperbolic form. If $m$ is odd, there are no basis elements that square to a nonzero multiple of $e_F$. If $m$ is even, we have that $$\left(X_0^{\frac{m-2}{2}}\cdots X_n^{\frac{m-2}{2}}\right)^2 = \frac{e_F}{m^{n+1}\prod_{i=0}^na_i}.$$ Using Theorem \ref{Final theorem Levine}, this gives rise to the form $\langle -m\prod_{i=0}^na_i\rangle$. One can see from \cite[Theorem 5.3]{LevineGB} that the rank of $\chi(X/k)$ is equal to $\deg(c_n(T_X))$ where $T_X$ is the tangent bundle on $X$. Putting everything together, we find that 
		$$\chi(X/k) = \begin{cases}
		A_{n,m}\cdot H &\text{ if } n \text{ even}\\
		A_{n,m}\cdot H + \langle m \rangle &\text{ if } n, m \text{ odd}\\
		A_{n,m}\cdot H + \langle m \rangle + \langle -m\prod_{i=0}^na_i\rangle  &\text{ if } n \text{ odd},m \text{ even}
		\end{cases} $$
		for integers $A_{n,m}\in\mathbb{Z}$ given by 
		$$A_{n,m} = \begin{cases}
		\frac{1}{2}\deg(c_{n}(T_{X})) &\text{ if $n$ even} \\
		\frac{1}{2}\deg(c_{n}(T_{X})) - 1 &\text{ if } n,m \text{ odd} \\
		\frac{1}{2}\deg(c_{n}(T_{X})) - 2 &\text{ if $n$ odd, $m$ even} 	
		\end{cases} $$
		This is also \cite[Theorem 11.1]{LevineAEGQF}, but there it is proven in a different way, namely using Levine's quadratic Riemann-Hurwitz formula from \cite{LevineAEGQF}. 
	\end{exa}
	
	\section{Setup, cohomology of differential forms and primitive cohomology}\label{section cohomology}
	In the next sections, we will be working with the following setup. 
	\begin{nota} 
		Let $n,m,r\in\mathbb{Z}_{\geq 1}$ be such that $n\geq r+2$ and $m\geq 2$. Assume that $m$ is coprime to the characteristic of $k$, if $\text{char}(k)$ is positive. Let $F_0,\cdots, F_r\in k[X_0,\cdots, X_n]$ be homogeneous polynomials of the same degree $m$. Let $X = V(F_0,\cdots, F_r)\subset\mathbb{P}^n$ be the intersection of the $V(F_i)$ and assume that this is a smooth complete intersection. We define $F= Y_0F_0+ \cdots + Y_rF_r$ and consider the hypersurface $$\mathcal{X} = V(F)\subset\mathbb{P}^r\times\mathbb{P}^n.$$ 
		We write $i: \mathcal{X}\to \mathbb{P}^r\times\mathbb{P}^n$ for the inclusion. Note that $\mathcal{X}$ is of bidegree $(1,m)$ and that it has dimension $n+r-1$. We note that $\frac{\partial F}{\partial Y_i} = F_i$ for $i\in\{0,\cdots, r\}$ and write $\bar{F}_j = \frac{\partial F}{\partial X_j}$ for $j\in\{0,\cdots, n\}$. 
	\end{nota} 
	\begin{nota}
		We denote the canonical projections by $\pi_n: \mathbb{P}^r\times\mathbb{P}^n\to\mathbb{P}^n$ and $\pi_r: \mathbb{P}^r\times\mathbb{P}^n\to\mathbb{P}^r$. 
	\end{nota}
	\begin{rem} 
		Note that $\mathcal{X}$ is smooth: suppose $(y_0,\cdots, y_r,x_0,\cdots, x_n)\in \mathcal{X}$ is a point where all $F_i$ and $\bar{F}_j$ vanish, then $(y_0,\cdots, y_r,x_0,\cdots, x_n)$ is in $\mathbb{P}^r\times X$. As $X$ is smooth, the vectors $$\left(\frac{\partial F_i}{\partial X_0}(x),\cdots, \frac{\partial F_i}{\partial X_n}(x)\right)$$ for $i\in\{0,\cdots, r\}$ are linearly independent for any $x\in X$. Now as $$\bar{F}_j(y_0,\cdots, y_r,x_0,\cdots, x_n) = \sum_{i=0}^ry_i\frac{\partial F_i}{\partial X_j}(x_0,\cdots, x_n) = 0$$ for all $j\in\{0,\cdots, n\}$, we have that 
		$$\sum_{i=0}^ry_i\left(\frac{\partial F_i}{\partial X_0}(x_0,\cdots, x_n),\cdots, \frac{\partial F_i}{\partial X_n}(x_0,\cdots, x_n)\right)=0 $$
		and so $y_0 = \cdots = y_r = 0$, but this is impossible. 
	\end{rem} 
	\begin{rem} 
		Note that we have the two Euler equations
		$$F = \sum_{i=0}^rY_iF_i\text{ and } mF = \sum_{j=0}^nX_j\bar{F}_j.$$
	\end{rem} 
	We now observe that we can compute $\chi(X/k)$ from $\chi(\mathcal{X}/k)$. 
	\begin{lem}
		We have that $$\chi(\mathcal{X}/k) = \chi(\mathbb{P}^{r-1}/k)\chi(\mathbb{P}^n/k) + \langle -1\rangle^r\chi(X/k).$$
	\end{lem}
	\begin{proof}
		Let $U$ be the complement of $X$ in~$\mathbb{P}^n$ and let $\pi_n|_{\mathcal{X}}:\mathcal{X}\to\mathbb{P}^n$ be the restriction of $\pi_n$ to $\mathcal{X}$. Then $(\pi_n|_{\mathcal{X}})^{-1}(X) = \mathbb{P}^{r}\times X$ and $(\pi_n|_{\mathcal{X}})^{-1}(U)\to U$ is a Zariski locally trivial $\mathbb{P}^{r-1}$-bundle. Using \cite[Proposition 1.4(4)]{LevineAEGQF}, we have that 
		$$\chi((\pi_n|_{\mathcal{X}})^{-1}(U)/k) = \chi(\mathbb{P}^{r-1}/k)\cdot\chi(U/k)$$
		and 
		$$\chi(\mathbb{P}^r\times X/k) = \chi(\mathbb{P}^r/k)\chi(X/k).$$
		By Proposition \ref{cut and paste proposition}, we also have that 
		$$\chi(\mathbb{P}^n/k) = \chi(U/k) + \langle -1 \rangle^{r+1} \chi(X/k).$$
		Recalling Example \ref{Euler char of projective space}, this yields
		\begin{align*}
		\chi(\mathcal{X}/k) &= \chi(\mathbb{P}^{r-1}/k)\chi(U/k) + \langle -1 \rangle^{r} \chi(\mathbb{P}^r/k)\chi(X/k)\\
		&= \chi(\mathbb{P}^{r-1}/k)\chi(\mathbb{P}^n/k) + (\langle -1 \rangle^{r} \chi(\mathbb{P}^r/k) - \langle -1\rangle^{r+1}\chi(\mathbb{P}^{r-1}/k))\chi(X/k) \\
		&=  \chi(\mathbb{P}^{r-1}/k)\chi(\mathbb{P}^n/k) + \left(\langle -1 \rangle^{r} \sum_{i=0}^r\langle -1\rangle^{i} - \langle -1\rangle^{r+1}\sum_{i=0}^{r-1}\langle -1\rangle^{i}\right)\chi(X/k) \\
		&= \chi(\mathbb{P}^{r-1}/k)\chi(\mathbb{P}^n/k) + \langle -1\rangle^r\chi(X/k)
		\end{align*}
		as desired. 
	\end{proof} 
	In the coming sections, the strategy will be to adapt the arguments of \cite{LevineECHWHH} to a hypersurface in $\mathbb{P}^r\times\mathbb{P}^n$, using inspiration from \cite{TerasomaIVH}. In this section, we start by introducing two exact sequences which we will use in what follows, and we study the cohomology groups of differential forms for a product of projective spaces, which will be needed later on. We also study first Chern classes of line bundles on $\mathbb{P}^r\times\mathbb{P}^n$ and primitive cohomology. 
	\begin{nota} 
		The Picard group of~$\mathbb{P}^r\times\mathbb{P}^n$ is isomorphic to~$\mathbb{Z}\oplus\mathbb{Z}$ with generators coming from the canonical sheaves~$\mathcal{O}_{\mathbb{P}^r}(a)$ and~$\mathcal{O}_{\mathbb{P}^n}(b)$ for $a,b\in\mathbb{Z}$. For a sheaf $\mathcal{F}$ on $\mathbb{P}^r\times\mathbb{P}^n$, we denote $$\mathcal{F}(a,b) =\mathcal{F}\otimes \pi_r^*\mathcal{O}_{\mathbb{P}^r}(a)\otimes \pi_n^*\mathcal{O}_{\mathbb{P}^n}(b).$$  There are thus canonical sheaves of the form~$\mathcal{O}(a,b)$ on~$\mathbb{P}^r\times\mathbb{P}^n$. \\
		For $a\in\mathbb{Z}_{\geq 0}$, we write $\mathcal{O}(a\mathcal{X})$ for the sheaf with sections having poles of order at most $a$ on $\mathcal{X}$, which are regular everywhere else.  We set $\mathcal{F}(a\mathcal{X}) = \mathcal{F}\otimes\mathcal{O}(a\mathcal{X})$. For $a<0$, we write $\mathcal{F}(a\mathcal{X}) = \mathcal{F}\otimes \mathcal{I}_{\mathcal{X}}^{-a}$ where $\mathcal{I}_{\mathcal{X}}$ is the ideal sheaf of $\mathcal{X}$. 
	\end{nota} 
	
	\subsection{Two exact sequences}
	Recall from e.g. \cite[Theorem II.8.17]{HartshorneAG} that there is an exact sequence 
	\begin{equation} \label{hyperplane exact sequence}
	0\to \mathcal{O}_{\mathcal{X}}(-\mathcal{X})\xrightarrow{dF/F\wedge(-)} i^*\Omega_{\mathbb{P}^r\times\mathbb{P}^n}\to \Omega_{\mathcal{X}}\to 0.
	\end{equation} 
	Here, the second map is the natural surjection. There is another useful exact sequence 
	\begin{equation} \label{residue exact sequence}
	0\to \Omega_{\mathbb{P}^r\times\mathbb{P}^n}^{p+1}\to \Omega_{\mathbb{P}^r\times\mathbb{P}^n}^{p+1}(\log(\mathcal{X}))\xrightarrow{res_X} i_*\Omega_{\mathcal{X}}^{p}\to 0
	\end{equation}
	for any $p\in\mathbb{Z}_{\geq 0}$ which is called the \textit{residue sequence}. We will need the following statement in the next sections. 
	\begin{lem}\label{ilowerstar is boundary}
		For $p,q\in\mathbb{Z}_{\geq 0}$, the boundary map $$\delta^{p,q}: H^q(\mathcal{X},\Omega^p_{\mathcal{X}})\to H^{q+1}(\mathbb{P}^r\times\mathbb{P}^n,\Omega_{\mathbb{P}^r\times\mathbb{P}^n}^{p+1})$$ induced from the long exact cohomology sequence of the exact sequence (\ref{residue exact sequence}) coincides with the pushforward map $i_*:H^q(\mathcal{X},\Omega^p_{\mathcal{X}})\to H^{q+1}(\mathbb{P}^r\times\mathbb{P}^n,\Omega_{\mathbb{P}^r\times\mathbb{P}^n}^{p+1})$, which is again the pushforward as defined in \cite{SrinivasGMCSHC}. 
	\end{lem}
	The proof works exactly the same as that of \cite[Lemma 2.2]{LevineECHWHH}. 
	\begin{rem}\label{remark wedge products for exact sequences}
		We will use in what follows that for an exact sequence of vector spaces $0\to V\to W\to Z\to 0$ with $V$ one dimensional, the induced sequence $0\to V\otimes \wedge^{k-1}Z\to \wedge^kW\to \wedge^kZ\to 0$ is again exact for any $k\in\mathbb{Z}_{\geq 1}$. For an exact sequence  $0\to V\to W\to Z\to 0$ with $Z$ a line bundle, we similarly have that the sequence $0\to \wedge^kV\to\wedge^k W\to \wedge^{k-1}V\otimes Z\to 0$ is again exact for any~$k\in\mathbb{Z}_{\geq 1}$. 
	\end{rem} 
	
	\subsection{Cohomology of differential forms}
	We will need Bott's theorem in what follows. 
	\begin{theorem}[Bott's theorem for projective space, see \cite{DolgachevWPV}, Theorem 2.3.2]\label{Botts theorem}
		Let $m\in\mathbb{Z}$. The cohomology of $\Omega_{\mathbb{P}^n}^q(m)$ satisfies $H^p(\mathbb{P}^n,\Omega_{\mathbb{P}^n}^q(m)) = 0$ for:
		\begin{itemize}
			\item $p>0$ and $m\geq q-n, m\neq 0$.
			\item $p>0, m=0$ and $p\neq q$.
			\item $p=0$ and $m\leq q$, except for $m = p = q = 0$. 
		\end{itemize}
	\end{theorem} 
	We will also make use of the following statement, which is useful to apply Bott's theorem to twisted sheaves of differentials on $\mathbb{P}^r\times\mathbb{P}^n$. 
	\begin{prop}\label{Botts theorem for general products}
		Let~$p,q\in \mathbb{Z}_{\geq 0}$ and let~$a,b\in\mathbb{Z}$. Then 
		\begin{align*} 
		H^p(\mathbb{P}^r\times\mathbb{P}^n, \Omega_{\mathbb{P}^r\times\mathbb{P}^n}^q(a,b)) &= \bigoplus_{i+j = q}\bigoplus_{k + l = p} H^k(\mathbb{P}^r, \Omega_{\mathbb{P}^r}^i(a))\otimes H^l(\mathbb{P}^n, \Omega_{\mathbb{P}^n}^j(b)).
		\end{align*} 
	\end{prop}
	\begin{proof}[Proof of Proposition \ref{Botts theorem for general products}]	
		We have that $\Omega_{\mathbb{P}^r\times\mathbb{P}^n} \cong \pi_r^*\Omega_{\mathbb{P}^r}\oplus\pi_n^*\Omega_{\mathbb{P}^n}$, by e.g. \cite[Exercise II.8.3]{HartshorneAG}. This implies that 
		$$\Omega_{\mathbb{P}^r\times\mathbb{P}^n}^q \cong \bigoplus_{i+j = q}\pi_r^*\Omega_{\mathbb{P}^r}^i\otimes\pi_n^*\Omega_{\mathbb{P}^n}^j$$ and tensoring with~$\mathcal{O}_{\mathbb{P}^r\times\mathbb{P}^n}(a,b) = \pi_r^*\mathcal{O}_{\mathbb{P}^r}(a)\otimes \pi_n^*\mathcal{O}_{\mathbb{P}^n}(b)$ yields that 
		$$\Omega_{\mathbb{P}^r\times\mathbb{P}^n}^q(a,b) \cong \bigoplus_{i+j = q}\pi_r^*\Omega_{\mathbb{P}^r}^i(a)\otimes\pi_n^*\Omega_{\mathbb{P}^n}^j(b).$$
		Now using \cite[Tag 0BED]{stacks-project}, we have that
		\begin{align*} 
		H^p(\mathbb{P}^r\times\mathbb{P}^n, \Omega_{\mathbb{P}^r\times\mathbb{P}^n}^q(a,b)) &= \bigoplus_{i+j = q}H^p(\mathbb{P}^r\times\mathbb{P}^n, \pi_r^*\Omega_{\mathbb{P}^r}^i(a)\otimes\pi_n^*\Omega_{\mathbb{P}^n}^j(b)) \\
		&= \bigoplus_{i+j = q}\bigoplus_{k + l = p} H^k(\mathbb{P}^r, \Omega_{\mathbb{P}^r}^i(a))\otimes H^l(\mathbb{P}^n, \Omega_{\mathbb{P}^n}^j(b)) 
		\end{align*} 
		which is the desired result. 
	\end{proof}
	
	\subsection{Primitive cohomology}  
	The primitive cohomology of $\mathcal{X}$ is defined in the same way as in the case of a smooth hypersurface in $\mathbb{P}^n$, see Definition \ref{definition primitive cohomology hypersurface}. 
	\begin{defn}
		Let $p,q\geq 0$ be such that $p+q = n+r-1$. The \textit{primitive cohomology} of $\mathcal{X}$ with respect to $p,q$ is $H^p(\mathcal{X}, \Omega_{\mathcal{X}}^q)_{\text{prim}} =\ker(i_*)\subset H^p(\mathcal{X}, \Omega_{\mathcal{X}}^q)$ where $i_*$ is the pushforward $i_*: H^p(\mathcal{X}, \Omega_{\mathcal{X}}^q)\to H^{p+1}(\mathbb{P}^r\times\mathbb{P}^n, \Omega_{\mathbb{P}^r\times\mathbb{P}^n}^{q+1})$ as defined in \cite{SrinivasGMCSHC}. 
	\end{defn}
	\begin{rem} 
		The same arguments as in Remark \ref{remark primitive cohomology} with $\mathcal{O}(1)$ replaced by $\mathcal{O}(1,m)$ show that this definition coincides with the classical definition of primitive cohomology. 
	\end{rem}
	The following result is probably standard, but we include a proof here for the reader's convenience. 
	\begin{lem}\label{ilowerstar surjective lemma}
		The map $i_*: H^p(\mathcal{X}, \Omega_{\mathcal{X}}^q)\to H^{p+1}(\mathbb{P}^r\times\mathbb{P}^n, \Omega_{\mathbb{P}^r\times\mathbb{P}^n}^{q+1})$ is surjective if either $p\neq q$ or $p=q$ and $p\geq r$. Furthermore, if $p\neq q$, we have that $H^{p+1}(\mathbb{P}^r\times\mathbb{P}^n, \Omega_{\mathbb{P}^r\times\mathbb{P}^n}^{q+1}) = 0$, so that $i_*$ is the zero map. 
	\end{lem}
	\begin{proof}
		Using Proposition \ref{Botts theorem for general products}, we have that 
		$$H^{p+1}(\mathbb{P}^r\times\mathbb{P}^n, \Omega_{\mathbb{P}^r\times\mathbb{P}^n}^{q+1})\cong \bigoplus_{i+j = p+1}\bigoplus_{k+l = q +1} H^i(\mathbb{P}^r,\Omega_{\mathbb{P}^r}^k)\otimes H^j(\mathbb{P}^n,\Omega_{\mathbb{P}^n}^l).$$
		Suppose $p\neq q$. Using Theorem \ref{Botts theorem} we see that all terms in the above sum are zero except for those with $i=k$ and $j=l$. However, for such terms we have that $i+j = k+l$ which cannot be true as $p\neq q$. This implies that $H^{p+1}(\mathbb{P}^r\times\mathbb{P}^n, \Omega_{\mathbb{P}^r\times\mathbb{P}^n}^{q+1}) = 0$. Therefore, $i_*$ is the zero map, so in particular, $i_*$ is surjective. \\		
		If $p=q$, then we see that
		$H^{p+1}(\mathbb{P}^r\times\mathbb{P}^n, \Omega_{\mathbb{P}^r\times\mathbb{P}^n}^{p+1})\cong k^{N_{p+1}}$ for some $N_{p+1}\in\mathbb{Z}_{\geq 0}$. A basis is given by $$\alpha_i^{p+1} = c_1(\mathcal{O}_{\mathbb{P}^r}(1))^i\otimes c_1(\mathcal{O}_{\mathbb{P}^n}(1))^j$$ where $i+j = p+1$, $i\leq r$ and $j\leq n$. 
		We now assume that $p\geq r$, then we have that $N_{p+1}\leq N_p$. As $$i_*i^*1 = c_1(\mathcal{O}(1,m)) = c_1(\mathcal{O}(1,0)) + mc_1(\mathcal{O}(0,1))$$ using the projection formula we see that 
		\begin{align*}
		i_*i^*\alpha_i^{p} &= i_*(i^*\alpha_i^p \otimes i^*1) \\
		&= \alpha_i^p \otimes (c_1(\mathcal{O}(1,0)) + mc_1(\mathcal{O}(0,1))) \\
		&= \alpha_{i+1}^{p+1} + m\alpha_i^{p+1}
		\end{align*}
		and the latter term is nonzero as $m$ is coprime to the characteristic of $k$. We see from this that the matrix of $i_*i^*$ has an $N_{p+1}\times N_{p+1}$ minor with determinant a power of $m$, hence invertible. We conclude from this that $i_*$ is surjective. 
	\end{proof}
	\begin{cor}
		For $p,q\in\mathbb{Z}_{\geq 0}$ such that $p+q=n+r-1$, we have that $H^p(\mathcal{X}, \Omega_{\mathcal{X}}^q)_{prim} = H^p(\mathcal{X}, \Omega_{\mathcal{X}}^q)$ as long as $p\neq q$. 
	\end{cor}
	
	\section{Isomorphism from the Jacobian ring to primitive cohomology}\label{section isomorphism from jacobian ring}
	We keep the notation that was set up at the beginning of the previous section. 
	\begin{defn} 
		The \textit{Jacobian ring} of $F$ is given by 
		$$J = k[Y_0,\cdots, Y_r, X_0,\cdots, X_n]/\left(F_0,\cdots, F_r, \bar{F}_0,\cdots, \bar{F}_n \right)$$
	\end{defn} 
	This definition is taken from \cite{TerasomaIVH}. Note that $J$ has a natural bigrading where $J^{a,b}$ has degree $a$ in the variables $Y_i$ and $b$ in the variables $X_i$. In this section, we will follow the argumentation of \cite{LevineECHWHH} with this Jacobian ring in order to show that certain graded pieces are isomorphic to certain primitive Hodge cohomology groups.  
	\begin{nota} 
		The line bundle $\Omega_{\mathbb{P}^r\times\mathbb{P}^n}^{n+r}(r+1,n+1)$ is isomorphic to $\mathcal{O}_{\mathbb{P}^r\times\mathbb{P}^n}$ with global generator $\omega\wedge\bar{\omega}$ where 
		$$\omega = \sum_{i=0}^r(-1)^i Y_idY^i$$ 
		where we write $dY^i = dY_0\wedge \cdots \wedge \hat{dY_i}\wedge\cdots \wedge dY_r$ and 
		$$\bar{\omega} = \sum_{j=0}^n(-1)^j X_jdX^j$$
		where we write $dX^j = dX_0\wedge \cdots \wedge \hat{dX_j}\wedge\cdots \wedge dX_n$. We will more generally write $dX^{i_0,\cdots, i_k}$ to mean the wedge product of all $dX_i$ with $dX_{i_0},\cdots,dX_{i_k}$ removed, and use the notation $dY^{i_0,\cdots, i_k}$ similarly.  
	\end{nota} 	
	\begin{nota} 
		We fix integers $p,q\in\mathbb{Z}_{\geq 0}$ satisfying $p+q = n+r-1$.  
	\end{nota} 
	\begin{constr} 
		Consider the exact sequence (\ref{hyperplane exact sequence}). For $j\in\{0,\cdots, q-1\}$, we take a wedge product and use Remark \ref{remark wedge products for exact sequences} to obtain the exact sequence 
		$$0\to \Omega_{\mathcal{X}}^{n+r-j-2}(-\mathcal{X}) \to i^*\Omega_{\mathbb{P}^r\times\mathbb{P}^n}^{n+r-1-j} \to \Omega_{\mathcal{X}}^{n+r-1-j}\to 0.$$
		Now twisting by the line bundle $i^*\mathcal{O}_{\mathbb{P}^r\times\mathbb{P}^n}((q-j)\mathcal{X})$, we obtain the exact sequences 
		\begin{equation} \label{Levine lemma 3.1(1) exact sequence proof}
		0\to \Omega_{\mathcal{X}}^{n+r-j-2}((q-j-1)\mathcal{X}) \to i^*\Omega_{\mathbb{P}^r\times\mathbb{P}^n}^{n+r-1-j}((q-j)\mathcal{X}) \to \Omega_{\mathcal{X}}^{n+r-1-j}((q-j)\mathcal{X})\to 0
		\end{equation} 
		where the first map is induced by $dF/F\wedge(-)$. Patching those together, we see that there is an exact sequence 
		\begin{align*}
		0\to \Omega_{\mathcal{X}}^p\xrightarrow{dF/F\wedge (-)} i^*\Omega_{\mathbb{P}^r\times\mathbb{P}^n}^{p+1}(\mathcal{X}) &\xrightarrow{dF/F\wedge (-)} i^*\Omega_{\mathbb{P}^r\times\mathbb{P}^n}^{p+2}(2\mathcal{X}) \xrightarrow{dF/F\wedge (-)}\cdots \\
		\cdots \to &i^*\Omega_{\mathbb{P}^r\times\mathbb{P}^n}^{n+r-1}(q\mathcal{X}) \xrightarrow{\pi_q}\Omega_{\mathcal{X}}^{n+r-1}(q\mathcal{X})\to 0
		\end{align*}
		of sheaves on $\mathcal{X}$. 
	\end{constr}
	\begin{nota} 
		Let $\mathcal{C}(p)$ be the complex 
		\begin{align*}
		0\to i^*\Omega_{\mathbb{P}^r\times\mathbb{P}^n}^{p+1}(\mathcal{X}) \xrightarrow{dF/F\wedge (-)} & i^*\Omega_{\mathbb{P}^r\times\mathbb{P}^n}^{p+2}(2\mathcal{X}) \xrightarrow{dF/F\wedge (-)}\cdots \\
		\cdots \to & i^*\Omega_{\mathbb{P}^r\times\mathbb{P}^n}^{n+r-1}(q\mathcal{X}) \xrightarrow{\pi_q}\Omega_{\mathcal{X}}^{n+r-1}(q\mathcal{X})\to 0
		\end{align*}
		where we put $i^*\Omega_{\mathbb{P}^r\times\mathbb{P}^n}^{p+1}(\mathcal{X})$ in degree zero. 
		This gives rise to the map 
		$$\delta: H^0(\mathcal{X},\Omega_{\mathcal{X}}^{n+r-1}(q\mathcal{X}))\to \mathbb{H}^q(\mathcal{X},\mathcal{C}(p))\cong H^q(\mathcal{X},\Omega_{\mathcal{X}}^p).$$ 
	\end{nota} 
	\begin{nota} 
		Note that a section $\xi$ of $\Omega_{\mathbb{P}^r\times\mathbb{P}^n}^{n+r}$ over $\mathbb{P}^r\times\mathbb{P}^n\setminus \mathcal{X}$ has a pole on $\mathcal{X}$ of order at most $a$ if and only if it is of the form $\xi = \frac{A\omega\wedge\bar{\omega}}{F^a}$ where $A\in k[Y_0,\cdots, Y_r,X_0,\cdots,X_n]$ is of bidegree $(a - (r+1), am - (n+1))$. This gives an isomorphism 
		\begin{align*} 
		\tilde{\psi_a}: k[Y_0,\cdots, Y_r,X_0,\cdots,X_n]^{a - (r+1), am - (n+1)} &\to H^0(\mathbb{P}^r\times\mathbb{P}^n, \Omega_{\mathbb{P}^r\times\mathbb{P}^n}^{n+r}(a\mathcal{X})) 
		\end{align*} 
		sending an element $A$ to $\frac{A\omega\wedge\bar{\omega}}{F^a}$. 
	\end{nota} 
	\begin{rem}
		Because $\Omega_{\mathbb{P}^r\times\mathbb{P}^n}^{n+r}(\log(\mathcal{X}))$ is isomorphic to $\Omega_{\mathbb{P}^r\times\mathbb{P}^n}^{n+r}(\mathcal{X})$, we can view the residue map in the exact sequence (\ref{residue exact sequence}) in this degree as a morphism $\Omega_{\mathbb{P}^r\times\mathbb{P}^n}^{n+r}(\mathcal{X}) \to i_*\Omega^{n+r-1}_{\mathcal{X}}$. 
	\end{rem}
	The purpose of this section is to prove the following statement, which is a generalization of \cite[Proposition 3.2]{LevineECHWHH} to products of projective spaces. The argument is taken from \cite{LevineECHWHH}, with some adaptations.   
	\begin{prop}\label{Levine Proposition 3.2}
		Suppose that $q\geq r$. The composition 
		\begin{align*} 
		k[Y_0,\cdots, Y_r, X_0,\cdots, X_n]^{(q+1)-(r+1), (q+1)m-(n+1)}&\xrightarrow{\tilde{\psi}_{q+1}} \\ H^0(\mathbb{P}^r\times\mathbb{P}^n, \Omega_{\mathbb{P}^r\times\mathbb{P}^n}^{n+r}((q+1)\mathcal{X})) 
		\xrightarrow{\text{res}} H^0(\mathcal{X}, \Omega_{\mathcal{X}}^{n+r-1}(q\mathcal{X})) &\xrightarrow{\delta} H^q(\mathcal{X}, \Omega^p_{\mathcal{X}})
		\end{align*} 
		descends to an isomorphism 
		$$\psi_q: J^{(q+1)-(r+1), (q+1)m-(n+1)}\to H^q(\mathcal{X}, \Omega^p_{\mathcal{X}})_{prim}.$$
	\end{prop}
	\begin{rem}\label{remark r=0}
		Let $X = V(F_0)\subset \mathbb{P}^n$ be a smooth hypersurface defined by a homogeneous polynomial $F_0\in k[X_0,\cdots, X_n]$ of degree $m$. Then we can form the hypersurface $\mathcal{X} = V(Y_0F_0)\subset\mathbb{P}^0\times\mathbb{P}^n$ which is isomorphic to $X$. The corresponding Jacobian ring is given by 
		$$J= k[Y_0,X_0\cdots, X_n]/\left(F_0,Y_0\frac{\partial F_0}{\partial X_0},\cdots, Y_0\frac{\partial F_0}{\partial X_n}\right).$$
		Consider the usual Jacobian ring 
		$$J_X = k[X_0,\cdots, X_n]/\left(\frac{\partial F_0}{\partial X_0},\cdots,\frac{\partial F_0}{\partial X_n}\right) $$
		then for $a,b\in\mathbb{Z}$ we have a natural map 
		$$J^{a,b}\to J_X^b, Y_0\mapsto 1, X_j\mapsto X_j.$$
		For fixed $a$, there is the section 
		$$g_a: J_X^{b}\to J^{a,b}, f\mapsto Y_0^af$$ which is an isomorphism. 
		Now let $p,q\in\mathbb{Z}_{\geq 0}$ be such that $p+q = n-1$. If we set $r=0$ (which is not possible with the assumptions we made, but we still do it for a moment) in Proposition \ref{Levine Proposition 3.2}, we find an isomorphism 
		$$\psi_q\circ g_q: J_X^{(q+1)m-n-1}\to H^q(X,\Omega^p_X)_{prim}.$$ We therefore find that the above statement is in accordance with \cite[Proposition 3.2]{LevineECHWHH}. 
	\end{rem}
	
	\subsection{An exact sequence relating $\delta$ and $\pi_q$} 
	In order to prove Proposition \ref{Levine Proposition 3.2}, we will need the following proposition, which is a generalization of \cite[Lemma 3.1(2)]{LevineECHWHH}. The proof is more or less the same, but included here for the reader's convenience.  
	\begin{prop}\label{Levine lemma 3.1(2)}
		The map $$\delta: H^0(\mathcal{X},\Omega_{\mathcal{X}}^{n+r-1}(q\mathcal{X}))\to \mathbb{H}^q(\mathcal{X},\mathcal{C}(p))\cong H^q(\mathcal{X},\Omega_{\mathcal{X}}^p)$$ gives rise to an exact sequence 
		\begin{equation}
		H^0(\mathcal{X}, i^*\Omega_{\mathbb{P}^r\times\mathbb{P}^n}^{n+r-1}(q\mathcal{X}))\xrightarrow{\pi_q} H^0(\mathcal{X}, \Omega_{\mathcal{X}}^{n+r-1}(q\mathcal{X}))\xrightarrow{\delta}H^q(\mathcal{X},\Omega_{\mathcal{X}}^p)
		\end{equation}
		where $\delta$ is surjective if $p\neq q$ and has image $H^p(\mathcal{X}, \Omega_{\mathcal{X}}^p)_{prim}$ in case that $p=q$. 
	\end{prop}
	In order to prove Proposition \ref{Levine lemma 3.1(2)}, we will use the hypercohomology spectral sequence 
	$$E_1^{a,b} = H^b(\mathcal{X}, \mathcal{C}(p)^a) \implies \mathbb{H}^{a+b}(\mathcal{X}, \mathcal{C}(p))\cong H^{a+b}(\mathcal{X}, \Omega^p_{\mathcal{X}}).$$
	We first prove two lemmas.  
	\begin{lem}\label{Levine lemma 3.1(2) helpful lemma}
		$E_1^{a,b} = 0$ for all $b>0$ and $a<q$, except for $a=0$ and $b=p$. 
	\end{lem}	
	\begin{proof}
		Let $a,b\in\mathbb{Z}$ be such that $b>0$ and $0\leq a<q$. Note that there is the standard exact sequence 
		$$0\to \mathcal{O}_{\mathbb{P}^r\times\mathbb{P}^n}(-\mathcal{X})\to\mathcal{O}_{\mathbb{P}^r\times\mathbb{P}^n}\to i_*\mathcal{O}_{\mathcal{X}}\to 0.$$
		Tensoring the above exact sequence with the sheaf $\Omega_{\mathbb{P}^r\times\mathbb{P}^n}^{p+a+1}((a+1)\mathcal{X})$ and noting that this sheaf is locally free gives the exact sequence 
		\begin{align} \label{exact sequence spectral sequence lemma}
		0\to \Omega_{\mathbb{P}^r\times\mathbb{P}^n}^{p+a+1}(a\mathcal{X})\to \Omega_{\mathbb{P}^r\times\mathbb{P}^n}^{p+a+1}((a+1)\mathcal{X}) 
		\to i_*i^*\Omega_{\mathbb{P}^r\times\mathbb{P}^n}^{p+a+1}((a+1)\mathcal{X})\to 0
		\end{align}
		where we used the projection formula to see that 
		\begin{align*}
		i_*\mathcal{O}_{\mathcal{X}}\otimes\Omega_{\mathbb{P}^r\times\mathbb{P}^n}^{p+a+1}((a+1)\mathcal{X}) &= i_*(\mathcal{O}_{\mathcal{X}}\otimes i^*\Omega_{\mathbb{P}^r\times\mathbb{P}^n}^{p+a+1}((a+1)\mathcal{X}))\\
		&= i_* i^*\Omega_{\mathbb{P}^r\times\mathbb{P}^n}^{p+a+1}((a+1)\mathcal{X}).
		\end{align*}
		Part of the long exact cohomology sequence of (\ref{exact sequence spectral sequence lemma}) is 
		\begin{align}\label{long exact sequence sslemma}
		&\cdots\to H^b(\mathbb{P}^r\times\mathbb{P}^n,  \Omega_{\mathbb{P}^r\times\mathbb{P}^n}^{p+a+1}(a\mathcal{X}))\to H^b(\mathbb{P}^r\times\mathbb{P}^n,  \Omega_{\mathbb{P}^r\times\mathbb{P}^n}^{p+a+1}((a+1)\mathcal{X}))\to \\
		&H^b(\mathbb{P}^r\times\mathbb{P}^n,  i_*i^*\Omega_{\mathbb{P}^r\times\mathbb{P}^n}^{p+a+1}((a+1)\mathcal{X}))\to H^{b+1}(\mathbb{P}^r\times\mathbb{P}^n,  \Omega_{\mathbb{P}^r\times\mathbb{P}^n}^{p+a+1}(a\mathcal{X}))\to \cdots  \nonumber 
		\end{align}
		We have by Proposition \ref{Botts theorem for general products} that 
		\begin{align*} 
		H^b(\mathbb{P}^r\times\mathbb{P}^n, &\Omega_{\mathbb{P}^r\times\mathbb{P}^n}^{p+a+1}((a+1)\mathcal{X})) = \\
		& \bigoplus_{i+j = p+a+1} \bigoplus_{k+l = b}H^k(\mathbb{P}^r, \Omega^i_{\mathbb{P}^r}(a+1))\otimes H^l(\mathbb{P}^n, \Omega^j_{\mathbb{P}^n}(m(a+1))). 
		\end{align*} 
		As $a+1 > 0$, by Theorem \ref{Botts theorem} one has that $H^k(\mathbb{P}^r, \Omega^i_{\mathbb{P}^r}(a+1)) =0$ provided that $k>0$ and $a+1 \geq i - r$ or $k=0$ and $a+1 \leq i$. We note that:
		\begin{itemize} 
			\item If $k>0$ and $a+1 < i-r$ we see in particular that $i-r>0$ and so $i>r$ implying that $\Omega_{\mathbb{P}^r}^i =0$ and so $H^k(\mathbb{P}^r, \Omega^i_{\mathbb{P}^r}(a+1)) =0$. 
			\item If $k=0$ and $a+1 > i$ then $l > 0$ as $b>0$. If $m(a+1) \geq j-n$ then $H^l(\mathbb{P}^n, \Omega^j_{\mathbb{P}^n}(m(a+1)))=0$. Otherwise, $m(a+1) < j-n$ implies that $j>n$ and so $H^l(\mathbb{P}^n, \Omega^j_{\mathbb{P}^n}(m(a+1)))=0$ in this case as well. 
		\end{itemize} 
		This proves that $H^b(\mathbb{P}^r\times\mathbb{P}^n, \Omega_{\mathbb{P}^r\times\mathbb{P}^n}^{p+a+1}((a+1)\mathcal{X})) =0$ for $b>0$ and $a\geq 0$. 
		Furthermore, using Proposition \ref{Botts theorem for general products} again we have that  
		\begin{align*} 
		H^{b+1}&(\mathbb{P}^r\times\mathbb{P}^n, \Omega_{\mathbb{P}^r\times\mathbb{P}^n}^{p+a+1}(a\mathcal{X})) = \\
		& \bigoplus_{i+j = p+a+1} \bigoplus_{k+l = b + 1}H^k(\mathbb{P}^r, \Omega^i_{\mathbb{P}^r}(a))\otimes H^l(\mathbb{P}^n, \Omega^j_{\mathbb{P}^n}(ma)).
		\end{align*} 
		We distinguish between two cases:
		\begin{itemize}
			\item If $a>0$ then $H^{b+1}(\mathbb{P}^r\times\mathbb{P}^n, \Omega_{\mathbb{P}^r\times\mathbb{P}^n}^{p+a+1}(a\mathcal{X})) = 0$ in a similar way as above. 
			\item If $a=0$ then using Theorem \ref{Botts theorem} we have that $H^k(\mathbb{P}^r, \Omega^i_{\mathbb{P}^r}) = 0$ provided that $k\neq i$ and $H^l(\mathbb{P}^n, \Omega^j_{\mathbb{P}^n}) = 0$ provided that $l\neq j$, so the only case where this is nonzero is if $k=i$ and $l=j$, i.e. if $p = b$. 
		\end{itemize}
		So if we assume that $p \neq b$ or $a>0$, then $H^{b+1}(\mathbb{P}^r\times\mathbb{P}^n, \Omega_{\mathbb{P}^r\times\mathbb{P}^n}^{p+a+1}(a\mathcal{X})) = 0$. From the sequence (\ref{long exact sequence sslemma}), we see that 
		\begin{align*} 
		E_1^{a,b} &= H^b(\mathcal{X}, i^*\Omega^{p+a+1}_{\mathbb{P}^r\times\mathbb{P}^n}((a+1)\mathcal{X}))= H^b(\mathbb{P}^r\times\mathbb{P}^n, i_*i^*\Omega^{p+a+1}_{\mathbb{P}^r\times\mathbb{P}^n}((a+1)\mathcal{X})) = 0
		\end{align*} which is precisely what we needed to show. 
	\end{proof}
	\begin{constr} 
		Now suppose that $a=0, b=p$ and $p>0$. \\
		We have seen in the proof of Lemma \ref{Levine lemma 3.1(2) helpful lemma} that 
		$$H^p(\mathbb{P}^r\times\mathbb{P}^n, \Omega_{\mathbb{P}^r\times\mathbb{P}^n}^{p+1}(\mathcal{X})) = H^{p+1}(\mathbb{P}^r\times\mathbb{P}^n, \Omega_{\mathbb{P}^r\times\mathbb{P}^n}^{p+1}(\mathcal{X})) = 0.$$ From the exact sequence (\ref{long exact sequence sslemma}) it follows that 
		\begin{equation} \label{identify E with Omega}
		E_1^{0,p} = H^{p}(\mathcal{X}, i^*\Omega_{\mathbb{P}^r\times\mathbb{P}^n}^{p+1}(\mathcal{X}))\cong H^{p+1}(\mathbb{P}^r\times\mathbb{P}^n, \Omega_{\mathbb{P}^r\times\mathbb{P}^n}^{p+1}).
		\end{equation} 
		Note that $E_\infty^{0,p}$ is nonzero. Indeed, $E_1^{a,b} = 0$ if $a+b = p+1$ and $b>0$, and one can compute that $E_1^{p+1,0} =0$. Also, we have that 
		$$H^{p+1}(\mathbb{P}^r\times\mathbb{P}^n, \Omega_{\mathbb{P}^r\times\mathbb{P}^n}^{p+1})\cong \bigoplus_{i+j =p+1}\bigoplus_{k+l = p+1} H^k(\mathbb{P}^r,\Omega_{\mathbb{P}^r}^i)\otimes H^l(\mathbb{P}^n,\Omega_{\mathbb{P}^n}^j)\cong k^{p+1}$$
		and so $E_1^{0,p}$ is nonzero by (\ref{identify E with Omega}). \\
		There is a surjection $H^p(\mathcal{X},\Omega_{\mathcal{X}}^p)\to E_\infty^{0,p}$. From this, we can define the map $$\beta: H^p(\mathcal{X},\Omega_{\mathcal{X}}^p)\to E_1^{0,p} = H^{p}(\mathcal{X}, i^*\Omega_{\mathbb{P}^r\times\mathbb{P}^n}^{p+1}(\mathcal{X})).$$
		Using (\ref{identify E with Omega}) again, we can view $\beta$ as a map $$H^{p}(\mathbb{P}^r\times\mathbb{P}^n, i_*\Omega_{\mathcal{X}}^p)\to H^{p+1}(\mathbb{P}^r\times\mathbb{P}^n, \Omega_{\mathbb{P}^r\times\mathbb{P}^n}^{p+1})$$ after identifying $H^p(\mathcal{X},\Omega_{\mathcal{X}}^p)$ with $H^{p}(\mathbb{P}^r\times\mathbb{P}^n, i_*\Omega_{\mathcal{X}}^p)$.  
	\end{constr} 
	\begin{lem}\label{beta is ilowerstar} 
		The morphism $\beta$ is precisely the coboundary map from the exact sequence~(\ref{residue exact sequence}). As a consequence, $\beta$ is surjective if $p\geq r$. 
	\end{lem} 
	\begin{proof} 
		Consider the diagram 
		\[
		\begin{tikzcd}
		0\arrow[r] & \Omega_{\mathbb{P}^r\times\mathbb{P}^n}^{p+1} \arrow[r] \arrow[d, equal] &  \Omega_{\mathbb{P}^r\times\mathbb{P}^n}^{p+1}(\log(\mathcal{X})) \arrow[r, "\text{res}"] \arrow[d] & i_*\Omega_{\mathcal{X}}^p\arrow[r] \arrow[d, "dF/F\wedge (-)"] & 0 \\
		0\arrow[r] & \Omega_{\mathbb{P}^r\times\mathbb{P}^n}^{p+1} \arrow[r] & \Omega_{\mathbb{P}^r\times\mathbb{P}^n}^{p+1}(\mathcal{X}) \arrow[r] & i_*i^*\Omega_{\mathbb{P}^r\times\mathbb{P}^n}^{p+1}(\mathcal{X}) \arrow[r] & 0
		\end{tikzcd} 
		\]
		with the top row coming from the exact sequence (\ref{residue exact sequence}) and the lower row coming from the exact sequence (\ref{exact sequence spectral sequence lemma}). \\
		This diagram commutes, because we know from the exact sequence (\ref{hyperplane exact sequence}) that the map $dF/F\wedge(-)$ is precisely the inclusion of $\mathcal{O}_{\mathcal{X}}(-\mathcal{X})$ into $i^*\Omega_{\mathbb{P}^r\times\mathbb{P}^n}$. By Lemma \ref{ilowerstar is boundary}, the coboundary map of the long exact sequence associated to the upper sequence is precisely $i_*$. Therefore, we have a commutative diagram 
		\[
		\begin{tikzcd} 
		H^p(\mathbb{P}^r\times\mathbb{P}^n, i_*\Omega_{\mathcal{X}}^p) \arrow[r, "i_*"] \arrow[d, "dF/F\wedge(-)"] \arrow[rd, "\beta"]& H^{p+1}(\mathbb{P}^r\times\mathbb{P}^n, \Omega_{\mathbb{P}^r\times\mathbb{P}^n}^{p+1}) \arrow[d, equal] \\
		H^p(\mathbb{P}^r\times\mathbb{P}^n, i_*i^*\Omega_{\mathbb{P}^r\times\mathbb{P}^n}^{p+1}(\mathcal{X})) \arrow[r, "\cong "] & H^{p+1}(\mathbb{P}^r\times\mathbb{P}^n, \Omega_{\mathbb{P}^r\times\mathbb{P}^n}^{p+1})
		\end{tikzcd} 
		\]
		Here we note that we constructed the spectral sequence from the exact sequence~(\ref{hyperplane exact sequence}), so that $\beta$ is by construction the map induced from $dF/F\wedge(-)$, composed with the isomorphism coming from the coboundary map in the sequence (\ref{exact sequence spectral sequence lemma}). It follows that $\beta = i_*$. We note that $i_*$ is surjective if $p\geq r$ by Lemma \ref{ilowerstar surjective lemma}, which gives the last part of the statement.  
	\end{proof} 
	\begin{lem}\label{outgoing differentials of 0,p}
		All differentials going into or out of $E_s^{0,p}$ are zero for $p\geq r$. 
	\end{lem}
	\begin{proof}
		Note that all incoming differentials to the terms $E_s^{0,p}$ are zero by reason of degree, and we have that $E_s^{0,p} = \ker(d_{s-1}^{0,p}) \subset E_{s-1}^{0,p}$ for all $s\geq 2$. Now the fact that the edge map $\beta$ is surjective for $p\geq r$ by Lemma \ref{beta is ilowerstar} shows that all outgoing differentials are zero. 
	\end{proof}
	\begin{proof}[Proof of Proposition \ref{Levine lemma 3.1(2)}] 
		Note that if $q=0$, then we have the exact sequence 
		$$H^0(\mathcal{X}, i^*\Omega_{\mathbb{P}^r\times\mathbb{P}^n}^{n+r-1})\to H^0(\mathcal{X}, \Omega_{\mathcal{X}}^{n+r-1})\to  H^0(\mathcal{X}, \Omega_{\mathcal{X}}^{n+r-1})$$ where the last map is the identity and 
		\begin{align*} 
		H^0(\mathcal{X}, i^*\Omega_{\mathbb{P}^r\times\mathbb{P}^n}^{n+r-1}) &= H^0(\mathbb{P}^r\times\mathbb{P}^n, i_*i^*\Omega_{\mathbb{P}^r\times\mathbb{P}^n}^{n+r-1}).
		\end{align*} 
		Note that \begin{align*} 
		H^0(\mathbb{P}^r\times\mathbb{P}^n, \Omega_{\mathbb{P}^r\times\mathbb{P}^n}^{n+r-1})&\cong H^0(\mathbb{P}^r,\Omega_{\mathbb{P}^r}^r)\otimes H^0(\mathbb{P}^n,\Omega_{\mathbb{P}^n}^{n-1})\\
		&\quad \oplus H^0(\mathbb{P}^r,\Omega_{\mathbb{P}^r}^{r-1})\otimes H^0(\mathbb{P}^n,\Omega_{\mathbb{P}^n}^{n})\\
		&= 0\end{align*} 
		using Proposition \ref{Botts theorem for general products} and Theorem \ref{Botts theorem} again. Similarly, we have that 
		\begin{align*} 
		H^1(\mathbb{P}^r\times\mathbb{P}^n, \Omega_{\mathbb{P}^r\times\mathbb{P}^n}^{n+r-1}(-\mathcal{X}))&\cong H^0(\mathbb{P}^r,\Omega_{\mathbb{P}^r}^r(-1))\otimes H^1(\mathbb{P}^n,\Omega_{\mathbb{P}^n}^{n-1}(-m))\\
		&\quad \oplus H^1(\mathbb{P}^r,\Omega_{\mathbb{P}^r}^{r-1}(-1))\otimes H^0(\mathbb{P}^n,\Omega_{\mathbb{P}^n}^{n}(-m))\\
		&\quad \oplus H^0(\mathbb{P}^r,\Omega_{\mathbb{P}^r}^{r-1}(-1))\otimes H^1(\mathbb{P}^n,\Omega_{\mathbb{P}^n}^{n}(-m))\\
		&\quad \oplus H^1(\mathbb{P}^r,\Omega_{\mathbb{P}^r}^{r}(-1))\otimes H^0(\mathbb{P}^n,\Omega_{\mathbb{P}^n}^{n-1}(-m))\\
		&= 0\end{align*} 
		so using the long exact sequence (\ref{long exact sequence sslemma}) we find that $H^0(\mathcal{X}, i^*\Omega_{\mathbb{P}^r\times\mathbb{P}^n}^{n+r-1})=0$. Therefore, we assume that $q>0$ from now on. \\
		First assume that $p\neq q$. We note that the contributions to $H^q(\mathcal{X}, \Omega_{\mathcal{X}}^p)$ come from all $E_1^{a,b}$ satisfying $a+b = q$. These are all zero except possibly for $E_1^{q,0}$, by Lemma \ref{Levine lemma 3.1(2) helpful lemma}. By Lemma  \ref{Levine lemma 3.1(2) helpful lemma}, we have that $E_1^{a,b}=0$ for $a+b=q-1$ and $0<a\leq q-1$, except possibly for $(a,b)=(q-1,0)$ and, if $p=q-1$, also $(a,b)=(0,q-1)$. Thus, the only possible non-zero incoming differentials to $E_*^{q,0}$ are $d_1^{q-1,0}$ and, in case $p=q-1$, also $d_q^{0,q-1}$. Also, there are no outgoing differentials out of $E_*^{q,0}$ by reason of degree. \\
		If $p = q-1$ and $q\neq 1$, then as $p+q = n+r-1$, we have that $2p = n+r-2\geq 2r$ and so $p\geq r$. Therefore by Lemma \ref{outgoing differentials of 0,p}, all outgoing differentials of $E_q^{0,q-1}$ are zero. If $q=1$, then $E_1^{0,q-1} = E_1^{0,0}$ and so it coincides with $E_1^{q-1,0}$.  Either way, we therefore have that the only possibly nonzero incoming differential comes from $E_1^{q-1,0}$, and so $E_\infty^{q,0}$ is equal to $E_2^{q,0}$, i.e. $E_1^{q,0}$ modulo the image of $E_1^{q-1,0}$. Since all the differentials leaving $E_*^{q,0}$ are zero, we have the edge homomorphism 
		$$E_1^{q,0} = H^0(\mathcal{X}, \Omega_{\mathcal{X}}^{n+r-1}(q\mathcal{X})) \to H^q(\mathcal{X}, \Omega_{\mathcal{X}}^p)$$ which is exactly the map $\delta$. Since $E_2^{q,0}=E_\infty^{q,0}$, this gives an exact sequence 
		$$0\to H^0(\mathcal{X}, i^*\Omega_{\mathbb{P}^r\times\mathbb{P}^n}^{n+r-1}(q\mathcal{X})) \to H^0(\mathcal{X}, \Omega_{\mathcal{X}}^{n+r-1}(q\mathcal{X})) \xrightarrow{\delta} H^q(\mathcal{X}, \Omega_{\mathcal{X}}^p).$$ Note that the map $\delta$ is surjective in this case since only $E_*^{q,0}$ can contribute. This completes the proof in case $p\neq q$.  \\
		If $p = q$ then the possibly nonzero term $E_1^{0, p}$ also contributes to $H^p(\mathcal{X}, \Omega_{\mathcal{X}}^p)$.  There are no incoming differentials to $E_*^{0,p}$ and since $p=q=n+r-1-p$, we have $p\ge r$, so there no outgoing differentials by Lemma \ref{outgoing differentials of 0,p}, and we have $E_1^{0,p}=E_\infty^{0,p}$. Thus, the edge homomorphism $$\beta: H^p(\mathcal{X},\Omega_{\mathcal{X}}^p)\to E_1^{0,p} = H^{p}(\mathcal{X}, i^*\Omega_{\mathbb{P}^r\times\mathbb{P}^n}^{p+1})$$ is surjective, and we find an exact sequence 
		\begin{align*} 
		0&\to H^0(\mathcal{X}, i^*\Omega_{\mathbb{P}^r\times\mathbb{P}^n}^{n+r-1}(q\mathcal{X})) \to H^0(\mathcal{X}, \Omega_{\mathcal{X}}^{n+r-1}(q\mathcal{X})) \xrightarrow{\delta} H^p(\mathcal{X}, \Omega_{\mathcal{X}}^p) \\
		&\xrightarrow{\beta} H^{p}(\mathcal{X}, i^*\Omega_{\mathbb{P}^r\times\mathbb{P}^n}^{p+1})\to 0
		\end{align*} 
		This means that the image of $\delta$ is equal to the kernel of $\beta$, which by Lemma~\ref{beta is ilowerstar} is the kernel of~$i_*$, which is by definition the primitive cohomology.
	\end{proof}
	
	\subsection{Proof of Proposition \ref{Levine Proposition 3.2}}
	\begin{proof}[Proof of Proposition \ref{Levine Proposition 3.2}]
		We first note that $H^1(\mathbb{P}^r\times\mathbb{P}^n, \Omega_{\mathbb{P}^r\times\mathbb{P}^n}^{n+r}(q\mathcal{X})) =0$. Indeed, it follows from Proposition \ref{Botts theorem for general products} that 
		\begin{align*}
		H^1(\mathbb{P}^r\times\mathbb{P}^n, \Omega_{\mathbb{P}^r\times\mathbb{P}^n}^{n+r}(q\mathcal{X})) = \bigoplus_{i+j = n+r} & H^0(\mathbb{P}^r, \Omega_{\mathbb{P}^r}^i(q))\otimes H^1(\mathbb{P}^n, \Omega_{\mathbb{P}^n}^j(mq)) \\
		&\oplus  H^1(\mathbb{P}^r, \Omega_{\mathbb{P}^r}^i(q))\otimes H^0(\mathbb{P}^n, \Omega_{\mathbb{P}^n}^j(mq)) 
		\end{align*}
		Using Theorem \ref{Botts theorem} we observe that: 
		\begin{itemize}
			\item If $q>0$ then $H^1(\mathbb{P}^n, \Omega_{\mathbb{P}^n}^j(mq)) =0$ provided that $mq \geq j-n$. But if $mq < j-n$ we have that $j>n$ and so $H^1(\mathbb{P}^n, \Omega_{\mathbb{P}^n}^j(mq))$ is always equal to zero. For similar reasons, $ H^1(\mathbb{P}^r, \Omega_{\mathbb{P}^r}^i(q))$ is always equal to zero and so $	H^1(\mathbb{P}^r\times\mathbb{P}^n, \Omega_{\mathbb{P}^r\times\mathbb{P}^n}^{n+r}(q\mathcal{X})) =0$ in this case. 
			\item If $q=0$ then $H^1(\mathbb{P}^n, \Omega_{\mathbb{P}^n}^j) =0$ provided that $j\neq 1$. However, in case $j=1$ we have that $i = n+r-1 > 0$ as $n\geq 2$. Now from Theorem \ref{Botts theorem} it follows that $H^0(\mathbb{P}^r,\Omega_{\mathbb{P}^r}^i) = 0$. Similarly the second factor always vanishes. 
		\end{itemize}
		This verifies the claim. \\
		This means that the exact sequence (\ref{exact sequence spectral sequence lemma}) gives rise to an exact sequence 
		\begin{align*} 
		0&\to H^0(\mathbb{P}^r\times\mathbb{P}^n, \Omega_{\mathbb{P}^r\times\mathbb{P}^n}^{n+r}(q\mathcal{X}))\to H^0(\mathbb{P}^r\times\mathbb{P}^n, \Omega_{\mathbb{P}^r\times\mathbb{P}^n}^{n+r}((q+1)\mathcal{X})) \\
		&\quad \to H^0(\mathcal{X}, i^*\Omega_{\mathbb{P}^r\times\mathbb{P}^n}^{n+r}((q+1)\mathcal{X}))\to 0.
		\end{align*} 
		We now consider $\text{res}:H^0(\mathbb{P}^r\times\mathbb{P}^n, \Omega_{\mathbb{P}^r\times\mathbb{P}^n}^{n+r}((q+1)\mathcal{X}))\to H^0(\mathcal{X},\Omega_{\mathcal{X}}^{n+r-1}(q\mathcal{X}))$. Noting that the  kernel of $\text{res}$ is $H^0(\mathbb{P}^r\times\mathbb{P}^n, \Omega_{\mathbb{P}^r\times\mathbb{P}^n}^{n+r}(q\mathcal{X}))$, we see that $\text{res}$ descends to a map
		$$\widetilde{\text{res}}: H^0(\mathcal{X}, i^*\Omega_{\mathbb{P}^r\times\mathbb{P}^n}^{n+r}((q+1)\mathcal{X}))\to H^0(\mathcal{X},\Omega_{\mathcal{X}}^{n+r-1}(q\mathcal{X}))$$ with the same image as $\text{res}$. Furthermore, we observe that the map $\tilde{\psi}_{q+1}$ now gives rise to an isomorphism 
		\begin{equation} \label{isomorphism tildepsiq}
		k[Y_0,\cdots, Y_r, X_0,\cdots, X_n]/(F)^{q-r, (q+1)m-(n+1)} \xrightarrow{\bar{\psi}_{q+1}} H^0(\mathcal{X}, i^*\Omega_{\mathbb{P}^r\times\mathbb{P}^n}^{n+r}((q+1)\mathcal{X})).
		\end{equation} 
		We now consider the following commutative diagram: 
		\begin{equation} \label{diagram section 3}
		\begin{tikzcd}
		H^0(\mathcal{X}, i^*\Omega_{\mathbb{P}^r\times\mathbb{P}^n}^{n+r}((q+1)\mathcal{X})) \arrow[rd, "\widetilde{\text{res}}"] & k[Y_0,\cdots, Y_r, X_0,\cdots, X_n]/(F)^\alpha \arrow[l, swap, "\bar{\psi}_{q+1}"] \arrow[d, "f"] \\
		H^0(\mathcal{X}, i^*\Omega_{\mathbb{P}^r\times\mathbb{P}^n}^{n+r-1}(q\mathcal{X})) \arrow[r, "\pi_q"] \arrow[u, "dF/F\wedge(-)"] & H^0(\mathcal{X}, \Omega_{\mathcal{X}}^{n+r-1}(q\mathcal{X})) 
		\end{tikzcd} 
		\end{equation} 
		where $\alpha = ((q+1)-(r+1), (q+1)m-(n+1))$ is the bidegree from the statement and $f = \widetilde{\text{res}}\circ\bar{\psi}_{q+1}$. For the rest of the proof, we will show that the image of the Jacobian ideal under $f$ is the same as the image of $\pi_q$, i.e. the kernel of $\delta$. This shows that the Jacobian ideal is precisely the kernel of the composition $\delta\circ f$. Note that the image of $\delta\circ f$ is the image of $\delta$ as $\text{res}$ is surjective and $\bar{\psi}_{q+1}$ is an isomorphism. The image of $\delta$ however, is precisely the primitive cohomology. Therefore, noting that $F$ is in the Jacobian ideal, we find an isomorphism 
		$$k[Y_0,\cdots, Y_r,X_0,\cdots, X_n]/(F_0,\cdots, F_r,\bar{F}_0,\cdots, \bar{F}_n)^\alpha\to H^p(\mathcal{X},\Omega_{\mathcal{X}}^p)_{prim}.$$
		This will complete the proof. \\
		We start by noting that the image of ideal $\left(F_0,\cdots, F_r, \bar{F}_0,\cdots, \bar{F}_n\right)/F$ under the isomorphism (\ref{isomorphism tildepsiq}) is generated by elements of the form $F_i\omega\wedge \bar{\omega}$ and $\bar{F}_j\omega\wedge \bar{\omega}$. Following \cite{TerasomaIVH}, we note that the sheaf $\Omega_{\mathbb{P}^r}^r(r+1)\otimes\Omega_{\mathbb{P}^n}^{n-1}(n)$ has global sections generated by the sections $\omega\wedge\tau_i$ for $i\in\{0,\cdots, n\}$ where 
		\begin{align*} 
		\tau_i &= \sum_{j<i}(-1)^jX_jdX^{j,i} + \sum_{i<j}(-1)^{j+1}X_jdX^{i,j}. 
		\end{align*} 
		So a section of this sheaf with a pole along $\mathcal{X}$ of order at most $q$ is of the form $\sum_{i=0}^n\frac{B_i\omega\wedge\tau_i}{F^q}$ where $B_i$ has bidegree $(q-(r+1), qm - n)$. We note that $dF = \sum_{i=0}^r F_idY_i + \sum_{j=0}^n \bar{F}_j dX_j$ and we also note that for $i,j\in\{0,\cdots, n\}$, we have that
		$$dX_j\wedge \tau_i= 
		\begin{cases} 
		X_j dX^i  &\text{ if } i\neq j\\
		\sum_{j<  i} (-1)^{j+i-1} X_jdX^j + \sum_{j >  i} (-1)^{j+i + 1} X_jdX^j & \text{ otherwise}
		\end{cases} $$
		So 
		\begin{align*} 
		\left(\sum_{j=0}^n \bar{F}_j dX_j \right)\wedge \tau_i &= \sum_{j\neq i}\bar{F}_jX_j dX^i + (-1)^{i-1}\bar{F}_i\sum_{j\neq i} (-1)^j X_jdX^j \\
		&= (-1)^{i-1}\bar{F}_i\bar{\omega} + mFdX^i.
		\end{align*} 
		It follows that modulo $F$, we have that
		\begin{align*} 
		\frac{dF}{F}\wedge \sum_{i=0}^n\frac{B_i\omega\wedge\tau_i}{F^q} &= \sum_{i=0}^n (-1)^{i-1}\frac{B_i\bar{F}_i\omega\wedge\bar{\omega}}{F^{q+1}}
		\end{align*} 
		We can do a similar thing for a section of $\Omega_{\mathbb{P}^r}^{r-1}(r)\otimes\Omega_{\mathbb{P}^n}^{n}(n+1)$.\\
		Using \cite[Exercise II.8.3]{HartshorneAG} we see that $\Omega_{\mathbb{P}^r\times\mathbb{P}^n}^{n+r-1} = \bigoplus_{i+j = n+r-1}\pi_r^*\Omega_{\mathbb{P}^r}^i\otimes \pi_n^*\Omega_{\mathbb{P}^n}^j$. As the map $dF/F\wedge(-)$ increases the bidegrees of the corresponding elements of $J$ by $(1, m-1)$, the only forms that will end up in the bidegree $\alpha$ when applying $dF/F\wedge(-)$ are either coming from $\Omega_{\mathbb{P}^r}^r\otimes\Omega_{\mathbb{P}^n}^{n-1}$ or from $\Omega_{\mathbb{P}^r}^{r-1}\otimes\Omega_{\mathbb{P}^n}^{n}$. Using the fact that $H^0(\mathbb{P}^r\times\mathbb{P}^n, \Omega_{\mathbb{P}^r\times\mathbb{P}^n}^{n+r}((q+1)\mathcal{X}))$ surjects onto $H^0(\mathcal{X}, i^*\Omega_{\mathbb{P}^r\times\mathbb{P}^n}^{n+r}((q+1)\mathcal{X}))$ the above computation shows that the image of the map $dF/F\wedge(-)$ in the diagram (\ref{diagram section 3}) is the same as the image of $\tilde{\psi}_{q+1}$, and so the images remain the same after applying $\widetilde{\text{res}}$. This shows that $f$ and $\pi_q$ indeed have the same image, which completes the proof.  
	\end{proof}
	
	\section{Comparing the two products}\label{section two products}
	The next goal is to compare the cup product on cohomology to the ring multiplication on the Jacobian ring. We keep all notation from the previous two sections. 
	\begin{nota} 
		Let  $\mathcal{U} = \{U_0,\cdots, U_r, \bar{U}_0,\cdots, \bar{U}_n\}$ be the open cover of $\mathbb{P}^r\times\mathbb{P}^n$ where $U_i = \{F_i\neq 0\}$ for  $i\in\{0,\cdots, r\}$ and $\bar{U}_j = \{\bar{F}_j\neq 0\}$ for $j\in\{0,\cdots, n\}$. Note that this is an open cover of $\mathbb{P}^r\times\mathbb{P}^n$ because $\mathcal{X}$ is smooth. We have that $\mathcal{U}$ restricts to a cover of $\mathcal{X}$, and we will use the same notation for both. Note that the $U_i$ are not affine and that the $F_i$'s and $\bar{F}_j$'s do not have the same bidegrees. We define an order on these open subsets as follows: $$U_0<U_1<\cdots <U_r<\bar{U}_0<\cdots < \bar{U}_n.$$ 
	\end{nota} 
	\begin{nota}
		For $i\in\{0,\cdots, r\}$, we let $K_i$ be inner multiplication with $\frac{\partial }{\partial Y_i}$, i.e. for $i\in\{0,\cdots, r\}$ we have that 
		\begin{align*} 
		&K_i(dY_{i_1}\wedge\cdots\wedge dY_{i_l}\wedge dX_{j_{1}}\wedge\cdots\wedge dX_{j_k})\\
		&= \begin{cases} 
		(-1)^{a-1}dY_{i_1}\wedge\cdots d\hat{Y}_{i_a}\cdots\wedge dY_{i_l}\wedge dX_{j_{1}}\wedge\cdots\wedge dX_{j_k} &\text{ if } i=i_a \\
		0 &\text{ if } i\notin\{i_1,\cdots, i_l\}
		\end{cases} 
		\end{align*} 
		Similarly, for $j\in\{0,\cdots, n\}$, we let $\tilde{K}_j$ be inner multiplication with $\frac{\partial }{\partial X_j}$, i.e. for $j\in\{0,\cdots, n\}$, we have that 	
		\begin{align*} 
		&\bar{K}_j(dY_{i_1}\wedge\cdots\wedge dY_{i_l}\wedge dX_{j_{1}}\wedge\cdots\wedge dX_{j_k})\\
		&= \begin{cases} 
		(-1)^{l+b-1}dY_{i_1}\wedge\cdots \wedge dY_{i_l}\wedge dX_{j_{1}}\wedge\cdots d\hat{X}_{j_b}\cdots \wedge dX_{j_k} &\text{ if } j=j_b \\
		0 &\text{ if } j\notin\{j_1,\cdots, j_k\}
		\end{cases} 
		\end{align*} 
		We have that  $K_i(\alpha\wedge\beta) = K_i(\alpha)\wedge\beta + (-1)^k\alpha\wedge K_i(\beta)$, if $\alpha$ is a $k$-form, and a similar formula holds for the $\bar{K}_j$'s. 
	\end{nota} 
	\begin{nota} 
		For subsets $I = {i_1,\cdots, i_l}\subset \{0,\cdots, r\}$ and $J\subset \{0,\cdots, n\}$, we set:
		\begin{itemize} 
			\item $U_{I,J} = \bigcap_{i\in I}U_i\cap \bigcap_{j\in J}\bar{U}_j$. 
			\item $\Omega_I = (\prod_{i\in I}K_i)(\omega) = (K_{i_l}\circ K_{i_{l-1}}\circ\cdots\circ K_{i_1})(\omega)$ and $\bar{\Omega}_J = (\prod_{j\in J}\bar{K}_j)(\bar{\omega})$.
			\item $F_I = \prod_{i\in I}F_i$ and $\bar{F}_J = \prod_{j\in J}\bar{F}_j$.  
		\end{itemize} 
	\end{nota}
	\begin{nota} 
		For a sheaf $\mathcal{F}$ on $\mathcal{X}$ or on $\mathbb{P}^r\times\mathbb{P}^n$, we have the group $C^i(\mathcal{U}, \mathcal{F})$ consisting of all families $\{s_{I,J}\in \mathcal{F}(U_{I,J})\}_{I,J}$ where $\#I + \#J = i+1$. These form the \v{C}ech complex 
		$$0\to C^0(\mathcal{U}, \mathcal{F})\to \cdots \to C^{n+r+1}(\mathcal{U},  \mathcal{F})\to 0$$ with differentials $\delta:C^i(\mathcal{U}, \mathcal{F})\to C^{i+1}(\mathcal{U}, \mathcal{F})$ given by $$\delta(\{s_{I,J}\}_{I,J}) = \Big\{\sum_{k=0}^{i+1}(-1)^ks_{(I,J)\setminus \{(I,J)_k\}}|_{U_{I,J}}\Big\}_{I,J}.$$ 
		Here, $(I,J)_k$ denotes the $k$'th element of the ordered set $(I,J)$. The cohomology groups of the above complex are denoted by $\check{H}^a(\mathcal{U}, \mathcal{F})$. Note that there are natural maps $\check{H}^a(\mathcal{U}, \mathcal{F})\to H^a(\mathcal{X},\mathcal{F})$ or $\check{H}^a(\mathcal{U}, \mathcal{F})\to H^a(\mathbb{P}^r\times\mathbb{P}^n,\mathcal{F})$, by \cite[Lemma II.4.4]{HartshorneAG}. 
	\end{nota} 
	\begin{nota} \label{construction of omega(I,J)}
		This notation is taken from \cite[page 222]{TerasomaIVH} with a small adaptation, see the remark below. Fix $p,q\in\mathbb{Z}_{\geq 0}$ such that $p+q = n+r-1$. Consider the bidegree 
		$$\rho = (n-r-1, (n+r+1)m - 2(n+1)).$$
		For subsets $I\subset\{0,\cdots, r\}$ and $J\subset\{0,\cdots, n\}$ such that $\#I + \#J = n+r$, we define an element $\Omega(I,J)\in H^0(\mathcal{X}\cap U_{I,J}, i^*\Omega_{\mathbb{P}^r\times\mathbb{P}^n}^{n+r-1}(-\rho))$ as follows:
		\begin{itemize}
			\item If we are not in the situation where $\#I = r$ and $\#J = n$, then $\Omega(I,J) = 0$. 
			\item If $I = \{i_0,\cdots, i_{r-1}\}\subset\{0,\cdots, r\}$ and $J = \{j_0,\cdots, j_{n-1}\}\subset\{0,\cdots, n\}$, then: 
			\begin{itemize} 
				\item If $q\leq r-1$, we write 
				$$I' = \begin{cases}
				\{i_0,\cdots, i_{q-1}\} &\text{ if } q>0\\
				\emptyset &\text{ if } q=0
				\end{cases}$$ and 
				$$I''=\begin{cases} 
				\{i_{q+1},\cdots, i_{r-1}\} &\text{ if } q<r-1\\
				\emptyset &\text{ if } q=r-1
				\end{cases} $$ 
				We define 
				$$\Omega(I,J) = \frac{(-1)^{nq}\Omega_{I',i_{q}}\wedge\bar{\omega}\wedge\Omega_{i_{q},I''}\wedge\bar{\Omega}_J}{F_{I'}F_{i_{q}}^2F_{I''}\bar{F}_J}.$$
				\item If $q\geq r$, then we write 
				$$J' = \begin{cases}
				\{j_0,\cdots, j_{q-r-1}\} &\text{ if } q>r\\
				\emptyset &\text{ if } q=r
				\end{cases}$$ and 
				$$J''= 
				\begin{cases} 
				\{j_{q-r+1},\cdots, j_{n-1}\}  &\text{ if } q < n+r-1\\
				\emptyset &\text{ if } q=n+r-1
				\end{cases}$$
				We define 
				$$\Omega(I,J) = \frac{(-1)^{r(n+q+r)}\Omega_{I}\wedge\bar{\Omega}_{J', j_{q-r}}\wedge\omega\wedge\bar{\Omega}_{j_{q-r}, J''}}{F_{I}\bar{F}_{J'}\bar{F}_{j_{q-r}}^2\bar{F}_{J''}}.$$
			\end{itemize}
		\end{itemize} 
	\end{nota} 
	\begin{rem}\label{Omega(I,J) = 0 with wrong cardinalities}
		Note that $\Omega(I,J)$ can also be defined without distinguishing the case $\#I = r$ and $\#J = n$, which is the definition used in \cite{TerasomaIVH}. We then still have that $\Omega(I,J) = 0$ unless $\#I = r$ and $\#J = n$. Namely, $\Omega_{\{0,\cdots, r\}} = 0$ and $\bar{\Omega}_{\{0,\cdots, n\}} = 0$, which implies that $\Omega(I,J) = 0$ for $\#I\geq r+1$ or $\#J\geq n+1$. As $\#I + \#J = n+r$, this means that $\Omega(I,J) = 0$ unless $\#I = r$ and $\#J = n$.  
	\end{rem} 
	In subsection \ref{subsection proof cup product}, we will prove the following statement, which is a generalization of \cite[Proposition 2.8]{TerasomaIVH} to other fields than $\mathbb{C}$. 
	\begin{prop}\label{Proposition comparison products}
		For $A\in J^{q-r, (q+1)m - (n+1)}$ and $B\in J^{p-r, (p+1)m - (n+1)}$, write $$\omega_A = \psi_q(A)\in H^q(\mathcal{X}, \Omega_{\mathcal{X}}^p)\text{ and }\omega_B = \psi_p(B)\in H^p(\mathcal{X}, \Omega_{\mathcal{X}}^q).$$ Then the cup product $\omega_A\cup\omega_B\in H^{n+r-1}(\mathcal{X}, \Omega_{\mathcal{X}}^{n+r-1})$ is represented by the \v{C}ech cochain 
		$$\{\pi(AB\Omega(I,J))\}_{I,J}\in C^{n+r-1}(\mathcal{U}, \Omega_{\mathcal{X}}^{n+r-1}).$$
	\end{prop} 
	The argument is taken almost directly from \cite{TerasomaIVH}, combined with some elements of the arguments in \cite{LevineECHWHH}. 
	\begin{rem}[Remark \ref{remark r=0} continued]\label{remark r=0 products}
		Let $V(F_0)\subset \mathbb{P}^n$ be a hypersurface defined by a homogeneous polynomial $F_0\in k[X_0,\cdots, X_n]$ of degree $m$. If we set $r=0$ then for $J = \{j_0,\cdots, j_{n-1}\} \subset \{0,\cdots, n\}$ and $I = \emptyset$, we have $$\Omega(I,J) = \frac{\bar{\Omega}_{j_0,\cdots, j_q}\wedge\bar{\Omega}_{j_q,\cdots, j_{n-1}}}{Y_0^{n}\frac{\partial F_0}{\partial X_q}\prod_{j=0}^n \frac{\partial F_0}{\partial X_j}}.$$ One can check that Proposition \ref{Proposition comparison products} then becomes \cite[Proposition 3.6(1)]{LevineECHWHH}. 
	\end{rem} 
	We can compute the element $i_*(\omega_A\cup\omega_B)$ as follows, generalizing \cite[Proposition 3.7(2)]{LevineECHWHH} to a product of projective spaces. 
	\begin{prop}\label{representation of image cech}
		We have that $i_*(\omega_A\cup\omega_B)\in H^{n+r}( \mathbb{P}^r\times\mathbb{P}^n, \Omega_{\mathbb{P}^r\times\mathbb{P}^n}^{n+r})$ is represented by the cochain in $\{s_0,\cdots, s_r,\bar{s}_0,\cdots, \bar{s}_n\}\in C^{n+r}(\mathcal{U}, \Omega_{\mathbb{P}^r\times\mathbb{P}^n}^{n+r})$ given by 
		$$s_v = \frac{(-1)^{v+r+1} mABY_vF_v\omega\wedge\bar{\omega}}{ \prod_{i=0}^rF_i\prod_{j=0}^n\bar{F}_{j}}$$ for $v\in\{0,\cdots, r\}$ corresponding to the intersection of all opens except for $U_v$ and  
		$$\bar{s}_w = \frac{(-1)^{w+1}ABX_w\bar{F}_w\omega\wedge\bar{\omega}}{ \prod_{i=0}^rF_i\prod_{j=0}^n\bar{F}_{j}}$$ for $w\in\{0,\cdots, n\}$ corresponding to the intersection of all opens except for $\bar{U}_w$. 
	\end{prop} 
	It will have the following consequence, which can be viewed as a generalization of \cite[Corollary 2.9]{TerasomaIVH} to other fields than $\mathbb{C}$. 
	\begin{cor}\label{corollary definition map Jrho to cohomology}
		Consider the morphism 
		$$\tilde{\phi}: k[Y_0,\cdots, Y_r, X_0,\cdots, X_n]^\rho\to C^{n+r}(\mathcal{U}, \Omega_{\mathbb{P}^r\times\mathbb{P}^n}^{n+r}), D\mapsto \{s_0,\cdots, s_r,\bar{s}_0,\cdots, \bar{s}_n\}$$ 
		where 
		$$s_v = \frac{(-1)^{v+r+1}mDY_vF_v\omega\wedge\bar{\omega}}{\prod_{i=0}^rF_i\prod_{j=0}^n\bar{F}_j} \text{ and } \bar{s}_w =  \frac{(-1)^{w+1}DX_w\bar{F}_w\omega\wedge\bar{\omega}}{\prod_{i=0}^rF_i\prod_{j=0}^n\bar{F}_j}$$
		for $v\in\{0,\cdots, r\}$ and $w\in\{0,\cdots, n\}$. This gives rise to a surjective morphism $\phi: J^\rho\to H^{n+r}(\mathbb{P}^r\times\mathbb{P}^n, \Omega_{\mathbb{P}^r\times\mathbb{P}^n}^{n+r})\cong k$, such that the diagram 
		\[
		\begin{tikzcd} 
		H^q(\mathcal{X}, \Omega_{\mathcal{X}}^p)_{prim}\otimes H^p(\mathcal{X}, \Omega_{\mathcal{X}}^q)_{prim} \arrow[r, "i_*\circ \cup"] & H^{n+r}(\mathbb{P}^r\times\mathbb{P}^n, \Omega_{\mathbb{P}^r\times\mathbb{P}^n}^{n+r}) \\
		J^{q-r,(q+1)m - (n+1)} \otimes J^{p-r,(p+1)m - (n+1)}  \arrow[r] \arrow[u, "\psi_p\otimes\psi_q"]  & J^\rho \arrow[u, "\phi"] 
		\end{tikzcd} 
		\]
		commutes. 
	\end{cor}  
	Then in subsection \ref{section onedimensionality of Jrho}, we use a slight variation on an argument from \cite{KonnoVTPCI} to find the following. 
	\begin{cor}\label{Corollary Jrho ondimensional}
		The map $\phi$ is an isomorphism, except possibly when $n$ is odd, $r=1$ and $m=2$. 
	\end{cor}
	\begin{rem} 
		Note that if $n$ is odd and $r=1$, we have that $n+r-1 = n$ is odd, so that $\mathcal{X}$ has odd dimension. By the Motivic Gauss Bonnet Theorem, see Theorem \ref{Motivic Gauss Bonnet}, we have that $\chi(\mathcal{X}/k)$ is hyperbolic in this case. Therefore, the one exception is not a problem for our purposes. Still, it is a good question why a complete intersection of two quadrics in an odd dimensional projective space is an exception. We do not have an explanation for this. 
	\end{rem} 
	We also introduce a variant of the Jacobian ring, namely, the ring 
	$$\tilde{J} = k[Y_0,\cdots, Y_r, X_0,\cdots, X_n]/(Y_0F_0,\cdots, Y_rF_r, X_0\bar{F}_0,\cdots, X_n\bar{F}_n)$$ and show the following statement. 
	\begin{prop}\label{tildeJrho onedimensional}
		$\tilde{J}^{\rho + (r+1,n+1)}$ is a one dimensional vector space over $k$. 
	\end{prop}  
	
	\subsection{Proof of Proposition \ref{Proposition comparison products}}\label{subsection proof cup product}
	In order to prove Proposition \ref{Proposition comparison products}, we first prove two lemmas and set up some notation. 
	\begin{constr}
		Over an open $U_i$ for $i\in\{0,\cdots, r\}$, one can define a splitting of the inclusion $dF/F\wedge(-): \mathcal{O}_{\mathcal{X}}(-\mathcal{X})\to i^*\Omega_{\mathbb{P}^r\times\mathbb{P}^n}$ by 
		$$H_i\left(\sum_{i=0}^r a_idY_i + \sum_{k=0}^n b_{k}dX_k \right) = a_i F\cdot F_i^{-1}$$ 
		extending to a map $H_i:  i^*\Omega_{\mathbb{P}^r\times\mathbb{P}^n}^{a}(b\mathcal{X})\to i^*\Omega_{\mathbb{P}^r\times\mathbb{P}^n}^{a-1}((b-1)\mathcal{X})$. Similarly, over $\bar{U}_j$ for $j\in\{0,\cdots, n\}$ one can define the splitting 
		$$\bar{H}_j\left(\sum_{i=0}^r a_idY_i + \sum_{k=0}^n b_{k}dX_k \right) =
		b_{j} F\cdot \bar{F}_j^{-1}$$ 
		extending to a map $\bar{H}_j:  i^*\Omega_{\mathbb{P}^r\times\mathbb{P}^n}^{a}(b\mathcal{X})\to i^*\Omega_{\mathbb{P}^r\times\mathbb{P}^n}^{a-1}((b-1)\mathcal{X})$. We get a map on \v{C}ech cochains 
		$$H: C^q(\mathcal{U}, i^*\Omega_{\mathbb{P}^r\times\mathbb{P}^n}^{a}(b\mathcal{X}))\to C^q(\mathcal{U}, i^*\Omega_{\mathbb{P}^r\times\mathbb{P}^n}^{a-1}((b-1)\mathcal{X})) $$ defined by 
		$$H(\{s_{I,J} \}_{I,J} ) = \begin{cases} 
		\{H_{(I,J)_0}(s_{I,J}) \}_{I,J} &\text{ if } I\neq \emptyset \\
		\{\bar{H}_{(I,J)_0}(s_{I,J}) \}_{I,J} &\text{ otherwise}
		\end{cases}$$ where $(I,J)_0$ denotes the first element in an ordered index $(I,J)$. We note that $H_i = F/F_i\cdot K_i$ for $i\in\{0,\cdots, r\}$ and $\bar{H}_j = F/\bar{F}_j \cdot \bar{K}_j$ for $j\in\{0,\cdots, n\}$. 
	\end{constr} 
	The following statement is a generalization of \cite[Lemma 3.4]{LevineECHWHH} to products of projective spaces. 
	\begin{lem}\label{Levine lemma 3.4}
		Let $A\in k[Y_0,\cdots, Y_r,X_0,\cdots, X_n]$ be a polynomial of bidegree $(b-(r+1), bm - (n+1))$ for some $b$. Then:
		\begin{enumerate}
			\item $\text{res}(A\omega\wedge\bar{\omega}/F^{b+1})\in H^0(\mathcal{X}, \Omega_{\mathcal{X}}^{n+r-1}(b\mathcal{X}))$ is represented by the element 
			$\{s_0,\cdots, s_r, \bar{s}_0,\cdots, \bar{s}_n\}\in C^0(\mathcal{U}, \Omega_{\mathcal{X}}^{n+r-1}(b\mathcal{X}))$ given by 
			$$s_i = \frac{A\Omega_i\wedge\bar{\omega}}{F_iF^b}\in \Omega_{\mathcal{X}}^{n+r-1}(b\mathcal{X})(U_i)$$ for $i\in\{0,\cdots, r\}$ and 
			$$\bar{s}_j = \frac{(-1)^rA\omega\wedge\bar{\Omega}_{j}}{\bar{F}_{j}F^b}\in \Omega_{\mathcal{X}}^{n+r-1}(b\mathcal{X})(\bar{U}_j)$$ for $j\in\{0,\cdots, n\}$. 
			\item For an element $\{\frac{A}{F_I\bar{F}_JF^b}\Omega_I\wedge \bar{\Omega}_J\}_{I,J}\in C^i(\mathcal{U}, i^*\Omega_{\mathbb{P}^r\times\mathbb{P}^n}^{a-1}(b-1)\mathcal{X})$ we have that 
			$$H\left(\Bigl\{(dF/F)\wedge (\frac{A}{{F_I\bar{F}_JF^b}}\Omega_I\wedge \bar{\Omega}_J)\Bigr\}_{I,J}\right) = \Bigl\{\frac{A}{F_I\bar{F}_JF^b}\Omega_I\wedge \bar{\Omega}_J \Bigr\}_{I,J}.$$
			\item Applying the \v{C}ech differential $\delta$ to $$\alpha = \Bigl\{\frac{(-1)^{(r+\#I)\#J}A}{F_I\bar{F}_JF^b}\Omega_I\wedge \bar{\Omega}_J \Bigr\}_{\#I + \#J = i+1}\in C^i(\mathcal{U}, \Omega_{\mathcal{X}}^{a}(b\mathcal{X}))$$ we have that 
			$$\delta(\alpha) = \Biggl\{(dF/F)\wedge \left((-1)^{i+1}\frac{(-1)^{(r+\#I')\#J'}A}{F_I\bar{F}_JF^{b-1}}\Omega_{I'}\wedge \bar{\Omega}_{J'}\right)\Biggr\}_{\#I'+\#J' = i+2}$$ in $ C^{i+1}(\mathcal{U},  \Omega_{\mathcal{X}}^{a}(b\mathcal{X}))$.
			\item Let $\pi: i^*\Omega_{\mathbb{P}^r\times\mathbb{P}^n}^{a-1}((b-1)\mathcal{X})\to \Omega_{\mathcal{X}}^{a-1}((b-1)\mathcal{X})$ be the canonical projection. We have that $\pi\circ H$ is a splitting to $$dF/F\wedge (-):  C^i(\mathcal{U},  \Omega_{\mathcal{X}}^{a-1}((b-1)\mathcal{X}))\to  C^i(\mathcal{U}, i^*\Omega_{\mathbb{P}^r\times\mathbb{P}^n}^{a}(b\mathcal{X})).$$
		\end{enumerate}
	\end{lem} 
	\begin{proof} 
		This follows the method of \cite{CarlsonIVHS} directly. Note that to check an identity on a sheaf of $p$-forms on some open subset $U$ of $\mathcal{X}$, it suffices to check on $f^{-1}(U)$ for $f:\mathcal{X}’\to \mathcal{X}$ any smooth morphism with $f^{-1}(U)$ nonempty. We will take $f$ to be the restriction of $(\mathbb{A}^{r+1}\setminus\{0\})\times (\mathbb{A}^{n+1}\setminus\{0\})\to\mathbb{P}^r\times\mathbb{P}^n$ to $\mathcal{X}$. The point is that one does not need to assume that all the terms involved in the computation arise as forms on $\mathcal{X}$: the individual terms need not satisfy the Euler equations. \\
		We start by proving (1). Write $E_n = \sum_{i=0}^{n}X_i\partial/\partial X_i$ and $E_r = \sum_{i=0}^{r}Y_i\partial/\partial Y_i$, and let $dV_n = dX_0\wedge\cdots\wedge dX_n$ and $dV_r = dY_0\wedge\cdots\wedge dY_r$. Note that interior multiplication $\iota(E_n)$ with $E_n$ gives $\iota(E_n)(dV_n) = \bar{\omega}$ and similarly we have that $\iota(E_r)(dV_r) = \omega$. We note that $dF\wedge dV_r\wedge dV_n = 0$, so 
		\begin{align*} 
		0 &= \iota(E_n)\iota(E_r)(dF\wedge dV_r\wedge dV_n)\\
		&= \iota(E_n)\left(FdV_r\wedge dV_n - dF\wedge\omega\wedge dV_n\right)\\
		&= (-1)^{r+1}FdV_r\wedge\bar{\omega} - mF\omega\wedge dV_n + (-1)^rdF\wedge\omega\wedge\bar{\omega}
		\end{align*} 
		Restricting to the affine cone over $\mathcal{X}$, we have that $F=0$ and so we find that 
		\begin{equation}\label{equation dF and omegas}
		dF\wedge\omega\wedge\bar{\omega} = 0.
		\end{equation}
		Note that for $i\in\{0,\cdots, r\}$, we have that 
		$$K_idF = K_i\left(\sum_{k=0}^r F_kdY_k + \sum_{l=0}^n \bar{F}_ldX_l\right)= F_i.$$
		Applying $K_i$ to (\ref{equation dF and omegas}) therefore yields that $F_i\omega\wedge\bar{\omega} = dF\wedge\Omega_i\wedge\bar{\omega}$. Similarly, applying $\bar{K}_j$ gives that $\bar{F}_j\omega\wedge\bar{\omega} = (-1)^rdF\wedge\omega\wedge\bar{\Omega}_j$ for $j\in\{0,\cdots, n\}$. We see from this that 
		$$\frac{A\omega\wedge\bar{\omega}}{F^{b+1}} = \begin{cases}  \frac{A\Omega_i\wedge\bar{\omega}\wedge dF/F}{F_iF^b} &\text{ for } i\in\{0,\cdots, r\}\\
		\frac{(-1)^rA\omega\wedge\bar{\Omega}_j\wedge dF/F}{\bar{F}_jF^b} &\text{ for } j\in\{0,\cdots, n\}
		\end{cases} $$
		Applying the residue map to the left hand side, recalling the diagram~(\ref{diagram section 3}) from the proof of Proposition \ref{Levine Proposition 3.2}, is the same as applying $\pi$ to $\frac{A\Omega_i\wedge\bar{\omega}}{F_iF^b} $ or~$\frac{(-1)^rA\omega\wedge\bar{\Omega}_j}{\bar{F}_jF^b}$. We find the result as desired. \\
		To prove (2), let $I = \{i_0,\cdots, i_k\}\subset \{0,\cdots, r\}$ and $J = \{j_0,\cdots, j_l\}\subset \{0,\cdots, n\}$ be such that $\#I + \#J = i+1$. If $I$ is nonempty, we have that 
		$K_{i_0}\Omega_I\wedge\bar{\Omega}_J = 0$ as $K_{i_0}K_I\omega = 0$ (one removes $X_{i_0}$ from $\omega$ twice), and so 		
		\begin{align*} 
		H_{i_0}\left(dF/F\wedge \left(\frac{A}{F_I\bar{F}_JF^b}\Omega_I\wedge \bar{\Omega}_J\right)\right) &= \frac{F}{F_{i_0}}K_{i_0}\left(dF/F\wedge \left(\frac{A}{F_I\bar{F}_JF^b}\Omega_I\wedge \bar{\Omega}_J\right)\right) \\
		&= \frac{A}{F_I\bar{F}_JF^b}\Omega_I\wedge \bar{\Omega}_J
		\end{align*} 
		as desired. If $I$ is empty, one replaces $H_{i_0}$ by $\bar{H}_{j_0}$, $F_{i_0}$ by $\bar{F}_{j_0}$ and $K_{i_0}$ by $\bar{K}_{j_0}$ and the proof works in the exact same way. \\
		To prove statement (3), we consider two subsets $I' = \{i_0,\cdots, i_k\}\subset \{0,\cdots, r\}$ and $J' = \{j_0,\cdots, j_{i-k}\}\subset \{0,\cdots, n\}$ such that $\#I' + \#J' = i+2$. We apply $K_{I'}$ to the identity (\ref{equation dF and omegas}), so that we find that 
		\begin{align*} 
		0 &= (K_{i_k}\circ\cdots \circ K_{i_0})(dF\wedge\omega\wedge\bar{\omega})\\
		&= (K_{i_k}\circ\cdots \circ K_{i_1})\left(F_{i_0}\omega\wedge\bar{\omega} - dF \wedge\Omega_{i_0}\wedge\bar{\omega}\right)\\
		&= (K_{i_k}\circ\cdots \circ K_{i_2})\left(F_{i_0}\Omega_{i_1}\wedge\bar{\omega} - F_{i_1}\Omega_{i_0}\wedge\bar{\omega} + dF \wedge\Omega_{i_0i_1}\wedge\bar{\omega}\right) \\
		&= \cdots \\
		&= \sum_{l=0}^k(-1)^lF_{i_l}\Omega_{I'\setminus \{i_l\}}\wedge\bar{\omega} + (-1)^{k+1}dF\wedge\Omega_{I'}\wedge\bar{\omega}
		\end{align*} 
		Applying $\bar{K}_{J'}$ gives 
		\begin{align*} 
		0 &= (\bar{K}_{j_{i-k}}\circ\cdots\circ\bar{K}_{j_0})\left(\sum_{l=0}^k(-1)^lF_{i_l}\Omega_{I'\setminus \{i_l\}}\wedge\bar{\omega} + (-1)^{k+1}dF\wedge\Omega_{I'}\wedge\bar{\omega} \right)\\
		&= (\bar{K}_{j_{i-k}}\circ\cdots\circ\bar{K}_{j_1})(\sum_{l=0}^k(-1)^{l + r-k}F_{i_l}\Omega_{I'\setminus \{i_l\}}\wedge \bar{\Omega}_{j_0} \\
		&\quad + (-1)^{k+1}\bar{F}_{j_0}\Omega_{I'}\wedge\bar{\omega} + (-1)^{k+1+r-k}dF\wedge\Omega_{I'} \wedge\bar{\Omega}_{j_0})\\
		&= (\bar{K}_{j_{i-k}}\circ\cdots\circ\bar{K}_{j_2})(\sum_{l=0}^k(-1)^{l + 2(r-k)}F_{i_l}\Omega_{I'\setminus \{i_l\}}\wedge \bar{\Omega}_{j_0j_1} \\
		&\quad + (-1)^{k+1+ r-k-1}\bar{F}_{j_0}\Omega_{I'}\wedge\bar{\Omega}_{j_1} +  (-1)^{k+1+r-k}\bar{F}_{j_1}\Omega_{I'} \wedge\bar{\Omega}_{j_0} \\
		&\quad + (-1)^{k+1+2(r-k)}dF\wedge\Omega_{I'}\wedge\bar{\Omega}_{j_0j_1})\\
		&= \cdots \\
		&= \sum_{l=0}^k(-1)^{l + (i-k+1)(r-k)}F_{i_l}\Omega_{I'\setminus\{i_l\}}\wedge\bar{\Omega}_{J'} \\
		&\quad  + (-1)^{k+1}\sum_{l=0}^{i-k}(-1)^{l + (i-k)(r-k-1)}\bar{F}_{j_l}\Omega_{I'}\wedge\bar{\Omega}_{J'\setminus \{j_l\}} \\
		&\quad + (-1)^{k+1+(i-k+1)(r-k)}dF\wedge\Omega_{I'}\wedge\bar{\Omega}_{J'}
		\end{align*} 
		so 
		$$\sum_{l=0}^k(-1)^{l}F_{i_l}\Omega_{I'\setminus\{i_l\}}\wedge\bar{\Omega}_{J'}
		+ (-1)^{k+1 + r-i}\sum_{l=0}^{i-k}(-1)^{l}\bar{F}_{j_l}\Omega_{I'}\wedge\bar{\Omega}_{J'\setminus \{j_l\}} = (-1)^{k}dF\wedge\Omega_{I'}\wedge \bar{\Omega}_{J'} $$
		We see that 
		\begin{align*} 
		&\left(\delta\left(\Big\{\frac{(-1)^{(r+\#I)\#J}A}{F_I\bar{F}_JF^b}\Omega_I\wedge \bar{\Omega}_J\Big\}_{I,J}\right)\right)_{I',J'} \\
		&= \frac{A}{F^b}\sum_{l=0}^{k} (-1)^{l+(r+k)(i-k+1)}\frac{F_{i_l}\Omega_{I'\setminus \{i_l\}}\wedge\bar{\Omega}_{J'}}{F_{I'}\bar{F}_{J'}} \\
		&\quad + (-1)^{k+1}\frac{A}{F^b}\sum_{l=0}^{i-k} (-1)^{l+(r+k+1)(i-k)}\frac{\bar{F}_{j_l}\Omega_{I'}\wedge\bar{\Omega}_{J'\setminus \{j_l\}}}{F_{I'}\bar{F}_{J'}} \\
		&= (-1)^{i+1}\frac{(-1)^{(r+k+1)(i-k-1)}A}{F_{I'}\bar{F}_{J'}F^b}dF\wedge\Omega_{I'}\wedge\bar{\Omega}_{J'}
		\end{align*}  which proves the claim. \\
		Finally, for (4), we note that as the original maps $H_j$ are splittings, the composition $\pi\circ H$ is one too, which completes the proof. 
	\end{proof} 
	We can use this to prove the following result, which is a generalization of \cite[Proposition 3.6]{LevineECHWHH} to products of projective spaces. 
	\begin{lem}\label{Levine 3.6}
		Let $A\in J^{q-r, (q+1)m - (n+1)}$, then $\psi_q(A)\in H^q(\mathcal{X}, \Omega_{\mathcal{X}}^{p})$ is represented by the \v{C}ech cochain $\{\pi((-1)^{(r+\#I)\#J}A\Omega_I\wedge\bar{\Omega}_J/F_I\bar{F}_J)\}_{I,J}$ in $C^q(\mathcal{U}, \Omega_{\mathcal{X}}^{p})$. 
	\end{lem} 
	\begin{proof}
		Recall that for $j\in\{0,\cdots, q-1\}$, there is the exact sequence (\ref{Levine lemma 3.1(1) exact sequence proof}) given by
		$$0\to \Omega_{\mathcal{X}}^{n+r-j-2}((q-j-1)\mathcal{X}) \to i^*\Omega_{\mathbb{P}^r\times\mathbb{P}^n}^{n+r-j-1}((q-j)\mathcal{X}) \to \Omega_{\mathcal{X}}^{n+r-j-1}((q-j)\mathcal{X})\to 0.$$ 
		We find coboundary maps 
		$$\delta_j: H^j(\mathcal{X}, \Omega_{\mathcal{X}}^{n+r-j-1}((q-j)\mathcal{X})\to   H^{j+1}(\mathcal{X}, \Omega_{\mathcal{X}}^{n+r-j-2}((q-j-1)\mathcal{X}).$$
		We will show by induction on $j$ that $\delta_{j-1}\circ\cdots\circ\delta_0 (\text{res}(A\omega\wedge\bar{\omega}/F^{q+1}))$ is represented by the \v{C}ech cocycle $$\{\pi((-1)^{j(j-1)/2+(r+\#I)\#J}A\Omega_I\wedge\bar{\Omega}_J/F_I\bar{F}_JF^{q-j}) \}_{I,J}\in C^j(\mathcal{U}, \Omega_{\mathcal{X}}^{n+r-j-1}(q-j)\mathcal{X}).$$ Then taking $j=q$ will give the desired result: using \cite[Remark 2.3]{LevineECHWHH}, we have that the element above represents $\psi_q(A)$ for $j=q$ up to a factor $(-1)^{q(q-1)/2}$. \\
		First of all, note that the case where $j=0$ is precisely Lemma \ref{Levine lemma 3.4}, part (1). Now assume that $\delta_{j-1}\circ\cdots\circ\delta_0 (\text{res}(A\omega\wedge\bar{\omega}/F^{q+1}))$ is represented by the \v{C}ech cocycle $\{\pi((-1)^{j(j-1)/2}(-1)^{(r+\#I)\#J}A\Omega_I\wedge\bar{\Omega}_J/F_I\bar{F}_JF^{q-j}) \}_{I,J}$ for some~$j$. Then $\delta_{j}\circ\cdots\circ\delta_0 (\text{res}(A\omega\wedge\bar{\omega}/F^{q+1}))$ is represented by the coboundary of $$\{\pi((-1)^{j(j-1)/2}(-1)^{(r+\#I)\#J}A\Omega_I\wedge\bar{\Omega}_J/F_I\bar{F}_JF^{q-j}) \}_{I,J}.$$ Using \cite[Remark 2.2]{LevineECHWHH}, this is defined by lifting to the cochain $$\{(-1)^{j(j-1)/2+(r+\#I)\#J}A\Omega_I\wedge\bar{\Omega}_J/F_I\bar{F}_JF^{q-j} \}_{I,J}\in C^j(\mathcal{U}, i^*\Omega_{\mathbb{P}^r\times\mathbb{P}^n}^{n+r-j-1}((q-j)\mathcal{X})$$ and applying the negative of the \v{C}ech coboundary operator $\delta$, and then viewing this as an element coming from $C^{j+1}(\mathcal{U}, \Omega_{\mathcal{X}}^{n+r-j-2}((q-j-1)\mathcal{X}))$ of which the inclusion is induced by the map $dF/F\wedge(-)$. Using Lemma \ref{Levine lemma 3.4}, part (3) we have that 
		\begin{align*}
		-\delta&\left(\Big\{\frac{(-1)^{j(j-1)/2}(-1)^{(r+\#I)\#J}A\Omega_I\wedge \bar{\Omega}_J}{F_I\bar{F}_JF^{q-j}}\Big\}_{I,J}\right) \\
		&= \Big\{(dF/F)\wedge \left(\frac{(-1)^{j(j+1)/2}(-1)^{(r+\#I)\#J}A\Omega_I\wedge \bar{\Omega}_J}{F_I\bar{F}_JF^{q-j-1}}\right)\Big\}_{I,J}
		\end{align*} 
		By part (4), $\pi\circ H$ provides a splitting to $dF/F\wedge(-)$ implying by part (2) that the desired element is 
		\begin{align*} 
		(\pi\circ H)&\left(\Big\{\frac{dF}{F}\wedge \left(\frac{(-1)^{j(j+1)/2}(-1)^{(r+\#I)\#J}A\Omega_I\wedge \bar{\Omega}_J}{F_I\bar{F}_JF^{q-j-1}}\right)\Big\}_{I,J}\right) \\
		&= \Big\{\pi \left((-1)^{j(j+1)/2}\frac{(-1)^{(r+\#I)\#J}A\Omega_I\wedge \bar{\Omega}_J}{F_I\bar{F}_JF^{q-j-1}}\right)\Big\}_{I,J}\end{align*}  completing the induction. 
	\end{proof}
	\begin{proof}[Proof of Proposition \ref{Proposition comparison products}]
		We know from Lemma \ref{Levine 3.6} that $\omega_A$ and $\omega_B$ are represented by the cochains 
		$$\{\pi((-1)^{(r+\#I)\#J}A\Omega_I\wedge\bar{\Omega}_J/F_I\bar{F}_J)\}_{I,J}\in C^q(\mathcal{U}, \Omega_{\mathcal{X}}^{n+r-1-q})$$ and $$\{\pi((-1)^{(r+\#I)\#J}B\Omega_I\wedge\bar{\Omega}_J/F_I\bar{F}_J)\}_{I,J}\in C^p(\mathcal{U}, \Omega_{\mathcal{X}}^{n+r-1-p}).$$ 
		Now let $I = \{i_0,\cdots, i_{r-1}\}\subset \{0,\cdots, r\}$ and $J = \{j_0,\cdots, j_{n-1}\}\subset\{0,\cdots, n\}$. First assume that $q\leq r-1$. Then by definition of the cup product on \v{C}ech cochains, we have that 
		\begin{align*}
		(\omega_A\cup\omega_B)_{I,J} &= \pi\left((-1)^{(2r-q)n}\frac{A\Omega_{i_0,\cdots, i_q}\wedge\bar{\omega}}{F_{i_0,\cdots, i_q}}\cdot \frac{B\Omega_{i_q,\cdots, i_{r-1}}\wedge\bar{\Omega}_{j_0,\cdots, j_{n-1}}}{F_{i_q,\cdots, i_{r-1}}\bar{F}_{j_0,\cdots, j_{n-1}}}\right) \\
		&=  \pi\left(AB\Omega(I,J)\right)
		\end{align*}
		Similarly, if $q\geq r$, we have that 
		\begin{align*}
		&(\omega_A\cup\omega_B)_{I,J}\\
		&= \pi\left((-1)^{2r(q-r+1) + r(n-q-r)}\frac{A\Omega_{i_0,\cdots, i_{r-1}}\wedge\bar{\Omega}_{j_0,\cdots, j_{q-r}}}{F_{i_0,\cdots, i_{r-1}}\bar{F}_{j_0,\cdots, j_{q-r}}}\cdot \frac{B\omega\wedge\bar{\Omega}_{j_{q-r},\cdots, j_{n-1}}}{\bar{F}_{j_{q-r},\cdots, j_{n-1}}}\right) \\
		&=  \pi\left(AB\Omega(I,J)\right)
		\end{align*}
		which proves the statement. 
	\end{proof} 
	
	\subsection{Proof of Proposition \ref{representation of image cech} and Corollary \ref{corollary definition map Jrho to cohomology}}
	In order to prove Proposition \ref{representation of image cech}, we first show the following lemma. 
	\begin{lem}\label{Carlsson Griffiths lemma}
		Consider subsets $$I = \{i_0,\cdots, i_{r-1}\}\subset\{0,\cdots, r\} \text{ and } J = \{j_0,\cdots, j_{n-1}\}\subset\{0,\cdots, n\}$$ such that $I = \{0,\cdots, r\}\setminus\{v\}$ and $J = \{0,\cdots, n\}\setminus\{w\}$. We have that $$dF\wedge \Omega(I,J) = \frac{(-1)^{v+w}Y_vF_vX_w\bar{F}_w\omega\wedge\bar{\omega}}{\prod_{i=0}^rF_i\prod_{j=0}^n\bar{F}_j}.$$
	\end{lem} 
	\begin{proof}
		If $q\leq r-1$, we have by \cite[Lemma on page 14]{CarlsonIVHS} that 
		$$\left(\sum_{i=0}^rF_idY_i\right)\wedge\Omega_{i_0,\cdots, i_q}\wedge\Omega_{i_q,\cdots, i_{r-1}} = (-1)^vY_vF_{i_q}\wedge\omega.$$
		Now noting that $\bar{\Omega}_{j_0,\cdots, j_{n-1}} = (-1)^wX_w$ and $dX_j\wedge\bar{\omega}=0$ for all $j\in\{0,\cdots, n\}$, we see that 
		\begin{align*}
		dF\wedge\Omega(I,J) &= \left(\sum_{i=0}^rF_idY_i\right)\wedge\frac{(-1)^{nq}F_v\bar{F}_w\Omega_{i_0,\cdots, i_q}\wedge\bar{\omega}\wedge\Omega_{i_q,\cdots, i_{r-1}}\wedge\bar{\Omega}_{j_0,\cdots, j_{n-1}}}{F_{i_q}\prod_{i=0}^rF_i\prod_{j=0}^n\bar{F}_n}\\
		&\quad + \left(\sum_{j=0}^n\bar{F}_jdX_j\right)\wedge\frac{(-1)^{nq}F_v\bar{F}_w\Omega_{i_0,\cdots, i_q}\wedge\bar{\omega}\wedge\Omega_{i_q,\cdots, i_{r-1}}\wedge\bar{\Omega}_{j_0,\cdots, j_{n-1}}}{F_{i_q}\prod_{i=0}^rF_i\prod_{j=0}^n\bar{F}_n}\\
		&= \left(\sum_{i=0}^rF_idY_i\right)\wedge\frac{(-1)^wX_wF_v\bar{F}_w\Omega_{i_0,\cdots, i_q}\wedge\Omega_{i_q,\cdots, i_{r-1}}\wedge\bar{\omega}}{F_{i_q}\prod_{i=0}^rF_i\prod_{j=0}^n\bar{F}_n}\\
		&\quad + \left(\sum_{j=0}^n\bar{F}_jdX_j\right)\wedge\frac{(-1)^{n(r-1)+w}X_wF_v\bar{F}_w\bar{\omega}\wedge\Omega_{i_0,\cdots, i_q}\wedge\Omega_{i_q,\cdots, i_{r-1}}}{F_{i_q}\prod_{i=0}^rF_i\prod_{j=0}^n\bar{F}_n}\\
		&= \frac{(-1)^{v+w}Y_vF_vX_w\bar{F}_w\omega\wedge\bar{\omega}}{\prod_{i=0}^rF_i\prod_{j=0}^n\bar{F}_n}
		\end{align*}
		If $q\geq r$, 
		we have by \cite[Lemma on page 14]{CarlsonIVHS} that 
		$$\left(\sum_{j=0}^n\bar{F}_jdX_j\right)\wedge\bar{\Omega}_{j_0,\cdots, j_{q-r}}\wedge\bar{\Omega}_{j_{q-r},\cdots, j_{n-1}} = (-1)^wX_w\bar{F}_{j_{q-r}}\bar{\omega}.$$
		Now as $\Omega_{i_0,\cdots, i_{r-1}} = (-1)^vY_v$ and $dY_i\wedge\omega=0$ for any $i\in\{0,\cdots, r\}$, we see that 
		\begin{align*}
		&dF\wedge\Omega(I,J)\\
		&= \left(\sum_{i=0}^rF_idY_i\right)\wedge\frac{(-1)^{r(n+q+r)}F_v\bar{F}_w\Omega_{i_0,\cdots, i_{r-1}}\wedge\bar{\Omega}_{j_0,\cdots, j_{q-r}}\wedge\omega\wedge\bar{\Omega}_{j_{q-r}, \cdots, j_{n-1}}}{\bar{F}_{j_{q-r}}\prod_{i=0}^rF_i\prod_{j=0}^n\bar{F}_j}\\
		&\quad + \left(\sum_{j=0}^n\bar{F}_jdX_j\right)\wedge\frac{(-1)^{r(n+q+r)}F_v\bar{F}_w\Omega_{i_0,\cdots, i_{r-1}}\wedge\bar{\Omega}_{j_0,\cdots, j_{q-r}}\wedge\omega\wedge\bar{\Omega}_{j_{q-r}, \cdots, j_{n-1}}}{\bar{F}_{j_{q-r}}\prod_{i=0}^rF_i\prod_{j=0}^n\bar{F}_j}\\
		&= \left(\sum_{i=0}^rF_idY_i\right)\wedge\frac{(-1)^{v + r}Y_vF_v\bar{F}_w\omega\wedge\bar{\Omega}_{j_0,\cdots, j_{q-r}}\wedge\bar{\Omega}_{j_{q-r}, \cdots, j_{n-1}}}{\bar{F}_{j_{q-r}}\prod_{i=0}^rF_i\prod_{j=0}^n\bar{F}_j}\\
		&\quad + \left(\sum_{j=0}^n\bar{F}_jdX_j\right)\wedge\frac{(-1)^{v + rn }Y_vF_v\bar{F}_w\bar{\Omega}_{j_0,\cdots, j_{q-r}}\wedge\bar{\Omega}_{j_{q-r}, \cdots, j_{n-1}}\wedge\omega}{\bar{F}_{j_{q-r}}\prod_{i=0}^rF_i\prod_{j=0}^n\bar{F}_j}\\
		&= \frac{(-1)^{v+w}Y_vF_vX_w\bar{F}_w\omega\wedge\bar{\omega}}{\prod_{i=0}^rF_i\prod_{j=0}^n\bar{F}_n}
		\end{align*}
		which completes the proof. 
	\end{proof}
	\begin{proof}[Proof of Proposition \ref{representation of image cech}]
		Using \cite[Remark 2.2]{LevineECHWHH} together with Lemma \ref{ilowerstar is boundary}, we can represent $i_*(\omega_A\cup\omega_B))$ by lifting to the section $AB\Omega(I,J)\wedge dF/F$ of $\Omega_{\mathbb{P}^r\times\mathbb{P}^n}^{n+r}(\log(\mathcal{X}))$ and then taking the negative of the \v{C}ech coboundary. Note that we use the diagram from the proof of Proposition \ref{Levine Proposition 3.2} again to see that this is really the lift. \\		
		Now note that $C^{n+r}(\mathcal{U},\Omega_{\mathbb{P}^r\times\mathbb{P}^n}^{n+r})$ has indices $(I',J')$ where either:
		\begin{itemize} 
			\item $I'= \{0,\cdots, r\}$ and $J'= \{0,\cdots, n\}\setminus \{w\}$ for a certain $w$
			\item $I'= \{0,\cdots, r\}\setminus\{v\}$ and $J'= \{0,\cdots, n\}$ for a certain $v$. 
		\end{itemize} 
		In the first case, we have using Lemma \ref{Carlsson Griffiths lemma} that
		\begin{align*} 
		\delta(\{AB&dF/F\wedge\Omega(I,J)\}_{I,J})_{I',J'} 
		\\
		&= AB\sum_{v=0}^r (-1)^v (-1)^{v+w}\frac{Y_vF_vX_w\bar{F}_{w}\omega\wedge\bar{\omega}}{F \prod_{i=0}^rF_i\prod_{j=0}^n\bar{F}_{j}} \\
		&= (-1)^{w}AB\frac{X_w\bar{F}_w\omega\wedge\bar{\omega}}{ \prod_{i=0}^rF_i\prod_{j=0}^n\bar{F}_{j}} 
		\end{align*} 
		If the index $(I',J')$ is of the second form, we similarly find that 
		\begin{align*} 
		\delta(\{AB (dF/F)&\wedge\Omega(I,J)\}_{I,J})_{I',J'} \\
		&= AB\sum_{w=0}^n (-1)^{v+w} (-1)^{w+r} \frac{Y_vF_{v}X_w \bar{F}_w\omega\wedge\bar{\omega}}{F \prod_{i=0}^rF_i\prod_{j=0}^n\bar{F}_{j}} \\
		&= (-1)^{v+r}mAB \frac{Y_vF_v\omega\wedge\bar{\omega}}{ \prod_{i=0}^rF_i\prod_{j=0}^n\bar{F}_{j}} 
		\end{align*} 
		This completes the proof.
	\end{proof} 
	\begin{proof}[Proof of Corollary \ref{corollary definition map Jrho to cohomology}]
		Let $D\in k[Y_0,\cdots, Y_r, X_0,\cdots, X_n]^{\rho}$. Note that $\tilde{\phi}(D)$ is a \v{C}ech cocycle, as we have that 
		\begin{align*} 
		\delta(\tilde{\phi}(D)) &= \frac{D\left(m\sum_{v=0}^r(-1)^{2v+r+1}Y_vF_v\omega\wedge\bar{\omega} + \sum_{w=0}^n(-1)^{2w+r+2}X_w\bar{F}_w\omega\wedge\bar{\omega}\right)} {\prod_{i=0}^rF_i\prod_{j=0}^n\bar{F}_j}\\
		&= \frac{(-1)^{r+1}D(mF - mF)\omega\wedge\bar{\omega}} {\prod_{i=0}^rF_i\prod_{j=0}^n\bar{F}_j}\\
		&= 0
		\end{align*} 
		This means that $\tilde{\phi}$ induces a map 
		$$\bar{\phi}: k[Y_0,\cdots, Y_r, X_0,\cdots, X_n]^{\rho}\to H^{n+r}(\mathbb{P}^r\times\mathbb{P}^n,\Omega_{\mathbb{P}^r\times\mathbb{P}^n}^{n+r}).$$ Using Proposition \ref{representation of image cech}, we find the following commutative diagram
		\[
		\begin{tikzcd} 
		H^q(\mathcal{X}, \Omega_{\mathcal{X}}^p)_{prim}\otimes H^p(\mathcal{X}, \Omega_{\mathcal{X}}^q)_{prim} \arrow[r, "i_*\circ \cup"] & H^{n+r}(\mathbb{P}^r\times\mathbb{P}^n, \Omega_{\mathbb{P}^r\times\mathbb{P}^n}^{n+r}) \\
		k[Y,X]^{q-r,(q+1)m - (n+1)} \otimes k[Y,X]^{p-r,(p+1)m - (n+1)}  \arrow[r] \arrow[u, "\tilde{\psi}_p\otimes\tilde{\psi}_q"]  & k[Y,X]^\rho \arrow[u, "\bar{\phi}"] 
		\end{tikzcd} 
		\]
		where we denote $k[Y,X] = k[Y_0,\cdots, Y_r,X_0,\cdots, X_n]$. By Proposition \ref{Levine Proposition 3.2}, we have that $\tilde{\psi}_q$ descends to an isomorphism
		$$J^{q-r,(q+1)m-n-1}\to H^q(\mathcal{X},\Omega_{\mathcal{X}}^p)_{prim}$$
		and similarly for $\tilde{\psi}_p$. \\
		We note that $\bar{\phi}$ maps the Jacobian ideal to zero. To see this, suppose that $D\in k[Y_0,\cdots, Y_r, X_0,\cdots, X_n]^\rho$ is a multiple of $F_i$ for some $i\in\{0,\cdots, r\}$. Write $\tilde{\phi}(D) = \{s_0,\cdots, s_r,\bar{s}_0,\cdots, \bar{s}_n\}$. Then $\{\xi_{I,J}\}_{I,J}\in C^{n+r-1}(\mathcal{U},\Omega_{\mathbb{P}^r\times\mathbb{P}^n}^{n+r})
		$ given by 
		$$\xi_{I,J}  = \begin{cases}
		(-1)^{i-1}s_v &\text{ if } I = \{0,\cdots, r\}\setminus \{v,i\}, J = \{0,\cdots, n\} \text{ and } v<i \\
		(-1)^{i}s_v &\text{ if } I = \{0,\cdots, r\}\setminus \{v,i\}, J = \{0,\cdots, n\} \text{ and } v>i \\
		(-1)^{i}\bar{s}_w &\text{ if } I = \{0,\cdots, r\}\setminus \{i\}, J = \{0,\cdots, n\}\setminus\{w\} \\
		0 &\text{ otherwise}
		\end{cases}$$
		satisfies $\delta(\xi_{I,J}) = \{s_0,\cdots, s_r,\bar{s}_0,\cdots, \bar{s}_n\}$. Indeed, for $I' = \{0,\cdots, r \}\setminus\{v\}$ with $v < i$ and $J' = \{0,\cdots, n\}$, we have that $\delta(\xi_{I,J})_{I',J'} = (-1)^{i-1}(-1)^{i-1}s_v = s_v$, and similarly for $v>i$. For  $I' = \{0,\cdots, r \}$ and $J' = \{0,\cdots, n\}\setminus\{w\}$, we have that $\delta(\xi_{I,J})_{I',J'} = (-1)^{i}(-1)^{i}\bar{s}_w = \bar{s}_w$. Finally, for $I' = \{0,\cdots, r \}\setminus\{i\}$ and $J' = \{0,\cdots, n\}$, we have that $$\delta(\xi_{I,J})_{I',J'} = \sum_{v\neq i}(-1)^{v+i-1}s_v + \sum_{w=0}^r(-1)^{i+w+r}\bar{s}_w = -(-1)^{i-1}(-1)^i s_i = s_i.$$
		Therefore, $\tilde{\phi}(D)$ is a coboundary. Similarly, if $D$ is a multiple of $\bar{F}_j$ for some $j\in\{0,\cdots, n\}$, then the element $\{\xi_{I,J}\}_{I,J}\in C^{n+r-1}(\mathcal{U},\Omega_{\mathbb{P}^r\times\mathbb{P}^n}^{n+r})
		$ given by 
		$$\xi_{I,J}  = \begin{cases}
		(-1)^{j+r}s_v &\text{ if } I = \{0,\cdots, r\}\setminus \{v\}, J = \{0,\cdots, n\}\setminus\{j\}\\
		(-1)^{j+r}\bar{s}_w &\text{ if } I = \{0,\cdots, r\}, J = \{0,\cdots, n\}\setminus\{w,j\} \text{ and } w<j \\
		(-1)^{j+r+1}\bar{s}_w &\text{ if } I = \{0,\cdots, r\}, J = \{0,\cdots, n\}\setminus\{w,j\} \text{ and } w>j \\
		0 &\text{ otherwise}
		\end{cases}$$
		satisfies $\delta(\xi_{I,J}) = \{s_0,\cdots, s_r,\bar{s}_0,\cdots, \bar{s}_n\}$. \\	 
		So $\bar{\phi}$ descends to a map $\phi: J^\rho\to H^{n+r}(\mathbb{P}^r\times\mathbb{P}^n, \Omega_{\mathbb{P}^r\times\mathbb{P}^n}^{n+r})$ which makes the diagram commute. As the cup product is non-degenerate, we have that $\phi$ is surjective. 
	\end{proof}
	
	\subsection{One dimensionality of $J^\rho$ and $\tilde{J}^{\rho + (r+1,n+1)}$}\label{section onedimensionality of Jrho}
	In this section, we will prove the following statement. 
	\begin{prop}\label{one dimensionality Jrho}
		$J^\rho$ is a one dimensional vector space over $k$, except possibly if $n$ is odd, $r=1$ and $m=2$. 
	\end{prop}
	Over $\mathbb{C}$, this is a special case of \cite[Lemma 6.3]{KonnoVTPCI} and the argument is partially the same. The idea to use the bundle $\Sigma_{\mathcal{L}}$ (see below) to give a description of the Jacobian ring and study its duality properties goes back to \cite[Section 2]{Green}. It will follow that the map from Corollary \ref{corollary definition map Jrho to cohomology} is an isomorphism whenever we are not in the situation where $n$ is odd, $r=1$ and $m=2$. 
	\begin{nota}
		Write $\mathcal{L} = \mathcal{O}(1,m)$ and let $\Sigma_{\mathcal{L}}$ be the bundle as defined in \cite[Section 2.1]{KonnoVTPCI}. There is a global presentation of $\Sigma_{\mathcal{L}}$ given by 
		\begin{align*} 
		0 &\to e_1\cdot\mathcal{O}_{\mathbb{P}^r}\oplus e_2\cdot\mathcal{O}_{\mathbb{P}^n}\to \\
		&\quad \text{id}_{\mathcal{L}}\cdot\mathcal{O}_{\mathbb{P}^r\times\mathbb{P}^n} \oplus \bigoplus_{i=0}^r \frac{\partial}{\partial Y_i}\cdot\mathcal{O}_{\mathbb{P}^r\times\mathbb{P}^n}(1,0)\oplus \bigoplus_{j=0}^n \frac{\partial}{\partial X_j}\cdot\mathcal{O}_{\mathbb{P}^r\times\mathbb{P}^n}(0,1)\to \Sigma_{\mathcal{L}}\to 0
		\end{align*} 
		where the first map is given by $$e_1\mapsto -\text{id}_{\mathcal{L}} + \sum_{i=0}^r Y_i\frac{\partial}{\partial Y_i}, e_2\mapsto -m\text{id}_{\mathcal{L}} + \sum_{j=0}^n X_j\frac{\partial}{\partial X_j}.$$ 
		We can map the above sequence into the Euler sequence 
		$$0\to \mathcal{O}_{\mathbb{P}^r}\oplus\mathcal{O}_{\mathbb{P}^n}\to \bigoplus_{i=0}^r\mathcal{O}(1,0)\oplus \bigoplus_{j=0}^n\mathcal{O}(0,1)\to T_{\mathbb{P}^r\times\mathbb{P}^n}\to 0$$ by sending $\text{id}_{\mathcal{L}}$ to zero. This yields the exact sequence 
		\begin{equation}\label{sigma tangent sequence}
		0\to \mathcal{O}_{\mathbb{P}^r\times\mathbb{P}^n}\to \Sigma_{\mathcal{L}}\to T_{\mathbb{P}^r\times\mathbb{P}^n}\to 0.
		\end{equation}
	\end{nota} 
	Consider the morphism $\Sigma_\mathcal{L}\to \mathcal{L}$ which sends a local section
	$$a\text{id}_{\mathcal{L}} + \sum_{i=0}^r b_i\frac{\partial }{\partial Y_i} + \sum_{j=0}^n c_j\frac{\partial }{\partial X_j}$$ to 
	$$aF + \sum_{i=0}^rb_iF_i  + \sum_{j=0}^nc_j\bar{F}_j.$$
	This gives rise to a surjective morphism $\Sigma_\mathcal{L}\otimes\mathcal{L}^{-1}\to \mathcal{O}_{\mathbb{P}^r\times\mathbb{P}^n}$. Form the associated Koszul complex 
	\begin{equation}
	0\to \Lambda^{n+r+1}\Sigma_{\mathcal{L}}\otimes \mathcal{L}^{-n-r-1}\to \cdots \to \Lambda^2\Sigma_\mathcal{L}\otimes\mathcal{L}^{-2}\to  \Sigma_\mathcal{L}\otimes\mathcal{L}^{-1}\to \mathcal{O}_{\mathbb{P}^r\times\mathbb{P}^n}\to 0.
	\end{equation}	
	Now applying the functor $\text{Hom}(-, \Omega_{\mathbb{P}^n\times\mathbb{P}^r}^{n+r})$ we find the exact sequence, 
	\begin{equation}\label{Koszul complex for J}
	0\to \Omega_{\mathbb{P}^n\times\mathbb{P}^r}^{n+r}\to \Sigma_{\mathcal{L}}^\vee\otimes\mathcal{L}\otimes\Omega_{\mathbb{P}^n\times\mathbb{P}^r}^{n+r} \to \cdots \to \Lambda^{n+r+1}\Sigma_{\mathcal{L}}^\vee\otimes \mathcal{L}^{n+r+1}\otimes\Omega_{\mathbb{P}^n\times\mathbb{P}^r}^{n+r}\to 0
	\end{equation}	
	defining a left resolution of $\Omega^{n+r}_{\mathbb{P}^r\times\mathbb{P}^n}$. Note that 
	\begin{align*}
	\Lambda^{n+r+1}\Sigma_{\mathcal{L}}^\vee\otimes \mathcal{L}^{n+r+1}\otimes\Omega_{\mathbb{P}^n\times\mathbb{P}^r}^{n+r} &= \mathcal{O}(\rho)
	\end{align*}
	as $\det(\Sigma_{\mathcal{L}}) = \det(T_{\mathbb{P}^r\times\mathbb{P}^n}) = \mathcal{O}(r+1,n+1)$ using the exact sequence (\ref{sigma tangent sequence}), $\mathcal{L}^{n+r+1} = \mathcal{O}(n+r+1,m(n+r+1))$ and $\Omega_{\mathbb{P}^n\times\mathbb{P}^r}^{n+r} = \mathcal{O}(-r-1,-n-1)$. \\
	Consider the hypercohomology spectral sequence 
	$$E_1^{p,q} = H^q(\mathbb{P}^r\times\mathbb{P}^n, \Lambda^{p+1}\Sigma_{\mathcal{L}}^\vee\otimes\mathcal{L}^{p+1}\otimes\Omega_{\mathbb{P}^n\times\mathbb{P}^r}^{n+r}) \implies H^{p+q}(\mathbb{P}^r\times\mathbb{P}^n, \Omega_{\mathbb{P}^n\times\mathbb{P}^r}^{n+r}).$$	
	\begin{lem}\label{spectral sequence lemma jrho onedimensional}
		We have that $E_1^{p,q} = 0$ for $q>0$ and either $p+q = n+r$ or $p+q = n+r-1$, except in the case when $n$ is odd, $r = 1$ and $m = 2$. 
	\end{lem}
	\begin{proof}[Proof of Lemma \ref{spectral sequence lemma jrho onedimensional}]
		This proof is partially taken from \cite{KonnoVTPCI}. Note that the sequence (\ref{sigma tangent sequence}) gives rise to exact sequences 
		$$0\to \Omega_{\mathbb{P}^n\times\mathbb{P}^r}^j\to \Lambda^j\Sigma_{\mathcal{L}}^\vee \to \Omega_{\mathbb{P}^n\times\mathbb{P}^r}^{j-1}\to 0$$
		so that it is enough to prove that 	for $1\leq s\leq n+r-1$, we have that 
		$$H^{n+r-s}(\mathbb{P}^r\times\mathbb{P}^n ,\Omega_{\mathbb{P}^r\times\mathbb{P}^n }^s\otimes\mathcal{L}^{s}\otimes \Omega_{\mathbb{P}^r\times\mathbb{P}^n }^{n+r}) = 0$$ and 
		$$H^{n+r-s}(\mathbb{P}^r\times\mathbb{P}^n ,\Omega_{\mathbb{P}^r\times\mathbb{P}^n }^{s-1}\otimes\mathcal{L}^{s}\otimes \Omega_{\mathbb{P}^r\times\mathbb{P}^n}^{n+r}) = 0$$
		and for $1\leq s\leq n+r$, we have that 
		$$H^{n+r + 1 -s}(\mathbb{P}^r\times\mathbb{P}^n ,\Omega_{\mathbb{P}^r\times\mathbb{P}^n }^s\otimes\mathcal{L}^{s}\otimes \Omega_{\mathbb{P}^r\times\mathbb{P}^n}^{n+r}) = 0 $$ and 
		$$H^{n+r + 1 -s}(\mathbb{P}^r\times\mathbb{P}^n ,\Omega_{\mathbb{P}^r\times\mathbb{P}^n }^{s-1}\otimes\mathcal{L}^{s}\otimes \Omega_{\mathbb{P}^r\times\mathbb{P}^n}^{n+r}) = 0.$$
		We show that the first condition holds except in the case when $n$ is odd, $r = 1$ and $m = 2$; the others are similar. That is, we will show that 
		\begin{equation}\label{condition to check}
		H^{n+r-s}(\mathbb{P}^r\times\mathbb{P}^n ,\Omega_{\mathbb{P}^r\times\mathbb{P}^n }^s(s-r-1,ms-n-1)) = 0
		\end{equation} for $1\leq s\leq n+r-1$ except in the case when $n$ is odd, $r = 1$ and $m = 2$. Using Proposition \ref{Botts theorem for general products}, we have that 
		\begin{align*}
		H^{n+r-s}&(\mathbb{P}^r\times\mathbb{P}^n ,\Omega_{\mathbb{P}^r\times\mathbb{P}^n }^s(s-r-1,ms-n-1)) \\
		&= \bigoplus_{i+j =s} \bigoplus_{k+l=n+r-s} H^k(\mathbb{P}^r,\Omega^i_{\mathbb{P}^r}(s-r-1))\otimes H^l(\mathbb{P}^n,\Omega^j_{\mathbb{P}^n}(ms-n-1))
		\end{align*}
		Note that by Theorem \ref{Botts theorem} we have that $H^k(\mathbb{P}^r,\Omega^i_{\mathbb{P}^r}(s-r-1)) = 0$ except possibly if we are in one of the following situations:
		\begin{enumerate}
			\item $k>0, j=0$ and $s\neq r+1$. Note that $H^l(\mathbb{P}^n,\mathcal{O}(ms-n-1))$ is zero except possibly for: 
			\begin{itemize}
				\item $l=0$. Then $k = n+r-s$. If $k>r$, then $H^k(\mathbb{P}^r,\Omega^i_{\mathbb{P}^r}(s-r-1)) = 0$. Otherwise, $n+r-s\leq r$ and so $s\geq n$ which implies that $i\geq n > r$ and so $H^k(\mathbb{P}^r,\Omega^i_{\mathbb{P}^r}(s-r-1)) = 0$ in this case as well. 
				\item $l=n$. Then as $ms-n-1 > -n-1$, we have that $$H^n(\mathbb{P}^n,\mathcal{O}(ms-n-1)) = \left(\frac{1}{X_0,\cdots, X_n} k[X_0^{-1},\cdots, X_n^{-1}]\right)_{ms-n-1} =0.$$
			\end{itemize}
			\item $s = r+1$ and $i=k$. We note that: 
			\begin{itemize}
				\item If $l=0$ then $k = n-1> r$ and so $H^k(\mathbb{P}^r,\Omega^i_{\mathbb{P}^r}(s-r-1)) = 0$. 
				\item If $l>0$ and $ms\neq n+1$ then $ms-r-1 \geq j-n$ as \begin{equation}\label{msr equation} ms = m(i+j) \geq 2j\geq j+1.\end{equation} It follows that $H^l(\mathbb{P}^n,\Omega^j_{\mathbb{P}^n}(ms-n-1)) =0$. 
				\item If $l>0$ and $ms = n+1$, then $H^l(\mathbb{P}^n,\Omega^j_{\mathbb{P}^n}) =0$ provided that $j\neq l$. Therefore, $H^{n+r-s}(\mathbb{P}^r\times\mathbb{P}^n ,\Omega_{\mathbb{P}^r\times\mathbb{P}^n }^s(s-r-1,ms-n-1))$ is possibly nonzero if $l=j$, that is, $r = i+j -1= k+l -1= n-2$. In this case, it follows from $m(r+1) = ms = n+1 = r+3$ that $r=1$ and $m=2$, so $n = 3$. 
			\end{itemize}
			\item $k=0$ and $j> r+1$. Then $l=n+r-s>0$ and:
			\begin{itemize}  
				\item If $ms\neq n+1$, we have that $ms-n-1\geq j-n$ by (\ref{msr equation}). This implies that $H^l(\mathbb{P}^n,\Omega^j_{\mathbb{P}^n}(ms-n-1)) = 0$. 
				\item If $ms = n+1$, then $H^l(\mathbb{P}^n,\Omega^j_{\mathbb{P}^n}) = 0$ unless $l=j$. But if $l=j$ then $s = i+j \geq j = n+r-s$ implies that $2s\geq n+r \geq n+1$ while on the other hand $2s \leq ms = n+1$. We see from this that  $H^{n+r-s}(\mathbb{P}^r\times\mathbb{P}^n ,\Omega_{\mathbb{P}^r\times\mathbb{P}^n }^s(s-r-1,ms-n-1))$ is possibly nonzero if $m=2, r=1$ and $n$ is odd. 
			\end{itemize} 
		\end{enumerate}
		We conclude that the statement holds.  
	\end{proof}
	\begin{proof}[Proof of Proposition \ref{one dimensionality Jrho}]
		We have that $H^{n+r}(\mathbb{P}^r\times\mathbb{P}^n, \Omega_{\mathbb{P}^n\times\mathbb{P}^r}^{n+r})\cong k$ is the final cohomology group of the sequence (\ref{Koszul complex for J}), by Lemma \ref{spectral sequence lemma jrho onedimensional}. Note that 
		$$H^0(\mathbb{P}^r\times\mathbb{P}^n, \mathcal{O}(\rho)) = k[Y_0,\cdots, Y_r,X_0,\cdots, X_n]^\rho.$$
		The image of the map 
		$$\Lambda^{n+r}\Sigma_{\mathcal{L}}^\vee\otimes \mathcal{L}^{n+r}\otimes\Omega_{\mathbb{P}^n\times\mathbb{P}^r}^{n+r} \to \Lambda^{n+r+1}\Sigma_{\mathcal{L}}^\vee\otimes \mathcal{L}^{n+r+1}\otimes\Omega_{\mathbb{P}^n\times\mathbb{P}^r}^{n+r}$$ is the ideal $(F_0,\cdots, F_r, \bar{F}_0,\cdots, \bar{F}_n)$ so  $H^{n+r}(\mathbb{P}^r\times\mathbb{P}^n, \Omega_{\mathbb{P}^n\times\mathbb{P}^r}^{n+r}) \cong J^\rho$. In particular, $J^\rho$ is one dimensional.   
	\end{proof} 
	\begin{rem} 
		As a consequence of this statement, the map $\phi$ from Corollary~\ref{corollary definition map Jrho to cohomology} is an isomorphism. This proves Corollary \ref{Corollary Jrho ondimensional}.  
	\end{rem} 
	\begin{rem} 
		An interesting question is whether $\phi$ coincides with the map $J^\rho\to H^{n+r}(\mathbb{P}^r\times\mathbb{P}^n,\Omega_{\mathbb{P}^r\times\mathbb{P}^n}^{n+r}$ which we get from the above argument. We have not been able to give a full answer, but suspect it might be true. A possible way to make the map from the above argument explicit might be as follows. Take a \v{C}ech resolution of all terms in the resolution (\ref{Koszul complex for J}). If we start with an element $\xi\in \mathcal{C}^{n+r-1}(\mathcal{U},\Omega_{\mathbb{P}^r\times\mathbb{P}^n}^{n+r})$, then using that $H^{n+r-1}(\mathbb{P}^r\times\mathbb{P}^n,\Omega_{\mathbb{P}^r\times\mathbb{P}^n}^{n+r}) = 0$, the element comes from some $\xi'\in \mathcal{C}^{n+r-2}(\mathcal{U},\Omega_{\mathbb{P}^r\times\mathbb{P}^n}^{n+r})$. Computing $\xi'$ and applying the horizontal map, we find an element of $\mathcal{C}^{n+r-2}(\mathcal{U},\Sigma_{\mathcal{L}}^\vee\otimes\mathcal{L}^{-n-r-1}\otimes\Omega_{\mathbb{P}^r\times\mathbb{P}^n}^{n+r})$. Repeating the above procedure and keeping track of everything, we end up with an element in $\mathcal{C}^0(\mathcal{U},\mathcal{O}(\rho))$. After that, we can check whether the resulting morphism is the same as $\phi$. 
	\end{rem} 
	We now prove Proposition \ref{tildeJrho onedimensional}. 
	\begin{nota}
		We write $\mathcal{E} = \bigoplus_{i=1}^{n+r+1}\mathcal{O}(-1,-m)$.
	\end{nota}
	There is a surjective morphism $\mathcal{E}\to\mathcal{O}_{\mathbb{P}^r\times\mathbb{P}^n}$ sending local generators $g_i$ for $i\in\{1,\cdots, r\}$ to $g_iY_iF_i$ and $g_i\in\{r+1,\cdots, n+r+1\}$ to $g_iX_{i-r-1}\bar{F}_{i-r-1}$. This gives rise to an exact Koszul complex 
	\begin{equation}
	0\to \Lambda^{n+r+1}\mathcal{E}\to \cdots \to \Lambda^2\mathcal{E}\to \mathcal{E}\to\mathcal{O}_{\mathbb{P}^r\times\mathbb{P}^n}\to 0.
	\end{equation}	
	and applying $\text{Hom}(-,\Omega_{\mathbb{P}^n\times\mathbb{P}^r}^{n+r})$ we find 
	\begin{equation}\label{Koszul complex for tildeJ}
	0\to \Omega_{\mathbb{P}^n\times\mathbb{P}^r}^{n+r}\to \mathcal{E}^\vee\otimes\Omega_{\mathbb{P}^n\times\mathbb{P}^r}^{n+r}\to \cdots \to \Lambda^{n+r+1}\mathcal{E}^\vee\otimes\Omega_{\mathbb{P}^n\times\mathbb{P}^r}^{n+r}\to 0.
	\end{equation}
	Note that $\Lambda^{n+r+1}\mathcal{E}^\vee\otimes\Omega_{\mathbb{P}^n\times\mathbb{P}^r}^{n+r} = \mathcal{O}(\rho - (r+1,n+1))$. 
	\begin{lem}\label{resolution Jtilderho acyclic lemma}
		The complex (\ref{Koszul complex for tildeJ}) defines an acyclic resolution of $\Omega_{\mathbb{P}^n\times\mathbb{P}^r}^{n+r}$. 
	\end{lem}
	\begin{proof}
		We have that $\Lambda^i\mathcal{E}^\vee$ is a direct sum of terms $\mathcal{O}(i-r-1,mi-n-1)$. For $q>0$ we have that 
		\begin{align}
		H^q(\mathbb{P}^r\times\mathbb{P}^n, \mathcal{O}&(i-r-1, mi-n-1))\\
		&= \bigoplus_{a+b=q}H^a(\mathbb{P}^r, \mathcal{O}(i-r-1))\otimes H^b(\mathbb{P}^n, \mathcal{O}(mi-n-1)).\nonumber 
		\end{align}
		Note that $H^a(\mathbb{P}^r, \mathcal{O}(i-r-1)) = 0$ unless possibly if $a\in\{0,r\}$. If $a=0$, we have that $H^b(\mathbb{P}^n, \mathcal{O}(mi-n-1)) = 0$ unless possibly if $b=n$. But in that case
		$$H^n(\mathbb{P}^n, \mathcal{O}(mi-n-1)) = \left(\frac{1}{X_0\cdots X_n}k\left[\frac{1}{X_0},\cdots, \frac{1}{X_n}\right] \right)^{mi-n-1} = 0$$ because $mi-n-1 \geq -n$ (as $mi-n-1 < -n$ would imply $mi < 1$). If $a=r$ then we similarly see that $H^r(\mathbb{P}^r, \mathcal{O}(i-r-1)) = 0$ as $i-r-1\geq -r$. 
	\end{proof}
	\begin{proof}[Proof of Proposition \ref{tildeJrho onedimensional}]
		By Lemma \ref{resolution Jtilderho acyclic lemma} we have that $H^{n+r}(\mathbb{P}^r\times\mathbb{P}^n, \Omega_{\mathbb{P}^n\times\mathbb{P}^r}^{n+r})$ is the final cohomology group of the sequence (\ref{Koszul complex for tildeJ}). The image of the map 
		$$\Lambda^{n+r}\mathcal{E}^\vee\otimes\Omega_{\mathbb{P}^n\times\mathbb{P}^r}^{n+r} \to \Lambda^{n+r+1}\mathcal{E}^\vee\otimes\Omega_{\mathbb{P}^n\times\mathbb{P}^r}^{n+r}$$ is the ideal $(Y_1F_1,\cdots, Y_rF_r, X_0\bar{F}_0,\cdots, X_n\bar{F}_n)$. Because of the Euler relations, this is equal to the ideal $(Y_0F_0,\cdots, Y_rF_r, X_0\bar{F}_0,\cdots, X_n\bar{F}_n)$  so that we find that $k\cong H^{n+r}(\mathbb{P}^r\times\mathbb{P}^n, \Omega_{\mathbb{P}^n\times\mathbb{P}^r}^{n+r}) \cong \tilde{J}^{\rho + (r+1,n+1)}$. 
	\end{proof} 
	
	\section{Computing the trace map}\label{section computing trace}
	Again, we keep the notation which was set up in the previous sections. For $A\in J^{p-r, (p+1)m-n-1}$ and $B\in J^{q-r, (q+1)m-n-1}$, consider their images 
	$$\omega_A = \psi_p(A)\in H^p(\mathcal{X},\Omega_{\mathcal{X}}^q)\text{ and }\omega_B = \psi_q(B)\in H^q(\mathcal{X},\Omega_{\mathcal{X}}^p)$$ under the isomorphism from Proposition \ref{Levine Proposition 3.2}. We note that as the trace map is compatible with pushforwards, we have that 
	$$\text{Tr}_{\mathcal{X}}(\omega_A\cup\omega_B) =  \text{Tr}_{\mathbb{P}^r\times\mathbb{P}^n}(i_*(\omega_A\cup\omega_B)).$$ 
	Therefore, we can compute $\text{Tr}(\omega_A\cup \omega_B)$ by representing $i_*(\omega_A\cup\omega_B)$ on the open cover $\mathcal{U}$ using Proposition \ref{representation of image cech} and comparing it with a representation of $c_1(\mathcal{O}(1,m))^{n+r}$. We will do so on a refinement of $\mathcal{U}$, which we now first construct. 
	\begin{nota} \label{extra assumptions}
		In this section we will make the following three extra assumptions: 
		\begin{itemize}
			\item $m+1$ is invertible in $k$. 
			\item $V(F_0),\cdots, V(F_r)$ are smooth hypersurfaces in $\mathbb{P}^n$ and they intersect transversally: $V(F_{i_0},\cdots, F_{i_s})$ is a smooth closed subscheme of $\mathbb{P}^n$ which is of codimension $s+1$ for all $\{i_0,\cdots, i_s\}\subset \{0,\cdots, r\}$. 
			\item The first assumption remains true after setting any proper subset of the $X_i$'s equal to zero and replacing $\mathbb{P}^n$ with the linear subspace defined by the vanishing of the chosen $X_i$'s.
		\end{itemize}
	\end{nota} 
	Now let $$G_0 = Y_0F_0, \cdots, G_r = Y_rF_r, G_{r+1} = X_0\bar{F}_0,\cdots, G_{n+r+1} = X_n\bar{F}_n.$$ The following is true because of the extra assumptions made in Notation \ref{extra assumptions}. 
	\begin{lem}
		The set of opens $\mathcal{V} = \{V_0, \cdots, V_{n+r+1}\}$ where $V_i = \{G_i\neq 0\}$ is an open cover of $\mathbb{P}^r\times\mathbb{P}^n$. 
	\end{lem}
	\begin{proof}
		Suppose that $x\in\mathbb{P}^r\times\mathbb{P}^n$ is not in $V_0\cup\cdots\cup V_r$. Then we have that $Y_0F_0 = \cdots = Y_rF_r = 0$ at $x$. As not all $Y_i$ can be zero at $x$, there is some $F_i$ which is zero. As $V(F_i)$ is smooth, we have that $\partial F_i/\partial X_j$ is nonzero at $x$ for some~$j$. So $\bar{F}_j = \sum_{i=0}^rY_i\partial F_i/\partial X_j$ is nonzero at $x$. If $X_j$ is nonzero at $x$, then $x\in V_{j+r+1}$. Otherwise, as $V(F_i)$ remains smooth after intersecting with $X_j = 0$, we have that $\partial F_i/\partial X_{j'}$ is nonzero at $x$ for some $j'\neq j$. We now repeat the above argument until we find some open in $\mathcal{V}$ containing $x$. 
	\end{proof}
	\begin{prop}
		If we remove $V_j$ for any $j\in\{0,\cdots, n+r+1\}$, this still results in an open cover. 
	\end{prop}
	\begin{proof}
		We have the two Euler equations $\sum_{i=0}^rY_iF_i = F$ and $\sum_{i=0}^nX_i\bar{F}_i = mF$
		and so $m\sum_{i=0}^rG_i - \sum_{i=r+1}^{n+r+1}G_i = 0$. 
	\end{proof}
	\begin{nota} 
		We write $\mathcal{W}$ for the cover $\mathcal{V}$ with $V_0$ removed, i.e.
		$$\mathcal{W} = \{V_1,\cdots, V_{n+r+1}\}.$$ 
	\end{nota} 
	\begin{rem} \label{remark refinements}
		Note that $\mathcal{V}$ is a refinement of $\mathcal{U}$ via the identity map on the index sets. The inclusion 
		$\{1,\cdots, n+r+1\}\to \{0,\cdots, n+r+1\}$ makes $\mathcal{W}$ into a refinement of~$\mathcal{V}$. This yields a composition of refinement maps $$C^{n+r}(\mathcal{U}, \Omega_{\mathbb{P}^r\times\mathbb{P}^n}^{n+r})\to C^{n+r}(\mathcal{V}, \Omega_{\mathbb{P}^r\times\mathbb{P}^n}^{n+r})\to  C^{n+r}(\mathcal{W}, \Omega_{\mathbb{P}^r\times\mathbb{P}^n}^{n+r})$$
		sending a cocycle $\{s_0,\cdots, s_r, \bar{s}_0,\cdots, \bar{s}_n\}\in  C^{n+r}(\mathcal{U}, \Omega_{\mathbb{P}^r\times\mathbb{P}^n}^{n+r})$ - where $s_i$ lives on the intersection of all opens except for $U_i$ and $\bar{s}_j$ on the intersection of all opens except for $\bar{U}_j$ - to $\{s_0 \} \in C^{n+r}(\mathcal{W}, \Omega_{\mathbb{P}^r\times\mathbb{P}^n}^{n+r})$. 
	\end{rem} 
	\begin{nota}
		Consider the matrix $M$ given by
		\begin{equation}\label{matrix} \setlength\arraycolsep{2.5pt}
		\begin{pmatrix} 
		F_0 & 0 & \cdots & 0 & Y_0\frac{\partial F_0}{\partial X_0} & Y_0\frac{\partial F_0}{\partial X_1 } & \cdots & Y_0\frac{\partial F_0}{\partial X_n} \\
		0 & F_1 & \cdots & 0 & Y_1\frac{\partial F_1}{\partial X_0} & Y_1\frac{\partial F_1}{\partial X_1 } & \cdots & Y_1\frac{\partial F_1}{\partial X_n} \\
		\vdots & \vdots & \ddots & \vdots & \vdots & \vdots & \cdots & \vdots \\
		0 & 0 & \cdots & F_r & Y_r\frac{\partial F_r}{\partial X_0} & Y_r\frac{\partial F_r}{\partial X_1 } & \cdots & Y_r\frac{\partial F_r}{\partial X_n}\\
		X_0\frac{\partial F_0}{\partial X_0} & X_0\frac{\partial F_1}{\partial X_0} & \cdots & X_0\frac{\partial F_r}{\partial X_0} & \bar{F}_0 + X_0\frac{\partial \bar{F}_0}{\partial X_0} & X_0\frac{\partial \bar{F}_0}{\partial X_1} & \cdots & X_0\frac{\partial\bar{F}_0}{\partial X_n} \\
		X_1\frac{\partial F_0}{\partial X_1} & X_1\frac{\partial F_1}{\partial X_1} & \cdots & X_1\frac{\partial F_r}{\partial X_1} &  X_1\frac{\partial \bar{F}_1}{\partial X_0} & \bar{F}_1 + X_1\frac{\partial \bar{F}_1}{\partial X_1} & \cdots & X_1\frac{\partial\bar{F}_1}{\partial X_n} \\
		\vdots & \vdots & \cdots & \vdots & \vdots & \vdots & \ddots & \vdots \\
		X_n\frac{\partial F_0}{\partial X_n} & X_n\frac{\partial F_1}{\partial X_n}  & \cdots & X_n\frac{\partial F_r}{\partial X_n} &  X_n\frac{\partial \bar{F}_n}{\partial X_0} & X_n\frac{\partial \bar{F}_n}{\partial X_1} & \cdots & \bar{F}_n+ X_n\frac{\partial\bar{F}_n}{\partial X_n} 
		\end{pmatrix} \end{equation}
		Let $M_{i|j}$ be the minor with the $i$'th row and $j$'th column left out. Note that $Y_i$ divides $\det(M_{0|i})$ for $i>0$ and $X_j$ divides $\det(M_{0|j})$. 
	\end{nota}
	We will prove the following statement in Section \ref{subsection ctilde}. 
	\begin{lem}\label{construction of ctilde}
		There exists a unique $\tilde{C} \in k[Y_0,\cdots, Y_r, X_0,\cdots, X_n]^{\rho + (r+1,n+1)}$ such that
		$$(m+1)Y_iX_j\tilde{C} = (-1)^j\det(M_{0|j+r+1})Y_i + (-1)^{r+i}\det(M_{0|i})X_j$$
		for $i\in\{0,\cdots, r\}$ and $j\in\{0,\cdots, n\}$. We have that $c_1(\mathcal{O}(1,m))^{n+r}$ is represented by 
		$$\frac{\tilde{C}\omega\wedge\bar{\omega}}{\prod_{i=1}^rY_iF_i\prod_{j=0}^nX_j\bar{F}_j}\in C^{n+r}(\mathcal{W},\Omega^{n+r}_{\mathbb{P}^r\times\mathbb{P}^n}).$$
	\end{lem} 
	\begin{rem}
		One way to view the situation: we can embed $\mathbb{P}^r\times\mathbb{P}^n$ into $\mathbb{P}^{n+r+1}$ using coordinates $Y_0F_0,\cdots, Y_rF_r, X_0\bar{F}_0,\cdots, X_n\bar{F}_n$. This then lands in a hyperplane $H$. Computing the first Chern class boils down to pulling back the generator $\omega_H$ of $\Omega_H^{n+r}$ to $\mathbb{P}^r\times\mathbb{P}^n$. 
	\end{rem}
	We will then prove the following statement in Section \ref{subsection trace}. 
	\begin{theorem}\label{theorem trace}
		Assume that we are not in the situation that $\dim(\mathcal{X})$ is odd, $r=1$ and $m=2$. Then the map
		$$\psi: J^\rho\to \tilde{J}^{\rho + (r+1,n+1)}, D\mapsto D\prod_{i=0}^rY_i\prod_{j=0}^nX_j$$ is an isomorphism. Therefore, for the element $\tilde{C}$ from Lemma \ref{construction of ctilde}, we have that $\tilde{C} = \psi(C)$ for a unique $C\in J^\rho$. Write $AB = \lambda C$ in $J^\rho$ for some $\lambda\in k^*$. Then $$\text{Tr}(\omega_A\cup\omega_B) = (-1)^{r+1}m^{n+1}\binom{n+r}{r}\lambda.$$ 
	\end{theorem}	
	We can then find an analogue of the Scheja-Storch generator from the classical case. 
	\begin{defn}\label{definition SS generator}
		Suppose that $\binom{n+r}{r}$ is invertible in $k$. Then the trace one element $e_F = \frac{C}{m^n\binom{n+r}{r}}\in J^\rho$ is called the \textit{Scheja-Storch generator} of $X$.  
	\end{defn}
	\begin{rem}
		We conjecture that the assumption that $\binom{n+r}{r}$ is invertible in $k$ is not necessary, i.e. that one can find a similar construction of the Scheja-Storch generator as in Construction \ref{construction SS generator}. 
	\end{rem}
	\begin{rem}[Continuation of Remark \ref{remark r=0} and Remark \ref{remark r=0 products}]
		If we take $r=0$, the formula from Theorem \ref{theorem trace} becomes 
		$$\text{Tr}(\omega_A\cup\omega_B) = -m\lambda$$ for $A,B\in J$ such that $AB = \lambda e_F\in J^\rho$ for some $\lambda\in k^*$. This is in accordance with \cite[Theorem 3.8]{LevineECHWHH}. 
	\end{rem} 
	In characteristic zero, one has an explicit formula for the Scheja-Storch generator in the case where $r=0$, used in \cite{LevineECHWHH}. In Section \ref{subsection r=0}, we show that under the map $Y_0\mapsto 1$ from $J$ to the classical Jacobian ring, the Scheja-Storch element from Definition \ref{definition SS generator} maps to this Scheja-Storch element up to a factor~$\frac{m^{n+1}}{m+1}$. 
	
	\subsection{Proof of Lemma \ref{construction of ctilde}}\label{subsection ctilde}
	Before proving Lemma \ref{construction of ctilde}, we first prove a lemma in a slightly more general setting. Write $Z_0 = Y_0, \cdots, Z_r = Y_r, Z_{r+1} = X_0, \cdots, Z_{n+r+1} = X_n$ and let $G_0,\cdots, G_{n+r+1}$ be homogeneous polynomials in the $Z_i$ of the same total degree~$m+1$. Consider the matrix $M = (m_{ij})_{ij}$ with $m_{ij} = \frac{\partial G_i}{\partial Z_j}$.
	\begin{lem}
		Let $i,j\in\{0,\cdots, n+r+1\}$ be two distinct elements. We have that 
		$$dG^{i,j} = \sum_{k<l} \det(M_{ij|kl}) dZ^{k,l}$$ where $M_{ij|kl}$ is the minor of $M$ with the $i$'th and $j$'th row and the $k$'th and $l$'th column removed. 
	\end{lem}
	\begin{proof}
		We have that 
		\begin{align*} 
		dG^{i,j} &= \bigwedge_{k\neq i,j}\left(\sum_{l=0}^{n+r+1}\frac{\partial G_k}{\partial Z_l}dZ_l\right)\\
		&= \sum_{(j_1,\cdots j_{n+r})}\prod_{k\neq i,j} m_{kj_k} dZ_{j_1}\wedge\cdots\wedge dZ_{j_{n+r}}\\
		&= \sum_{k<l}\sum_{\sigma\in S_{n+r}}\text{sign}(\sigma)\prod_{p\neq i,j} m_{p\sigma(p)} dZ^{k,l}\\
		&= \sum_{k<l}\det(M_{i j| k l}) dZ^{k,l}
		\end{align*} 
		where $S_{n+r}$ is the symmetric group on $n+r$ elements and we used the Leibniz definition of a determinant. 
	\end{proof}
	Now for $k\in\{0,\cdots, n+r+1\}$ we write $$\tau_k = \sum_{i<k}(-1)^iZ_idZ^{i,k} + \sum_{i>k}(-1)^{i+1}Z_idZ^{k,i}.$$  
	Note that this definition is similar to the definition of generators $\tau_i$ of $\Omega_{\mathbb{P}^l}^{l-1}$ used in \cite{TerasomaIVH}, which was also used in the argument of Proposition \ref{Levine Proposition 3.2}. 
	\begin{lem}\label{lemma Cech representation}
		Let $i\in \{0,\cdots, n+r+1\}$. Then on the intersection of the opens $\{G_j\neq 0\}$ for $j\neq i$, we have that
		\begin{align*}
		d\log\left(\frac{G_1}{G_0}\right)\wedge & \cdots \wedge d\log\left(\frac{G_{i+1}}{G_{i-1}}\right)\wedge \cdots \wedge d\log\left(\frac{G_{n+r+1}}{G_{n+r}}\right) \\
		&= \frac{G_i}{(m+1)\prod_{j=0}^{n+r+1}G_j}\sum_{k=0}^{n+r+1}\det(M_{i|k})\tau_k
		\end{align*}
	\end{lem}
	\begin{proof}
		We can compute that 
		\begin{align*}
		d&\log\left(\frac{G_1}{G_0}\right)\wedge \cdots \wedge d\log\left(\frac{G_{i+1}}{G_{i-1}}\right)\wedge \cdots \wedge d\log\left(\frac{G_{n+r+1}}{G_{n+r}}\right)\\
		&= \left(\frac{dG_{i+1}}{G_{i+1}} - \frac{dG_{i-1}}{G_{i-1}}\right)\prod_{j\neq i, i-1}\left(\frac{dG_{j+1}}{G_{j+1}} - \frac{dG_j}{G_j}\right)\\
		&= \frac{G_i\left(\sum_{j < i}(-1)^jG_jdG^{j,i} + \sum_{j> i}(-1)^{j+1}G_jdG^{i,j}\right)}{\prod_{p=0}^{n+r+1}G_{p}}\\
		&= \frac{G_i\sum_{j < i}(-1)^jG_j\sum_{k<l}\det(M_{ij|kl})dZ^{k,l}}{\prod_{p=0}^{n+r+1}G_{p}} \\
		&\quad + \frac{G_i\sum_{j > i}(-1)^{j+1}G_j\sum_{k<l}\det(M_{ij|kl})dZ^{k,l}}{\prod_{p=0}^{n+r+1}G_{p}} 
		\end{align*}
		Now recall the Euler equation $(m+1)G_j = \sum_{p=0}^{n+r+1} Z_{p}\frac{\partial G_j}{\partial Z_{p}} = \sum_{p=0}^{n+r+1} Z_{p} m_{jp}$. We see that 
		\begin{align*}
		&(m+1)\left(\sum_{j < i}(-1)^jG_j\sum_{k<l}\det(M_{ij|kl})dZ^{k,l} + \sum_{j > i}(-1)^{j+1}G_j\sum_{k<l}\det(M_{ij|kl})dZ^{k,l}\right) \\
		&= \sum_{j < i}(-1)^j\sum_{p= 0}^{n+r+1}Z_{p}m_{jp}\sum_{k<l}\det(M_{ij|kl})dZ^{k,l} \\
		&\quad + \sum_{j > i}(-1)^{j+1}\sum_{p=0}^{n+r+1}Z_{p}m_{jp}\sum_{k<l}\det(M_{ij|kl})dZ^{k,l} \\
		&= \sum_p(-1)^pZ_p\big(\sum_{j < i}(-1)^{j+p}m_{jp} +  \sum_{j > i}(-1)^{j+p+1}m_{jp}\big) \sum_{k < l}\det(M_{ij|kl})dZ^{k,l} \\
		&= \sum_p(-1)^pZ_p\Big(\sum_{k<p}\big(\sum_{j < i}(-1)^{j+p}m_{jp} +  \sum_{j > i}(-1)^{j+p+1}m_{jp}\big)\det(M_{ij|kp})dZ^{k,p} \\
		&\quad + \sum_{p<l}\big(\sum_{j < i}(-1)^{j+p}m_{jp} +  \sum_{j > i}(-1)^{j+p+1}m_{jp}\big) \det(M_{ij|pl})dZ^{p,l}\\ 
		&\quad + \sum_{\substack{k<l\\ k,l\neq p}}\big(\sum_{j < i}(-1)^{j+p}m_{jp} +  \sum_{j > i}(-1)^{j+p+1}m_{jp}\big) \det(M_{ij|kl})dZ^{k,l}\Big) \\
		&= \sum_p(-1)^pZ_p\big(\sum_{k<p}\det(M_{i|k})dZ^{k,p} - \sum_{p>l}\det(M_{i|l})dZ^{p,l}\big) \\
		&= \sum_{k=0}^{n+r+1}\det(M_{i|k}) \big(\sum_{p<k}(-1)^pZ_pdZ^{p,k} + \sum_{p>k}(-1)^{p+1}Z_pdZ^{k,p}\big)\\
		&= \sum_k \tau_k \det(M_{i|k})	
		\end{align*} 
		as desired. Note that the seventh line is zero, as this is the determinant of $M_{i|l}$ with the $k$'th row removed and replaced by the $p$'th one. Also, note that the sixth line picks up an extra minus as we have to jump over an extra column when computing the determinant. 
	\end{proof}
	Now let $$G_0 = Y_0F_0, \cdots, G_r = Y_rF_r, G_{r+1} = X_0\bar{F}_0,\cdots, G_{n+r+1} = X_n\bar{F}_n$$ as in the introduction of this section. Then those are all of bidegree $(1,m)$ and the total degree is $m+1$. The matrix $M$ is exactly given by (\ref{matrix}).
	\begin{proof}[Proof of Lemma \ref{construction of ctilde}]
		We can represent $c_1(\mathcal{O}(1,m))\in H^1(\mathbb{P}^r\times\mathbb{P}^n,\Omega_{\mathbb{P}^r\times\mathbb{P}^n}^1)$ on the cover $\mathcal{W}$ by 
		$$\Bigl\{d\log\left(\frac{G_j}{G_i}\right)\Bigr\}_{i,j}\in C^1(\mathcal{W},\Omega_{\mathbb{P}^r\times\mathbb{P}^n}^1).$$ Taking the cup product $n+r$ times following the rules for a cup product on \v{C}ech cohomology, we see that $c_1(\mathcal{O}(1,m))^{n+r}$ is represented by 
		$$d\log\left(\frac{G_2}{G_1}\right)\wedge \cdots \wedge d\log\left(\frac{G_{n+r+1}}{G_{n+r}}\right) \in C^{n+r}(\mathcal{W},\Omega_{\mathbb{P}^r\times\mathbb{P}^n}^{n+r}).$$
		By Lemma \ref{lemma Cech representation} this is equal to 
		$$\frac{\sum_{k=0}^{n+r+1}\det(M_{0|k})\tau_k}{(m+1)\prod_{i=1}^rY_iF_i\prod_{j=0}^nX_j\bar{F}_j}.$$
		The numerator $\sum_{k=0}^{n+r+1}\det(M_{0|k})\tau_k$ is a global section of the twisted sheaf $\Omega^{n+r}_{\mathbb{P}^r\times \mathbb{P}^n}(n+r+1, m(n+r+1))$. Because $\Omega^{n+r}_{\mathbb{P}^r\times\mathbb{P}^n}(r+1,n+1)$ is a trivial line bundle which has $\omega\wedge\bar{\omega}$ as a global generator, there exists a unique rational function $\tilde{C}\in k[Y_0,\cdots, Y_r,X_0,\cdots, X_n]^{n, (n+r+1)m-n-1}$ such that 
		$$\frac{\sum_{k=0}^{n+r+1}\det(M_{0|k})\tau_k}{(m+1)\prod_{i=1}^rY_iF_i\prod_{j=0}^nX_j\bar{F}_j} = \frac{\tilde{C}\omega\wedge\bar{\omega}}{\prod_{i=1}^rY_iF_i\prod_{j=0}^nX_j\bar{F}_j}.$$
		In order to find $\tilde{C}$, consider the affine patch $\{Y_i\neq 0\}, \{X_j\neq 0\}$ with coordinates $y_k = \frac{Y_k}{Y_i}, x_k = \frac{X_k}{X_j}$. We see that 
		$$\tau_k = (-1)^rdY^k\wedge \bar{\omega} = \begin{cases}
		0 &\text{ if } k\neq i \\
		(-1)^{r+i}\omega\wedge\bar{\omega} &\text{ if } k = i
		\end{cases} $$
		for $k\in\{0,\cdots, r\}$ and 
		$$\tau_k = \omega\wedge dX^{k-r-1} = \begin{cases}
		0 &\text{ if } k\neq j + r + 1 \\
		(-1)^j\omega\wedge\bar{\omega} &\text{ if } k = j + r + 1
		\end{cases} $$
		for $k\in \{r+1,\cdots, n+r+1\}$. Therefore, $\sum_{k=0}^{n+r+1}\det(M_{0|k})\tau_k$ reduces to $(-1)^j\det(M_{0|j+r+1}) + (-1)^{r+i} \det(M_{0|i})\omega\wedge\bar{\omega}$. We have that $\tilde{C}$ is of bidegree $\rho + (r+1,n+1)$. Homogenizing again and comparing coefficients of $\tau_k$ in $\omega\wedge\bar{\omega}$, we get
		$$(m+1)Y_iX_j\tilde{C} = (-1)^j\det(M_{0|j+r+1})Y_i + (-1)^{r+i}\det(M_{0|i})X_j$$ as desired. 
	\end{proof}
	
	\subsection{Proof of Theorem \ref{theorem trace}}\label{subsection trace}
	\begin{proof}[Proof of Theorem \ref{theorem trace}]
		Note that $\psi$ is well defined, because an element of the Jacobian ideal $(F_0,\cdots, F_r,\bar{F}_0,\cdots, \bar{F}_n)$ will be mapped to zero. \\
		Let $\phi: J^\rho\to H^{n+r}(\mathbb{P}^r\times\mathbb{P}^n,\Omega_{\mathbb{P}^r \times \mathbb{P}^n}^{n+r})$ be the map from Corollary \ref{corollary definition map Jrho to cohomology}. 
		Composing the map $J^\rho\to \check{H}^{n+r}(\mathcal{U}, \Omega_{\mathbb{P}^r\times\mathbb{P}^n}^{n+r})$ that gives rise to $\rho$ with the refinement map from Remark \ref{remark refinements}, we find the morphism 
		$$\psi_J: J^\rho\to \check{H}^{n+r}(\mathcal{W}, \Omega_{\mathbb{P}^r\times\mathbb{P}^n}^{n+r}), D\mapsto \frac{(-1)^{r+1}mDY_0F_0\omega\wedge\bar{\omega}}{\prod_{i=0}^rF_i\prod_{j=0}^n\bar{F}_j}.$$
		Because $\mathcal{W}$ is affine, the \v{C}ech cohomology of this cover computes the usual cohomology. This implies that $\psi_J$ is surjective. Using Proposition \ref{one dimensionality Jrho}, we conclude that $\psi_J$ is an isomorphism. \\
		Now consider the morphism 
		\begin{align*} 
		k[Y_0,\cdots, Y_r,X_0,\cdots, X_n]^{\rho + (r+1,n+1)}&\to C^{n+r}(\mathcal{W},\Omega^{n+r}_{\mathbb{P}^r\times\mathbb{P}^n}),\\ D&\mapsto  \frac{(-1)^{r+1}mD\omega\wedge\bar{\omega}}{\prod_{i=1}^rY_iF_i\prod_{j=0}^nX_j\bar{F}_j}.
		\end{align*} 	
		Note that coboundaries on $\mathcal{W}$ are precisely coming from the ideal generated by the $G_i$. We therefore find an induced morphism $$\psi_{\tilde{J}}: \tilde{J}^{\rho + (r+1,n+1)}\to \check{H}^{n+r}(\mathcal{W},\Omega^{n+r}_{\mathbb{P}^r\times\mathbb{P}^n}).$$
		By Lemma \ref{construction of ctilde}, we have that $\tilde{C}$ maps to the nonzero element $c_1(\mathcal{O}(1,m))^{n+r}$, meaning that $\psi_{\tilde{J}}$ is surjective. Using Proposition \ref{tildeJrho onedimensional}, we see that $\psi_{\tilde{J}}$ is an isomorphism. \\ 
		We now have the commutative diagram 
		\[
		\begin{tikzcd}
		J^\rho \arrow[rr, "\psi"]\arrow[rd, swap, "\psi_J"] & & \tilde{J}^{\rho + (r+1,n+1)} \arrow[ld, "\psi_{\tilde{J}}"]\\
		& \check{H}^{n+r}(\mathcal{W}, \Omega^{n+r}_{\mathbb{P}^r\times\mathbb{P}^n}) & 
		\end{tikzcd} 
		\]
		From this, we see that $\psi$ has to be an isomorphism, which proves the first part of the statement.  \\
		Now using Proposition \ref{representation of image cech} and applying the refinement morphisms, we have that $i_*(\omega_A\cup\omega_B)$ is represented by 
		$$\frac{(-1)^{r+1}ABmY_0\omega\wedge\bar{\omega}}{ \prod_{i=1}^rF_i\prod_{j=0}^n\bar{F}_{j}}\in C^{n+1}(\mathcal{W}, \Omega_{\mathbb{P}^r\times\mathbb{P}^n}^{n+r}).$$
		By Lemma \ref{construction of ctilde}, we have that $c_1(\mathcal{O}(1,m))^{n+r}$ is represented by 
		$$\frac{\tilde{C}\omega\wedge\bar{\omega}}{\prod_{i=1}^rY_iF_i\prod_{j=0}^nX_j\bar{F}_j}\in C^{n+r}(\mathcal{W},\Omega^{n+r}_{\mathbb{P}^r\times\mathbb{P}^n}).$$
		As $\psi$ is an isomorphism, there exists a $C\in J^\rho$ such that $\tilde{C} = \psi(C)$, from which we see that $C$ maps to $c_1(\mathcal{O}(1,m))^{n+r}\in H^{n+r}(\mathbb{P}^r\times\mathbb{P}^n,\Omega_{\mathbb{P}^r\times\mathbb{P}^n}^{n+r})$. Now using Proposition \ref{one dimensionality Jrho}, we have that $AB = \lambda C$ for some $\lambda\in k$. Using that the trace of $c_1(\mathcal{O}(1,m))^{n+r}$ is equal to $\binom{n+r}{r}m^{n}$ we obtain the desired result. 
	\end{proof}
	
	\subsection{The Scheja-Storch generator in characteristic zero for $r=0$}\label{subsection r=0}
	\begin{nota} 
		Assume in this section that $\text{char}(k) =0$. 	
	\end{nota} 
	\begin{nota} 
		As in Remark \ref{remark r=0}, let $X = V(F)\subset\mathbb{P}^n$ be a smooth hypersurface, defined by a homogeneous polynomial $F\in k[X_0,\cdots, X_n]$ of degree $m$. Form the hypersurface $\mathcal{X} = V(Y_0F)\subset\mathbb{P}^0\times\mathbb{P}^n$ and let $F_i = \frac{\partial F}{\partial X_i}$ and $F_{ij} = \frac{\partial F_i}{\partial X_j}$ for $i,j\in\{0,\cdots, n\}$. Note that we have Euler equations 
		\begin{equation}\label{Euler equation r=0}
		(m-1)F_i = \sum_{j =0}^nX_jF_{ij}. 
		\end{equation}
	\end{nota} 
	The Jacobian ring of $\mathcal{X}$ is the bigraded ring
	$$J = k[Y_0,X_0,\cdots, X_n]/(F, Y_0F_0,\cdots, Y_0F_n).$$ Furthermore, we set
	$$\tilde{J} = k[Y_0,X_0,\cdots, X_n]/(Y_0F, Y_0X_0F_0,\cdots, Y_0X_nF_n).$$ 
	Let 
	$$J_X= k[X_0,\cdots, X_n]/(F_0,\cdots, F_n)$$ be the Jacobian ring of $X$ as defined in \cite{LevineECHWHH}, which is a usual graded ring. We have the map $f:J^{a,b}\to J_X^b, Y_0\mapsto 1, X_j\mapsto X_j$ for all $a,b\in\mathbb{Z}$. Let 
	$$e_F = \frac{\det(\text{Hess}(F))}{(m-1)^{n+1}}\in J_X^{(n+1)(m+2)}$$ be the classical Scheja-Storch element of $X$, used in the proof of \cite[Lemma 3.7]{LevineECHWHH}. In this section we will prove the following result. 
	\begin{prop} \label{proposition r=0}
		We have that $f(\frac{(m+1)C}{m^{n+1}}) = e_F$ in $J_X^{(n+1)(m+2)}$.  
	\end{prop} 
	Applying Theorem \ref{theorem trace}, we find the $(n+2)\times (n+2)$ matrix 
	$$M = \left(\begin{matrix} 
	F & Y_0F_0 & Y_0F_1 & \cdots & Y_0F_n \\
	X_0F_0 & Y_0(F_0 + X_0F_{00}) & Y_0X_0F_{01} & \cdots & Y_0X_0F_{0n} \\
	\vdots & \vdots & \vdots & \ddots & \vdots \\
	X_nF_n & Y_0X_nF_{n0} & Y_0X_nF_{n1} & \cdots & Y_0(F_n + X_nF_{nn})
	\end{matrix} \right)$$
	and we note that this matrix has rank $n+1$, as the sum of the last $n+1$ rows is $m$ times the first row. 
	\begin{lem}\label{lemma detM for r=0}
		We have that $\det(M_{0|0}) = Y_0^{n+1}\frac{m^{n+1}}{(m-1)^{n+1}}\det(\text{Hess}(F)) \prod_{i=0}^nX_i$. 
	\end{lem} 
	\begin{proof} 
		Let $\det(\text{Hess}(F))_{j_0,\cdots, j_r}$ be the minor where the rows and columns $(j_0,\cdots, j_r)$ have been removed. We claim that 
		\begin{align*} 
		Y_0^{n+1}\det(\text{Hess}(F))_{j_0,\cdots, j_r} \prod_{i\notin (j_0,\cdots, j_r)} X_i \prod_{i\in (j_0,\cdots, j_r)} F_i &= Y_0^{n+1}\frac{\det(\text{Hess}(F))}{(m-1)^{j+1}}\prod_{i=0}^nX_i.
		\end{align*} 
		Without loss of generality, we can assume that $(j_0,\cdots, j_r) = (k+1,\cdots, n)$ for some $k$. The proof of the claim proceeds by induction on $k$. For $k =n$, the result is clear. Now suppose that 
		\begin{align*} 
		Y_0^{n+1}\det(\text{Hess}(F))_{k+1,\cdots, n}\cdot  \prod_{i=0}^{k} X_i \cdot \prod_{i=k+1}^nF_i &= Y_0^{n+1}\frac{\det( \text{Hess}(F))}{(m-1)^{n-k}}\prod_{i=0}^nX_i
		\end{align*} 
		for some $k$. Denote $H = \text{Hess}(F)_{k+1,\cdots, n}$ and write $H_{i,j}$ for the minor of $H$ with the $i$'th row and $j$'th column removed. Note that 
		$$0=(m-1)Y_0X_jF_j = Y_0\sum_{l=0}^nX_jX_lF_{lj}$$ in $\tilde{J}$  and so
		\begin{align*}
		((-1)^{k+j}X_k&\det(H_{j,k}) - X_j\det(H_{k,k}))Y_0 \prod_{i=0}^{k-1}X_i\\
		&= \left(-\sum_{i=0}^k(-1)^{i+j} (X_kF_{ki} + X_jF_{ji}) \det(H_{kj,ki})\right) Y_0\prod_{i=0}^{k-1}X_i \\
		&= \left(-\sum_{i=0}^k(-1)^{i+j} (X_kF_{ki} - X_kF_{ki} - \sum_{l\neq j,k} X_lF_{li}) \det(H_{kj,ki})\right)Y_0 \prod_{i=0}^{k-1}X_i \\
		&= \left(\sum_{i=0}^k(-1)^{i+j}\sum_{l\neq j,k} X_lF_{li}\det(H_{kj,ki})\right)Y_0\prod_{i=0}^{k-1}X_i  \\
		&= 0
		\end{align*}
		as the second sum on the fourth line is the determinant of $H_{k,k}$ with the $j$'th row replaced by the $i$'th row, which is zero. \\
		We now have that 
		\begin{align*}
		Y_0^{n+1}\frac{\det(\text{Hess}(F))}{(m-1)^{n-k}}\prod_{i=0}^nX_i &= \det(\text{Hess}(F))_{k+1,\cdots, n}\  \prod_{i=0}^{k} X_i \prod_{i=k+1}^nF_i\\ &=  Y_0^{n+1}\left(\sum_{j=0}^k(-1)^{k+j}F_{kj}\det(H_{k,j})\right)\prod_{i=0}^{k} X_i\prod_{i=k+1}^nF_i \\
		&=   Y_0^{n+1}\left(\sum_{j=0}^k X_jF_{kj}\det(H_{k,k})\right)\prod_{i=0}^{k-1} X_i \prod_{i=k+1}^nF_i\\
		&= (m-1)Y_0^{n+1}\det(\text{Hess}(F))_{k,\cdots, n}\prod_{i=0}^{k-1} X_i \prod_{i=k}^nF_i 
		\end{align*}
		This completes the proof of the claim. \\
		Let $\tilde{M}_i$ be the matrix given by
		$$\left(\begin{matrix} \setlength\arraycolsep{2pt}
		X_0F_{00} & \cdots & X_0F_{0,i-1} & X_0F_{0,i+1} & \cdots & X_0F_{0n} \\
		\vdots & \ddots & \vdots & \vdots & \ddots & \vdots \\
		X_{i-1}F_{i-1,i-1} & \cdots & X_{i-1}F_{i-1,i-1} &  X_{i-1}F_{i-1,i+1} & \cdots & X_{i-1}F_{i-1,n} \\
		X_{i+1}F_{i+1,0} & \cdots &  X_{i+1}F_{i+1,i-1} & F_{i+1} + X_{i+1}F_{i+1,i+1} & \cdots & X_nF_{i+1,n}\\
		\vdots & \ddots & \vdots & \vdots & \ddots & \vdots \\
		X_nF_{n0} & \cdots & X_nF_{n,i-1} & X_nF_{n,i+1} & \cdots & F_n+X_nF_{nn}
		\end{matrix} \right) $$ 
		We find that 
		\begin{align*} 
		\det(M_{0|0}) &= Y_0^{n+1}\det\left(\begin{matrix} 
		X_0F_{00} & \cdots & X_0F_{0n} \\
		\vdots & \ddots & \vdots \\
		X_nF_{n0} & \cdots & X_nF_{nn}
		\end{matrix} \right)\\ 
		&\quad + Y_0^{n+1}\sum_{i=0}^n F_i\det(\tilde{M}_i)\\
		&= Y_0^{n+1}\det(\text{Hess}(F))\prod_{j=0}^nX_j \\
		&\quad + Y_0^{n+1}\sum_{(j_0,\cdots, j_r)}\det(\text{Hess}(F))_{j_0,\cdots, j_r}\prod_{i\notin (j_0,\cdots, j_r)} X_i\prod_{i\in (j_0,\cdots, j_r)}F_i \\
		&= Y_0^{n+1}\left(\frac{1}{(m-1)^{n+1}}\sum_{j=0}^{n+1} \binom{n+1}{j}(m-1)^{n+1-j}\right)\det(\text{Hess}(F)) \prod_{i=0}^nX_i\\
		&= Y_0^{n+1}\frac{m^{n+1}}{(m-1)^{n+1}}\det(\text{Hess}(F)) \prod_{i=0}^nX_i
		\end{align*} 
		as desired. 
	\end{proof} 
	\begin{proof}[Proof of Proposition \ref{proposition r=0}]
		By expanding to the first column, one can show that $\det(M_{0|n+2}) = 0$, as $Y_0X_iF_i=0$ in $\tilde{J}$. Therefore, using Lemma \ref{lemma detM for r=0} and Lemma \ref{construction of ctilde} we have that $$(m+1)\tilde{C} = Y_0^{n+1}\frac{m^{n+1}}{(m-1)^{n+1}}\prod_{i=0}^nX_i\det(\text{Hess}(F))$$ and so using Theorem \ref{theorem trace}, we find that $$C = Y_0^{n}\frac{m^{n+1}}{(m+1)(m-1)^{n+1}}\det(\text{Hess}(F))\in J^\rho.$$ This implies that $f(C) = \frac{m^{n+1}}{m+1}e_F$ in $J_X^{(n+1)(m+2)}$. 
	\end{proof} 
	
	\section{Example: intersecting two generalized Fermat hypersurfaces of the same degree}\label{section Fermat}
	To see an application of Theorem \ref{theorem trace}, we compute the quadratic Euler characteristic of a complete intersection of two generalized Fermat hypersurfaces of the same degree. 
	\begin{nota} 
		Let $m\geq 2$ be coprime to $\text{char}(k)$, assume $m+1$ is invertible in $k$ and let $F_0 = \sum_{i=0}^na_iX_i^m$ and $F_1 = \sum_{i=0}^nb_iX_i^m$. Let $X = V(F_0,F_1)$ be their complete intersection. Furthermore, assume that $a_ib_j - a_jb_i\neq 0$ for all $j\neq i$. Write $L_i = (a_iY_0 + b_iY_1)$. Then $V(F_0)$ and $V(F_1)$ are both smooth, and so is $X$, and these conditions still hold when we set any subset of the $X_i$ equal to zero. We have that $$F = Y_0F_0 + Y_1F_1 = \sum_{i=0}^n(a_iY_0 + b_iY_1)X_i^m = \sum_{i=0}^nL_iX_i^m.$$ Again, we write $\mathcal{X} = V(F)$. 
	\end{nota} 
	We will prove the following result.
	\begin{prop}\label{qec of L}
		Define 
		$$A_{n,m} = \begin{cases}
		\frac{1}{2}\deg(c_n(T_{\mathcal{X}})) &\text{ if $n$ or $m$ odd} \\
		\frac{1}{2}\deg(c_n(T_{\mathcal{X}})) - n - 1 &\text{ if } n,m \text{ even} 
		\end{cases} $$
		The quadratic Euler characteristic of $\mathcal{X}$ is equal to 
		$$\chi(\mathcal{X}/k) = \begin{cases}
		A_{n,m}H & \text{ if $n$ or $m$ odd} \\
		A_{n,m}H + \sum_{k=0}^n \langle \prod_{i=0, i\neq k}^n (a_kb_i - a_ib_k) \rangle &\text{ if } n,m \text{ even} 
		\end{cases}$$
	\end{prop}
	This will imply the following.
	\begin{cor}\label{qec of X}
		Define 
		$$B_{n,m} = \begin{cases}
		\frac{1}{2}\deg(c_{n-2}(T_{X})) &\text{ if $n$ odd} \\
		\frac{1}{2}\deg(c_{n-2}(T_{X})) - 1 &\text{ if $n$ even, $m$ odd} \\
		\frac{1}{2}\deg(c_{n-2}(T_{X})) - n - 1 &\text{ if } n,m \text{ even} 
		\end{cases} $$
		The quadratic Euler characteristic of $X$ is equal to 
		$$\chi(X/k) = \begin{cases}
		B_{n,m}H &\text{ if $n$ odd} \\
		B_{n,m}H + \langle 1 \rangle &\text{ if $n$ even, $m$ odd} \\
		B_{n,m}H + \langle 1 \rangle +  \sum_{k=0}^n  \langle \prod_{i=0, i\neq k}^n (a_kb_i - a_ib_k) \rangle &\text{ if } n,m \text{ even} 
		\end{cases}$$
	\end{cor}
	
	\subsection{The case where $n=2$}
	The case where $n=2$ is special, so we treat that argument here first. In this case, $X = V(F_0,F_1)$ is the intersection of two Fermat curves $V(F_0)$ and $V(F_1)$ in~$\mathbb{P}^2$ with 
	$$F_0 = a_0X_0^m + a_1X_1^m + a_2X_2^m$$
	and
	$$F_1 = b_0X_0^m + b_1X_1^m + b_2X_2^m$$
	where the $a_i,b_i\in k^*$ satisfy $a_ib_j - a_jb_i \neq 0$ for all $i\neq j$. In order to calculate the corresponding quadratic Euler characteristic, we will need that for a separable field extension $k\subset L$, the natural map $\pi: \text{Spec}(L)\to \text{Spec}(k)$ induces a morphism $\pi_*: \text{GW}(L)\to \text{GW}(k)$ where for a form $\langle u \rangle\in \text{GW}(L)$, we have that $\pi_*\langle u \rangle$ is given by the composition 
	$$L\times L \xrightarrow{\langle u \rangle} L \xrightarrow{\text{Tr}_{L/k}} k.$$
	By \cite[Theorem 1.9]{HoyoisQRGLVTF} we have that $\chi(\text{Spec}(L)/k) = \pi_*(\langle 1 \rangle )$. The following is a standard fact about quadratic forms, but we include a proof for the sake of completeness. 
	\begin{lem}\label{trace of extension by mth root}
		Let $K$ be a perfect field of characteristic coprime to $2m$ and let $a\in K^*$. Consider the field extension $K(\alpha) = K[X]/(X^m + a)$ of $K$ and let~$u\in K(\alpha)^*$ be a unit. Then 
		$$\text{Tr}_{K(\alpha)/K}(\langle u\rangle) = \begin{cases}
		\frac{m-1}{2}H + \langle um\rangle &\text{ if~$m$ is odd} \\
		\frac{m-2}{2}H + \langle um\rangle + \langle -aum \rangle &\text{ if~$m$ is even}
		\end{cases}$$
	\end{lem}
	\begin{proof}
		Note that $K(\alpha)$ has the basis $1,\alpha,\alpha^2,\cdots, \alpha^{m-1}$ over $K$. 
		We have that 
		$$\text{Tr}_{K(\alpha)/K}(u\alpha^{i+j}) = \begin{cases}
		um &\text{ if } i = j = 0\\
		- aum &\text{ if } i+j = m\\
		0 &\text{ otherwise}
		\end{cases} $$
		Namely, if $i+j = m$, the multiplication by $u\alpha^m = -au$ corresponds to the diagonal matrix with $-au$ as its entries, and this has trace $-aum$. If $i=j=0$, multiplication by $u$ is the diagonal matrix with $u$ on the diagonal, which has trace $um$. If $i$ or $j$ is not zero and $i+j\neq m$, we are taking the trace of the matrix 
		\[
		\left(\begin{matrix}
		0 & 0 & \cdots & 0 & -au & 0 & \cdots & 0 \\
		0 & 0 &  \cdots & 0 & 0 & -au & \cdots & 0 \\
		\vdots & \vdots & \cdots & \vdots & \vdots & \vdots & \ddots & \vdots \\
		0 & 0 &  \cdots & 0 & 0 & 0 & \cdots & -au\\
		u & 0 &  \cdots & 0 & 0 & 0 & \cdots & 0\\
		0 & u &  \cdots & 0 & 0 & 0 & \cdots & 0\\
		\vdots & \vdots & \ddots & \vdots & \vdots & \vdots & \vdots & \vdots \\
		0 & 0 &  \cdots & u & 0 & 0 & \cdots & 0\\
		\end{matrix}
		\right)
		\] 
		which has trace zero. \\
		Therefore, the quadratic form $ \text{Tr}_{K(\alpha)/K}(\langle u\rangle) $ corresponds to the symmetric bilinear form with matrix 
		\[\left(
		\begin{matrix}
		um & 0 & 0 & \cdots & 0 & -uma \\
		0 & 0 & 0 & \cdots & -uma & 0 \\
		\vdots & \vdots & \vdots & \ddots & \vdots & \vdots \\
		0 & -uma & 0 & \cdots & 0 & 0\\
		-uma & 0 & 0 & \cdots & 0 & 0 
		\end{matrix}\right)
		\]
		which gives the form from the statement. 
	\end{proof}
	\begin{prop} 
		The quadratic Euler characteristic of~$X$ equals 
		$$\chi(X/k) = \begin{cases}
		\frac{(m+1)(m-1)}{2}H + \langle 1 \rangle  &\text{ if~$m$ is odd} \\
		\frac{(m+2)(m-2)}{2}H + \langle 1 \rangle + \sum_{i=0}^2\langle \prod_{j\neq i}(a_ib_j-a_jb_i)\rangle  &\text{ if~$m$ is even}
		\end{cases}$$
	\end{prop} 
	\begin{proof} 
		Without loss of generality, we can assume that~$X = V(F_0,F_1)$ lies inside the affine patch where~$X_2\neq 0$; otherwise, we change coordinates. Choosing coordinates~$x = \frac{X_0}{X_2}$ and~$y = \frac{X_1}{X_2}$ on~$\mathbb{A}^2$, we have that~$X$ is the zero set of the ideal~$$(a_0x^m + a_1y^m + a_2, b_0x^m + b_1y^m + b_2).$$ 
		Let~$K$ be the residue field of~$X$, that is: 
		$$K = k[x,y]/(a_0x^m + a_1y^m + a_2, b_0x^m + b_1y^m + b_2).$$
		Define~$$e = \frac{a_0b_2 - a_2b_0}{a_1b_0 - a_0b_1}\text{ and } f = \frac{a_1b_2 - a_2b_1}{a_0b_1 - a_1b_0}.$$ 
		Then we can view the extension~$k\subset K$ as one which takes place in two steps: 
		$$k\subset k(\alpha) = k[t]/(t^m + e) \subset K = k(\alpha)[s]/(s^m + f).$$
		Indeed, the system of equations 
		$$a_0x^m + a_1y^m + a_2 = 0 \text{ and } b_0x^m + b_1y^m + b_2 = 0$$ implies that 
		$$(a_1b_0 -a_0b_1)y^m + a_2b_0 - a_0b_2 = 0\text{ and } (a_0b_1 - a_1b_0)x^m + a_2b_1 - a_1b_2 = 0.$$
		We see from Lemma \ref{trace of extension by mth root} that for odd~$m$, we have that
		\begin{align*}
		\text{Tr}_{K/k}(\langle 1 \rangle) &= \text{Tr}_{k(\alpha)/k}(\text{Tr}_{K/k(\alpha)}\langle 1 \rangle) \\
		&= \text{Tr}_{k(\alpha)/k}\left(\frac{m-1}{2}H + \langle m\rangle \right)\\
		&= \frac{m(m-1)}{2}H + \frac{m-1}{2}H + \langle m^2\rangle  \\
		&= \frac{(m+1)(m-1)}{2}H + \langle 1 \rangle 
		\end{align*}
		and for even~$m$ we compute 
		\begin{align*}
		\text{Tr}_{K/k}(\langle 1 \rangle) &= \text{Tr}_{k(\alpha)/k}(\text{Tr}_{K/k(\alpha)}\langle 1 \rangle) \\
		&= \text{Tr}_{k(\alpha)/k}\left(\frac{m-2}{2}H + \langle m\rangle + \langle -fm \rangle \right) \\
		&= \frac{m(m-2)}{2}H + \langle m^2\rangle + \langle -m^2e\rangle + \langle -m^2f\rangle + \langle m^2ef\rangle + (m-2)H \\
		&= \frac{(m+2)(m-2)}{2}H  + \langle 1 \rangle + \langle -e \rangle + \langle -f \rangle + \langle ef \rangle \\
		&= \frac{(m+2)(m-2)}{2}H  + \langle 1 \rangle + \langle (a_0b_1-a_1b_0)(a_0b_2 - a_2b_0) \rangle \\
		&\quad + \langle (a_1b_0-a_0b_1)(a_1b_2 - a_2b_1) \rangle + \langle (a_2b_0 - a_0b_2)(a_2b_1 - a_1b_2) \rangle 
		\end{align*}
		which is the desired result.
	\end{proof} 
	
	\subsection{The Jacobian ring}\label{subsection J}
	In this situation, we can give a very explicit proof of the one dimensionality of the bidegree $\rho$ part of the Jacobian ring 
	$$J = k[Y_0,Y_1,X_0,\cdots, X_n]/(F_0,F_1,mL_0X_0^{m-1},\cdots, mL_nX_n^{m-1})$$
	and also give generators and understand their relations. This is following \cite[Section 4 and Section 5.1]{TerasomaIVH}.  
	\begin{prop}\label{Lj in terms of Li Lk}
		Let $i,j,k\in\{0,\cdots, n\}$ be distinct. We can write $L_i$ as a linear combination of $L_j$ and $L_k$, more precisely, we have that 
		$$L_i = \frac{a_kb_i - a_ib_k}{a_kb_j - a_jb_k}L_j +  \frac{a_ib_j - b_ia_j}{a_kb_j - a_jb_k}L_k.$$
	\end{prop}
	\begin{proof}
		The expression $L_i = aL_j + bL_k$ leads to the system of equations 
		\begin{align*}
		aa_j + ba_k &= a_i\\
		ab_j + bb_k &= b_i
		\end{align*}
		These imply that $ba_kb_j - bb_ka_j = a_ib_j - b_ia_j$ and so $b = \frac{a_ib_j - b_ia_j}{a_kb_j - a_jb_k}$ implying that $$a = a_j^{-1}(a_i - a_kb) = \frac{a_ia_kb_j - a_ja_ib_k - a_ka_ib_j + a_kb_ia_j}{a_j(a_kb_j - a_jb_k)} = \frac{a_kb_i - a_ib_k}{a_kb_j - a_jb_k}.$$ This proves the statement. 
	\end{proof} 
	\begin{nota} 
		Let $k,l\in \{0,\cdots, n\}$ be distinct. 
	\end{nota} 
	\begin{prop}\label{generators of the Jacobian ring}
		The graded piece $J^\rho$ is generated by the elements $$A_j = X_j^{m}\cdot X_0^{m-2}\cdots X_n^{m-2}\prod_{i\neq j,k,l}(a_iY_0 + b_iY_1)$$ for $j\in\{0,\cdots, n\}\setminus\{k,l\}$. 
	\end{prop}
	The statement will follow from two lemmas. 
	\begin{lem}\label{generators Jacobian ring first lemma}
		Consider a term $$A = X_0^{i_0}X_1^{i_1}\cdots X_n^{i_n}(a_{j_1}Y_0 + b_{j_1}Y_1)\cdots (a_{j_{n-2}}Y_0 + b_{j_{n-2}}Y_1)$$ where $i_0 + \cdots + i_n = (n+1)(m-2) + m$. If $i_k, i_l\geq m-1$ for $k,l\in \{0,\cdots, n\}$ distinct, then $A=0$. 
	\end{lem}
	\begin{proof}
		Assume without loss of generality that $i_0,i_1\geq m-1$. By Proposition~\ref{Lj in terms of Li Lk}, we can write any $L_i$ for $i\geq 2$ as a linear combination of $L_0$ and $L_1$. This implies that $A$ can be written as a linear combination of terms of the form $cX_0^{i_0}X_1^{i_1}\cdots X_n^{i_n}L_0^pL_1^q$ where $p+q = n-2$ and $c\in k$ is some constant. But $$(a_0Y_0 + b_0Y_1)X_0^{m-1} = 0 \text{ and } (a_1Y_0 + b_1Y_1)X_1^{m-1} =0$$ in $J$ and so $A=0$. 
	\end{proof}
	\begin{lem}\label{generators Jacobian ring second lemma}
		Let $A$ be as in Lemma \ref{generators Jacobian ring first lemma}. Then $\max_{k=0,\cdots n}i_k \leq 2m-2$.  
	\end{lem}
	\begin{proof}
		Suppose, without loss of generality, that $i_0 > 2m-2$. We have that $X_0^m = -\frac{1}{a_0}(a_1X_1^m + \cdots + a_nX_n^m)$ and so 
		\begin{align*} 
		A &= -\frac{1}{a_0}(a_1X_0^{i_0-m}X_1^{i_1+m}\cdots X_n^{i_n}(a_{j_1}Y_0 + b_{j_1}Y_1)\cdots (a_{j_{n-2}}Y_0 + b_{j_{n-2}}Y_1) + \\
		&\quad \cdots + a_nX_0^{i_0 - m}X_1^{i_1}\cdots X_n^{i_n + m}(a_{j_1}Y_0 + b_{j_1}Y_1)\cdots (a_{j_{n-2}}Y_0 + b_{j_{n-2}}Y_1))
		\end{align*} 
		Note that $i_0 - m \geq 2m-1-m = m-1$ and $i_k + m \geq m-1$ for all $k\in\{1,\cdots, n\}$. By Lemma \ref{generators Jacobian ring first lemma} this implies that $A=0$. 
	\end{proof}
	\begin{proof}[Proof of Proposition \ref{generators of the Jacobian ring}]
		First, note that $J^\rho$ is generated by all the terms~$A$ as in Lemma \ref{generators Jacobian ring first lemma}. Now take such a term $A$ and assume that it is nonzero. Then by Lemma~\ref{generators Jacobian ring first lemma}, there can only be one $j\in\{0,\cdots, n\}$ such that $i_j\geq m-1$, and by Lemma \ref{generators Jacobian ring second lemma}, all $i_l$ are smaller than $2m-2$. But as $$i_0 + \cdots + i_n = (n+1)(m-2) + m$$ the only possible way in which this can happen is if $i_j = 2m-2$ and all other $i_l$ are equal to $m-2$. Therefore, we can generate $J^\rho$ by all terms of the form 
		$$X_j^m\cdot X_0^{m-2}\cdots X_n^{m-2} L_{j_1}\cdots L_{j_{n-2}}.$$
		By Proposition \ref{Lj in terms of Li Lk}, we can choose the generators so that $j_i\notin \{j,k,l\}$ for all $i\in\{1,\cdots, n-2\}$. Furthermore, we can choose all $j_i$ to be distinct as the total degree has to be $n-2$. Finally, we can exclude the terms $A_k$ and $A_l$ using the relations $F_0 = 0$ and $F_1 = 0$. 
	\end{proof} 
	\begin{lem}[See \cite{TerasomaIVH}, Lemma 4.9]\label{new jacobian ring lemma 2}
		Let $p,q,r\in \{0,\cdots, n\}$ be distinct. Then $$X_q^{m}\prod_{i\neq p,q,r}L_i = - \frac{a_pb_r-a_rb_p}{a_pb_q - b_pa_q}X_r^m\prod_{i\neq p,q,r}L_i$$ in $J$. 
	\end{lem}
	\begin{proof} 
		We first note that 
		$$\sum_{i=0}^n (a_pb_i - a_ib_p)X_i^m = a_pF_1 - b_pF_0 = 0$$ and so multiplying by $\prod_{i\neq p,q,r}L_i$ we see that 
		$$(a_pb_q - a_qb_p)X_q^m\prod_{i\neq p,q,r}L_i + (a_pb_r - a_rb_p)X_r^m\prod_{i\neq p,q,r}L_i =0$$
		as desired. 
	\end{proof}
	\begin{cor}\label{relation between Aks}
		Let $j,j'\in\{1,\cdots, n-1\}$ be distinct. We have that $$A_j =  \frac{(a_{j'}b_k - a_{k}b_{j'})(a_lb_{j'} - a_{j'}b_l)}{(a_jb_k - a_kb_j)(a_lb_j - a_jb_l)}A_{j'}$$ in $J^\rho$. In particular, $J^\rho$ is one dimensional. 
	\end{cor} 
	\begin{proof}
		Using Proposition \ref{Lj in terms of Li Lk} we see that
		$$L_{j'} = \frac{a_{j'}b_j - a_jb_{j'}}{a_jb_k - a_kb_j}L_k + \frac{a_kb_{j'} - a_{j'}b_k}{a_jb_k - a_kb_j}L_j$$
		so using Lemma \ref{new jacobian ring lemma 2} we have that 
		\begin{align*} 
		A_j &= X_j^{m}\cdot X_0^{m-2}\cdots X_n^{m-2}\prod_{i\neq j,k,l}L_i \\
		&= \frac{a_{j'}b_j - a_{j}b_{j'}}{a_jb_k - a_kb_j}X_j^{m}\cdot X_0^{m-2}\cdots X_n^{m-2}\prod_{i\neq {j'},j,l}L_i \\
		&= -\frac{(a_{j'}b_j - a_{j}b_{j'})(a_lb_{j'} - a_{j'}b_l)}{(a_jb_k - a_kb_j)(a_lb_j - a_jb_l)}X_{j'}^{m}\cdot X_0^{m-2}\cdots X_n^{m-2}\prod_{i\neq {j'},j,l}L_i \\
		&= \frac{(a_{j'}b_k - a_{k}b_{j'})(a_lb_{j'} - a_{j'}b_l)}{(a_jb_k - a_kb_j)(a_lb_j - a_jb_l)}X_{j'}^{m}\cdot X_0^{m-2}\cdots X_n^{m-2}\prod_{i\neq {j'},k,l}L_i 
		\end{align*} 
		as desired. 
	\end{proof}
	
	\subsection{Computing the trace of multiples of the generators}
	Moving to the setting of Theorem \ref{theorem trace}, we note that in this case, we have that 
	$$M = \left(\begin{matrix} 
	F_0 & 0 & ma_0Y_0X_0^{m-1} & ma_1Y_0X_1^{m-1} & \cdots & ma_nY_0X_n^{m-1} \\
	0 & F_1 & mb_0Y_1X_0^{m-1} & mb_1Y_1X_1^{m-1} & \cdots & mb_nY_1X_n^{m-1} \\
	ma_0X_0^m & mb_0X_0^m & m^2L_0X_0^{m-1} & 0 & \cdots & 0\\
	ma_1X_1^m & mb_1X_1^m & 0 & m^2L_1X_1^{m-1} & \cdots & 0\\
	\vdots & \vdots & \vdots & \vdots & \ddots & \vdots \\
	ma_nX_n^m & mb_nX_n^m & 0 & 0 & \cdots & m^2L_nX_n^{m-1}
	\end{matrix} \right) $$
	As in Section \ref{subsection J}, let $k,l\in\{0,\cdots, n\}$ be distinct. 
	\begin{lem}\label{trace of Aj}
		Let $A_j = X_j^m\cdot X_0^{m-2}\cdots X_n^{m-2}\prod_{i\neq j,k,l}(a_iY_0+b_iY_1)$ be a generator as in Proposition \ref{generators of the Jacobian ring} and assume that $n$ is even. Let $A,B\in J$. If $AB = \lambda A_j$ in $J^\rho$ for some $\lambda\in k^*$, we have that
		$$\text{Tr}(\omega_{A}\cup\omega_B) = m^{3n+2}(n+1)^2(a_jb_k - a_kb_j)(a_jb_l-a_lb_j)\lambda.$$
	\end{lem}
	\begin{proof}
		One can show that
		$$\det(M_{0|1}) = -m^{2(n+1)}Y_1\left(\prod_{i=0}^nL_iX_i^{m-1}\right)\sum_{i=0}^n\frac{a_ib_iX_i^m}{L_i}$$
		and check that 
		$$\det(M_{0|n+2}) = (-1)^{n+1}m^{2n + 1}X_n\left(\prod_{i=0}^nL_iX_i^{m-1}\right)\sum_{i=0}^n\frac{a_ib_iX_i^m}{L_i}.$$
		It follows that
		$$(m+1)\tilde{C} = -(m+1)m^{2n+1}\left(\prod_{i=0}^nL_iX_i^{m-1}\right)\sum_{i=0}^n\frac{a_ib_iX_i^m}{L_i}.$$
		We note that for all $i$ we have
		\begin{align*}
		a_ib_iX_i^m & X_0^{m-2}\cdots X_n^{m-2}\prod_{p\neq i}L_p \\
		&= A_i\cdot a_ib_iL_kL_l\\
		&= A_i\cdot a_ib_i(a_ka_lY_0^2 + (a_kb_l + a_lb_k)Y_0Y_1 + b_kb_lY_1^2) \\
		&= A_i \cdot (a_ia_ka_lb_iY_0^2 + (-(a_ib_k - a_kb_i)(a_ib_l - a_lb_i)\\
		&\quad  + a_i^2b_kb_l + b_i^2a_ka_l)Y_0Y_1 + a_ib_kb_lb_iY_1^2)\\
		&= A_i\cdot (-(a_ib_k - a_kb_i)(a_ib_l - a_lb_i)Y_0Y_1 + a_ia_ka_lb_iY_0^2\\
		&\quad  + (a_i^2b_kb_l + b_i^2a_ka_l)Y_0Y_1 + a_ib_kb_lb_iY_1^2)
		\end{align*}
		and we have that
		\begin{align*}
		(a_ia_ka_lb_iY_0^2& + (a_i^2b_kb_l + b_i^2a_ka_l)Y_0Y_1 + a_ib_kb_lb_iY_1^2)X_i^m \\
		&= (a_ia_ka_lb_iY_0^2 - a_ib_ib_kb_lY_1^2 - a_ia_ka_lb_iY_0^2 + a_ib_kb_lb_iY_1^2)X_i^m\\
		&= 0
		\end{align*}
		in $\tilde{J}^{\rho + (r+1,n+1)}$. From this, we see that  	
		\begin{align*} 
		\tilde{C} &= m^{2n+1}Y_0Y_1X_0\cdots X_n\Big(\sum_{i\neq k,l} (a_ib_k-a_kb_i)(a_ib_l-a_lb_i) A_i \\
		&\quad + (a_kb_{j'} - a_{j'}b_k)(a_kb_l - a_lb_k)\tilde{A_k}+ (a_lb_{j'} - a_kb_{j'})(a_lb_k - a_kb_l)\tilde{A_l}\Big)
		\end{align*} 	
		for some $j'\notin \{j,k,l\}$, where 
		$$\tilde{A_k} = X_k^mX_0^{m-2}\cdots X_n^{m-2}\prod_{i\neq k,l,j'}L_i\text{ and }\tilde{A_l} = X_l^mX_0^{m-2}\cdots X_n^{m-2}\prod_{i\neq k,l,j'}L_i.$$ 
		We note that 
		\begin{align*} 
		(a_lb_i-a_ib_l)&(a_kb_i-a_ib_k)A_i \\
		&= (a_lb_i-a_ib_l)(a_kb_i-a_ib_k)\frac{(a_lb_j-a_jb_l)(a_kb_j - a_jb_k)}{(a_lb_i-a_ib_l)(a_kb_i-a_ib_k)}A_j\\
		&= (a_jb_l - a_lb_j)(a_jb_k-a_kb_j)A_j
		\end{align*} for $i\neq j$ and that 
		\begin{align*}
		(a_kb_{j'} - a_{j'}b_k)&(a_kb_l - a_lb_k)\tilde{A_k} \\
		&= (a_kb_{j'} - a_{j'}b_k)(a_kb_l - a_lb_k)\frac{(a_jb_{j'} - a_{j'}b_j)(a_jb_l - a_lb_j)}{(a_kb_{j'} - a_{j'}b_k)(a_kb_l-a_lb_k)}\tilde{A_j}\\
		&= \frac{(a_jb_k - a_kb_j)(a_jb_{j'} - a_{j'}b_j)(a_jb_l - a_lb_j)}{a_jb_{j'} - a_{j'}b_j}A_j \\
		&= (a_jb_k - a_kb_j)(a_jb_l - a_lb_j)A_j 
		\end{align*}
		and similarly, we have that 
		$$(a_lb_{j'} - a_kb_{j'})(a_lb_k - a_kb_l)\tilde{A_l} = (a_jb_k - a_kb_j)(a_jb_l - a_lb_j)A_j.$$
		Putting this together, we see that 
		$$C = m^{2n+1}(n+1)(a_jb_k - a_kb_j)(a_jb_l - a_lb_j)A_j.$$
		By Theorem \ref{theorem trace}, we have that $$\text{Tr}(\omega_{A}\cup\omega_B) = m^{3n+2}(n+1)^2(a_jb_k - a_kb_j)(a_jb_l-a_lb_j)\lambda$$
		as desired. 
	\end{proof}
	
	\subsection{The quadratic Euler characteristic} 
	\begin{nota} 
		Assume that $n = 2p$ is even.
	\end{nota} 
	In order to prove Proposition \ref{qec of L}, we will need to compute the form $Q$ from Theorem \ref{Motivic Gauss Bonnet}, i.e. \cite[Corollary 8.7]{LevineGB}, given by 
	$$H^p(\mathcal{X}, \Omega^p_{\mathcal{X}})\times H^p(\mathcal{X}, \Omega^p_{\mathcal{X}})\xrightarrow{\cup} H^n(\mathcal{X}, \Omega^n_{\mathcal{X}})\xrightarrow{\text{Tr}} k.$$ 
	The result from the previous section will allow us to do so on primitive cohomology, but we will also need to understand the form $Q$ on the complement. 
	\begin{constr} 
		Using Proposition \ref{Botts theorem for general products}, we have that $H^p(\mathbb{P}^1\times\mathbb{P}^n,\Omega_{\mathbb{P}^1\times\mathbb{P}^n}^{p})$ has rank two over $k$. Generators are given by 
		$$\alpha = c_1(\mathcal{O}(1,0))\cup c_1(\mathcal{O}(0,1))^{p-1}\text{ and }\beta = c_1(\mathcal{O}(0,1))^{p}.$$  
		Also, one can show that $H^{p+1}(\mathbb{P}^1\times\mathbb{P}^n,\Omega_{\mathbb{P}^1\times\mathbb{P}^n}^{p+1})$ has rank two. Generators are given by 
		$$\alpha' = c_1(\mathcal{O}(1,0))\cup c_1(\mathcal{O}(0,1))^{p}\text{ and }\beta' = c_1(\mathcal{O}(0,1))^{p+1}.$$  
	\end{constr}
	\begin{lem}
		The complement to $k\cdot i^*\alpha\oplus k\cdot i^*\beta$ inside $H^p(\mathcal{X}, \Omega^p_{\mathcal{X}})$ under the trace pairing is precisely $H^p(\mathcal{X}, \Omega^p_{\mathcal{X}})_{prim}$. 
	\end{lem}
	\begin{proof} 
		Let $\gamma \in H^p(\mathcal{X}, \Omega^p_{\mathcal{X}})$ be an arbitrary element. Using the projection formula, we note that 
		\begin{align*}
		\text{Tr}(i^*\alpha\cup\gamma) &= \text{Tr}_{\mathbb{P}^1\times\mathbb{P}^n}(i_*(i^*\alpha\cup\gamma))\\
		&= \text{Tr}_{\mathbb{P}^1\times\mathbb{P}^n}(\alpha\cup i_*\gamma). 
		\end{align*}
		If $\gamma\in \ker(i_*)$, this implies that $\text{Tr}(i^*\alpha\cup\gamma) = 0$, and a similar argument shows that $\text{Tr}(i^*\beta\cup\gamma)=0$, and so $\ker(i_*) = H^p(\mathcal{X}, \Omega^p_{\mathcal{X}})_{prim}$ is contained in the complement of $k\cdot i^*\alpha\oplus k\cdot i^*\beta$. \\
		We now show that the other inclusion holds. As $\alpha'$ and $\beta'$ are generators of $H^{p+1}(\mathbb{P}^1\times\mathbb{P}^n,\Omega_{\mathbb{P}^1\times\mathbb{P}^n}^{p+1})$, we have that $i_*\gamma = a\alpha'+ b\beta'$ for certain $a,b\in k$. Note that $\alpha\cup\alpha' = 0 = \beta\cup \beta'$, while $\alpha \cup \beta'$ and $\beta\cup\alpha'$ give the generator of $H^{n+1}(\mathbb{P}^1\times\mathbb{P}^n, \Omega_{\mathbb{P}^1\times\mathbb{P}^n}^{n+r})$. It follows that 
		$$\text{Tr}_{\mathbb{P}^1\times\mathbb{P}^n}(\alpha\cup i_*\gamma) =  b\text{Tr}_{\mathbb{P}^1\times\mathbb{P}^n}(\alpha\cup \beta') = b$$ and similarly 
		$$\text{Tr}_{\mathbb{P}^1\times\mathbb{P}^n}(\beta\cup i_*\gamma) =  a\text{Tr}_{\mathbb{P}^1\times\mathbb{P}^n}(\beta\cup \alpha') = a.$$ 
		If these are both zero then $a = b = 0$, i.e. we have $\gamma\in\ker(i_*) = H^p(\mathcal{X}, \Omega_{\mathcal{X}}^p)_{prim}$. 
	\end{proof} 
	\begin{prop}\label{contribution primitive cohomology}
		We have that $\text{Tr}(i^*\alpha\cup i^*\alpha) = 0$ and $\text{Tr}(i^*\beta\cup i^*\beta) = 1$. 
	\end{prop}
	\begin{proof}
		We have that 
		\begin{align*} 
		\text{Tr}(i^*\alpha\cup i^*\alpha) &= \text{Tr}_{\mathbb{P}^1\times\mathbb{P}^n}(\alpha\cup i_*i^*\alpha) \\
		&= \text{Tr}_{\mathbb{P}^1\times\mathbb{P}^n}(\alpha^2\cup c_1(\mathcal{O}(1,m))) \\
		&= 0
		\end{align*} as $\alpha^2 = 0$ and $i_*i^*\alpha = \alpha\cup c_1(\mathcal{O}(1,m))$. Similarly 
		\begin{align*} 
		\text{Tr}(i^*\beta\cup i^*\beta) &= \text{Tr}_{\mathbb{P}^1\times\mathbb{P}^n}(\beta^2\cup c_1(\mathcal{O}(1,m)))\\
		&= \text{Tr}_{\mathbb{P}^1\times\mathbb{P}^n}(c_1(\mathcal{O}(0,1))^n\cup c_1(\mathcal{O}(1,m))) \\
		&= 1
		\end{align*}  
		as desired. 
	\end{proof} 
	\begin{proof}[Proof of Proposition \ref{qec of L}]
		Using Theorem \ref{Motivic Gauss Bonnet}, we know that $\chi(\mathcal{X}/k)$ is hyperbolic for $n$ odd. For $n$ even, the quadratic Euler characteristic is equal to a hyperbolic form plus the trace form. Therefore, assume from now on that $n$ is even. \\
		In order to compute the trace form, we evaluate it on basis elements of $J$. Choose a generator $A_j = X_j^mX_0^{m-2}\cdots X_n^{m-2}\prod_{i\neq j,k,l}L_i$ of $J^\rho$ as in Proposition~\ref{generators of the Jacobian ring}. Note that if $AB = \lambda A_j$ for some $\lambda\in k^*$ and two distinct basis elements $A,B$, then $BA = \lambda A_j$ and one can check that this yields a hyperbolic form. If $m$ is odd, there are no basis elements that square to a nonzero multiple of $A_j$. If $m=2q$ is even, then $\rho = (n-2, 2q(n+2) - 2(n+1))$ is divisible by $2$. For each subset $\{i_0,\cdots, i_{\frac{n-2}{2}}\}\subset \{0,\cdots, n\}\setminus \{j,k,l\}$, we find the element
		$$A_{i_0,\cdots, i_{\frac{n-2}{2}}} = X_j^{q}X_0^{q-1}\cdots X_n^{q-1}\prod_{i\in \{i_0,\cdots, i_{\frac{n-2}{2}}\}}L_i$$
		of $J^{\frac{\rho}{2}}$ such that, using Lemma \ref{Lj in terms of Li Lk} again, we have that
		\begin{align*} 
		A_{i_0,\cdots, i_{\frac{n-2}{2}}}^2 &= X_j^{m}X_0^{m-2}\cdots X_n^{m-2}\prod_{i\in \{i_0, \cdots, i_{\frac{n-2}{2}}\}}L_i^2 \\
		&= \left(\frac{\prod_{i\in \{i_0,\cdots, i_{\frac{n-2}{2}}\}}(a_ib_j-a_jb_i)}{\prod_{i\notin \{j,k,l,i_0, \cdots, i_{\frac{n-2}{2}}\}}(a_ib_j-a_jb_i)}\right)A_j
		\end{align*} 
		We note that all such $A_{i_0,\cdots, i_{\frac{n-2}{2}}}$ are multiples of each other in $J^{\frac{\rho}{2}}$, because $X_j$ has degree $2q-1 = m-1$, so we can use the same argument as for Proposition~\ref{generators of the Jacobian ring}. Also, if $j, j'\in \{0,\cdots, n\}$ are distinct and $A,A'\in J^{\frac{\rho}{2}}$ are such that $A^2 = \lambda A_j$ and $(A')^2 = \lambda'A_{j'}$ for some $\lambda,\lambda'\in k^*$, then $A$ and $A'$ are distinct elements of $J^{\frac{\rho}{2}}$. Therefore, for each $j$ there is exactly one basis element that squares to a nonzero multiple of $A_j$, and we choose this basis element to be such that $j, k$ and $l$ lie in the complement of the $i_j$. \\
		Using Lemma \ref{trace of Aj}, this gives rise to the term $\sum_{j=0}^n\langle \prod_{i\neq j} (a_jb_i - a_ib_j) \rangle$. Also, we note that by Proposition \ref{contribution primitive cohomology}, the contribution coming from primitive cohomology is the form with matrix 
		$$\left( \begin{matrix}
		0 & m \\ m & 1 
		\end{matrix} \right)$$
		which is hyperbolic. Finally, as the rank of $\chi(\mathcal{X}/k)$ is equal to $\deg(c_n(T_{\mathcal{X}}))$ by \cite[Theorem 5.3]{LevineGB}, we see that the coefficient of $H$ is equal to $A_{n,m}$ as desired (also in the case where $n$ is odd). 
	\end{proof}
	We can now also deduce Corollary \ref{qec of X}. 
	\begin{proof}[Proof of Corollary \ref{qec of X}]
		Using Proposition \ref{cut and paste proposition}, we have that  
		$$\langle -1 \rangle \chi(X/k) = \chi(\mathcal{X}/k) - \chi(\mathbb{P}^n/k)$$ and $\chi(\mathbb{P}^n/k) = \sum_{i=0}^n\langle -1 \rangle^i$ as we saw in Example \ref{Euler char of projective space}, i.e. \cite[Proposition 1.4(4)]{LevineAEGQF}. So for odd $n$, we have that $\chi(\mathbb{P}^n/k)$ is hyperbolic, and for even $n$ we get an extra $\langle 1 \rangle$-term. This gives the desired statement. 
	\end{proof} 
	
	\subsection{Checking the answer using the quadratic Riemann-Hurwitz formula}
	There is another way to compute $\chi(\mathcal{X}/k)$: using the quadratic Riemann-Hurwitz formula from \cite{LevineAEGQF}. We will now do this and see that we recover Proposition \ref{qec of L}. 
	\begin{nota}
		Because we know that $\chi(\mathcal{X}/k)$ is hyperbolic if $n$ is odd by Theorem \ref{Motivic Gauss Bonnet}, we assume throughout that $n$ is even. 
	\end{nota}
	Note that the natural projection map $\mathbb{P}^1\times\mathbb{P}^n\to\mathbb{P}^1$ yields a projective morphism $f: \mathcal{X}\to \mathbb{P}^1$. The fiber of $f$ over a point $y\in\mathbb{P}^1$ is isomorphic to the zero locus of $\sum_{i=0}^nL_i(y)X_i^m$. This is smooth if $L_i(y)\neq 0$ for all $i\in\{0,\cdots, n\}$. If there is a $j\in\{0,\cdots, n\}$ such that $L_j(y) = 0$, we note that by Lemma \ref{Lj in terms of Li Lk}, we must have that $L_i(y)\neq 0$ for all $i\neq j$. The fiber of $f$ over $y$ is now the cone over the zero locus of $\sum_{i\neq j}L_i(y)X_i^m$ with vertex $[e_j]$ given by the unit vector that has a one on the $j$'th spot and zero's everywhere else. The vertex $[e_j]$ is the only singular point. In particular, this fiber is smooth inside the $\mathbb{P}^{n-1}$ defined by setting $X_j = 0$. \\
	To apply the quadratic Riemann-Hurwitz formula, the first thing which we need to do is to identify the set $\mathfrak{c}(f)$ of critical points of $f$, i.e. the locus of points of $\mathcal{X}$ where $df = 0$.
	\begin{prop}
		The critical locus $\mathfrak{c}(f)$ of $f$ has $n+1$ elements, and consists of those points satisfying $L_j = 0$ for some $j\in\{0,\cdots, n\}$ and $X_i = 0$ for all~$i\neq j$. 
	\end{prop}
	\begin{proof} 
		Let $j\in\{0,\cdots n\}$ and consider the affine patch of $\mathbb{P}^1\times\mathbb{P}^n$ given by $X_j\neq 0$ and $Y_0\neq 0$, with coordinates $y = \frac{Y_1}{Y_0}$ and $x_i = \frac{X_i}{X_j}$ for $i\neq j$. Here, $\mathcal{X}$ is given by the equation 
		$$\sum_{i\neq j}a_ix_i^m + y\sum_{i\neq j}b_ix_i^m + a_j + yb_j = 0.$$
		This implies in particular that 
		$$m\sum_{i\neq j}(a_i + yb_i)x_i^{m-1}dx_i + \left(\sum_{i\neq j}b_ix_i^m + b_j\right)dy =0.$$
		We have that $f$ is given by $f(y,x_0,\cdots, x_{j-1}, x_{j+1},\cdots x_n) = y$ and so 
		\begin{align*}
		df &= dy = -\frac{m}{\sum_{i\neq j}b_ix_i^m + b_j}\sum_{i\neq j}(a_i + yb_i)x_i^{m-1}dx_i.
		\end{align*}
		This implies that $\sum_{i\neq j}L_i(y)x_i^{m-1}dx_i = 0$ for those points. This gives us two possibilities for a critical point. \\
		First, we can have that $L_i(y)\neq 0$ for $i\neq j$, so that we must have that $x_i = 0$ for all $i\neq j$. We need in addition that $a_j + yb_j = L_j(y) = 0$, as the critical point also has to lie on $\mathcal{X}$. \\
		Secondly, we can have that there is some $k\neq j$ such that $L_k(y) = 0$. In this case, we have that $x_i = 0$ for all $i\neq k$. But for such a point to lie on $\mathcal{X}$, we need the condition that $L_j(y) = 0$ again as well, which yields a contradiction, as $a_kb_j - a_jb_k \neq 0$. \\
		Repeating this construction for other choices of $j$, we deduce the desired statement. 
	\end{proof} 
	\begin{rem}
		All critical points may not lie in the same affine patch, but all critical values (so all $y\in\mathbb{P}^1$ such that $y = f(p)$ for $p$ a critical point) do lie in the same affine patch of $\mathbb{P}^1$. Namely, if there would be a critical value with $Y_0 = 0$ then $L_j = 0$ would imply that $Y_1 = 0$ (as all $a_i,b_i\in k^*$). 
	\end{rem} 
	\begin{nota} 
		Let $y\subset \mathcal{X}$ be the subscheme of critical points of $f$. Consider the closed point $y'$ of $y$ given by $L_j = 0$. Consider the affine patch of $\mathbb{P}^1\times\mathbb{P}^n$ given by $Y_0\neq 0$ and $X_j\neq 0$ as before. We have that $\mathcal{O}_{\mathcal{X},y'}$ is a regular local ring and we choose parameters $x_0,\cdots, x_{j-1}, x_{j+1},\cdots, x_n$ which generate the maximal ideal $\mathfrak{m}_{y'}$ (we will not need $L_j$ as an additional generator because we already work on $\mathcal{X}$). Furthermore, let $x$ be the subscheme of $\mathbb{P}^1$ defined by $L_j = a_j + b_jy = 0$. By \cite[Remark 10.9]{LevineAEGQF} we have that 
		$$t_x = \frac{a_j + yb_j}{b_j}$$ is a normalized parameter. Define $s_i = -\frac{m}{\sum_{k\neq j}b_kx_k^m + b_j}(a_i + yb_i)x_i^{m-1}$ for $i\neq j$ and let $$[B_{y'}] := [B_{s_*,x_*}]\in \text{GW}(k(y'))$$ be the corresponding Scheja-Storch form (see \cite[Theorem 4.1 (3)]{LevineAEGQF}). 
	\end{nota} 
	By the quadratic Riemann-Hurwitz formula, see \cite[Corollary 10.6]{LevineAEGQF}, we have that 
	$$\chi(\mathcal{X}/k) = \sum_{y'\in \mathfrak{c}(f)}\text{Tr}_{k(y')/k}([B_{y'}]) - D(f)\cdot H$$
	where $D(f)\in\mathbb{Z}$ (see \cite[Theorem 10.2]{LevineAEGQF}). Also, note that for all $y'\in \mathfrak{c}(f)$, we have that $k(y') = k$ and so the trace doesn't have any effect. It therefore remains to compute $[B_{y'}]$ for all $y'\in \mathfrak{c}(f)$. 
	\begin{prop} 
		We have that $$[B_{y'}] = \begin{cases} 
		\frac{1}{2}(m-1)^n\cdot H  & \text{ if $m$ is odd} \\
		(\frac{1}{2}((m-1)^n - 1)\cdot H + \langle \prod_{i\neq j}(a_ib_j - a_jb_i)\rangle) &\text{ if $m$ is even} 
		\end{cases}$$
	\end{prop} 
	\begin{proof} 
		Let $y'$ be a critical point again. We note that $x_0,\cdots, x_n$ is a local framing for $(\det(\mathfrak{m}_{y'}/\mathfrak{m}_{y'}^2)^\vee)^{\otimes 2}$. The section we have is not diagonalizable, but we note that if we set $$\lambda = -\frac{m}{\sum_{k\neq j}b_kx_k^m + b_j}$$ and do a change of coordinates where we switch $dx_i$ with $\lambda dx_i$ for all $i$, it is, and it will only change the determinant by $\lambda^n$ which will be a square as $n$ is assumed to be even. We can therefore apply \cite[Example 4.5]{LevineAEGQF} combined with \cite[Corollary 4.3]{LevineAEGQF} to see that  
		\begin{align*}
		[B_{y'}] &= \langle (-\frac{m}{\sum_{i\neq j}b_ix_i^m + b_j})^n \prod_{i\neq j}(a_i + yb_i)\rangle \sum_{i=0}^{(m-1)^n - 1}\langle -1 \rangle^i.
		\end{align*} 
		We note that for $m$ odd, $(m-1)^n$ is even and so $$\sum_{i=0}^{(m-1)^n - 1}\langle -1 \rangle^i = \frac{1}{2}(m-1)^n\cdot H.$$ For even $m$, we have that $$\sum_{i=0}^{(m-1)^n - 1}\langle -1 \rangle^i  = \frac{1}{2}((m-1)^n - 1)\cdot H + \langle 1 \rangle.$$ Furthermore, the term $ (-\frac{m}{\sum_{i\neq j}b_ix_i^m + b_j})^n$ is a square as we assumed that $n$ is even. \\
		Finally, as $L_j = 0$ we have that $y = -\frac{a_j}{b_j}$ and so 
		$$\langle \prod_{i\neq j}(a_i + yb_i)\rangle = \langle \frac{1}{b_j^n} \prod_{i\neq j}(a_ib_j - a_jb_i)\rangle = \langle \prod_{i\neq j}(a_ib_j - a_jb_i)\rangle.$$
		So 
		$$[B_{y'}] = \begin{cases} 
		\frac{1}{2}(m-1)^n\cdot H  \langle \prod_{i\neq j}(a_ib_j - a_jb_i)\rangle &\text{ if $m$ is odd} \\
		(\frac{1}{2}(((m-1)^n - 1)\cdot H + \langle 1 \rangle) \langle \prod_{i\neq j}(a_ib_j - a_jb_i)\rangle  &\text{ if $m$ is even} 
		\end{cases}$$ which proves the statement. 	
	\end{proof} 
	Applying the quadratic Riemann-Hurwitz formula, we see from this that 
	$$\chi(\mathcal{X}/k)  = \begin{cases} 
	A_{n,m}\cdot H &\text{ if $m$ is odd}\\
	A_{n,m}\cdot H + \sum_{j=0}^n \langle \prod_{i\neq j}(a_ib_j - a_jb_i)\rangle &\text{ if $m$ is even}
	\end{cases}$$
	which coincides with the result of Proposition \ref{qec of L}. This therefore gives the same quadratic Euler characteristic of $X$ as we had before. 
	
	\phantomsection 
	\addcontentsline{toc}{section}{Bibliography}
	\bibliographystyle{plain}
	\bibliography{NewThesisSources}
	
\end{document}